\documentclass[10pt]{amsart}      % for a regular run

\usepackage[utf8]{inputenc}
\usepackage{amsmath}
\usepackage{appendix}
\usepackage{xcolor}
\usepackage{amsfonts}

\usepackage{todonotes}
\usepackage{amssymb}
\usepackage[normalem]{ulem}
\usepackage{subfigure}    
\usepackage{wrapfig}
\usepackage{mathtools}
\DeclareUnicodeCharacter{221E}{ }

\usepackage{booktabs}
\newtheorem{theo}{Theorem}[section] % use one numbering scheme for all these objects 

\newtheorem{prop}[theo]{Proposition} \newtheorem{lemm}[theo]{Lemma}
\newtheorem{rema}[theo]{Remark} 
\newtheorem{coro}[theo]{Corollary} 
 \newtheorem{assum}[theo]{Assumption}

\newenvironment{tightitemize}{%
    \list{{\textup{$\bullet$}}}{\settowidth\labelwidth{{\textup{\qquad}}}
    \leftmargin\labelwidth \advance\leftmargin\labelsep
    \parsep 0pt plus 0pt minus 0pt \topsep 0pt \itemsep 0pt
    }}{\endlist}
    
\newcommand{\Diag}{{\mbox{Diag}}}

\newcommand{\Cml}{\mathcal{C}_\mu^L(\mathbb{R}^d)}
\newcommand{\beq}{\begin{eqnarray}}
\newcommand{\eeq}{\end{eqnarray}}
\newcommand{\beqs}{\begin{eqnarray*}}
\newcommand{\eeqs}{\end{eqnarray*}}

\newcommand{\mg}[1]{{\color{black} #1}}

\newcommand{\txi}{\tilde{\xi}}

\newcommand{\rhogd}{\rho_{\tiny\mbox{GD}}}
\usepackage{hyperref}
\usepackage{cleveref}
\graphicspath{{./figures/}}
 
\usepackage{caption}
\usepackage{mwe}
\usepackage{amsaddr}

\makeatletter
\renewcommand{\email}[2][]{%
  \ifx\emails\@empty\relax\else{\g@addto@macro\emails{,\space}}\fi%
  \@ifnotempty{#1}{\g@addto@macro\emails{\textrm{(#1)}\space}}%
  \g@addto@macro\emails{#2}%
}
\makeatother

\calclayout
\definecolor{winered}{rgb}{0.5,0,0}
\definecolor{darkgreen}{rgb}{0.2,0.8,0.1}
\hypersetup
{
    pdfauthor={Dr. Mert G\"urb\"uzbalaban},
    pdfsubject={Robustly Stable Accelerated Momentum Methods},
    pdftitle={Robustly Stable Accelerated Momentum Methods},
    pdfkeywords={LaTeX, PDF, hyperlinks}
    pdfborder={0 0 0},
colorlinks=true,
 linkcolor=blue,
urlcolor={winered},
filecolor={orange},
citecolor={orange}%,
}

\begin{document}
%%%%%%%%%%%%%%%%

\address{Dedicated to Professor Michael Overton on his seventieth birthday}
\vspace{-0.6in}
\author{Mert G\"urb\"uzbalaban} % Department of Management Science and 
\address{Department of Management Science and Information Systems\\ 
\& Department of Electrical and Computer Engineering\\
\& Department of Statistics, Rutgers University, NJ, USA}
\email{mg1366@rutgers.edu}
\date{}

%Information Systems\\ Rutgers University
% Full title. Sample:
% \TITLE{Bundling Information Goods of Decreasing Value}
% Enter the full title:
\title[Robustly Stable Accelerated Momentum Methods]{Robustly Stable Accelerated Momentum Methods with a Near-Optimal %Pareto-optimal 
$L_2$ Gain and {$H_\infty$} Performance}

\maketitle

\begin{abstract}
 % Enter your abstract
We consider the problem of minimizing a strongly convex smooth function %in the inexact gradient setting 
where the gradients are subject to additive worst-case deterministic errors that are square-summable. We  study the trade-offs between the convergence rate and robustness to gradient errors when designing the parameters of a first-order algorithm. We focus on a general class of momentum methods (GMM) with constant stepsize and two momentum parameters which can recover gradient descent (GD), Nesterov's accelerated gradient (NAG), the heavy-ball (HB) and the triple momentum methods (TMM) as special cases. We measure the robustness of an algorithm in terms of the cumulative suboptimality over the iterations normalized by the squared $\ell_2$ norm of the gradient errors. This quantity can be interpreted as the (squared) $\ell_2$ gain of a dynamical system that represents the GMM iterations where the input is the gradient error sequence and the output is a weighted distance to the optimum. For quadratic objectives, we compute the $\ell_2$ gain explicitly leveraging its representation as the $H_\infty$ norm of the GMM system in the frequency domain %dynamical system corresponding to GMM %which is a robust stability measure 
and construct gradient errors that lead to worst-case performance explicitly. %by a closed-form formula. 
We also study the stability of GMM with respect to multiplicative errors by characterizing the structured real and stability radius of the GMM system through their connections to the $H_\infty$ norm. This allows us to compare GD, HB, NAG methods in terms of robustness, and argue that HB is not as robust as NAG despite being the fastest in terms of the rate. %We then characterize the real and complex stability radius of the GMM system, which allows us to quantify the stability of GMM with respect to multiplicative noise. 
We then develop the robustly stable heavy ball method that can be faster than NAG while being at the best robustness level possible. We also propose the robustly stable gradient descent that is the fastest version of GD with constant stepsize while being at the best robustness level. Finally, we extend our framework to general strongly convex smooth objectives, providing non-asymptotic rate results for inexact GMM methods and bounds on the %induced 
$\ell_2$ gain where we can choose the GMM parameters to systematically trade off the rate to robustness in a computationally tractable framework.\looseness=-1
\end{abstract}

\section{Introduction.}\label{sec-intro}
First-order methods are workhorse methods in optimization with a long history. %with a long that rely on gradient information are workhorse methods in optimization. % and have favorable scalability properties to large
Their favorable scalability properties to large dimensions and their potential to generate low to medium accuracy solutions at a low computational cost make them the preferred method for many applications. %Among first-order methods, gradient descent and its accelerated versions such as Nesterov's accelerated gradient (NAG), heavy ball or triple momentum methods admit 
The classical convergence theory of fundamental first-order methods such as gradient descent (GD) and its momentum-based accelerated versions \cite{nesterov2018lectures,polyakintroduction} assume that the gradients can be computed exactly, and provides convergence rates to the solution for the iterates. In this setting, accelerated methods such as Nesterov's accelerated gradient (NAG), Polyak's heavy-ball (HB) method or triple momentum method (TMM) \cite{gannot2022frequency} are superior to GD in the sense that they admit optimal rates in various convex and strongly convex settings improving upon the convergence rate of GD (see e.g. \cite{nesterov2018lectures,polyakintroduction,scoy-tmm-ieee}).

On the other hand, in many applications, the gradients are inexact, containing errors that can be of additive or multiplicative nature, also known as \emph{absolute errors} and \emph{relative errors} \cite{polyakintroduction}. These gradient errors can be random as in the case of stochastic gradient-like methods, or they can be deterministic as in the case of incremental gradient methods or in problems when evaluating the gradients require solving a subproblem whose solutions are computed inexactly \cite{zhang2022sapd+,devolder2013exactness,daspremont}. There are also many other settings and applications when the (inexactness) errors in the gradients admit a deterministic worst-case nature, see e.g. \cite{luo1993error,bertsekas2000,devolder2013exactness}. Such errors get accumulated and may potentially be amplified over time, %he iterations may deviate from their regular trajectories when subject to erros and the convergence rate may deteriorate %
causing the algorithm's trajectory to deviate from its regular (errorless) trajectory, potentially resulting in divergence or slower convergence. The total amount of deviation depends on the robustness to gradient errors of the underlying optimization algorithm. %For example, accelerated gradient methods are based on an extrapolation step
%based on momentum averaging where the steps are extrapolated, this results in potentially larger directional moves where the gradient errors can be amplified. 
This makes robustness to gradient errors another key performance metric in addition to (convergence) rates when designing a first-order method. In particular, accelerated methods such as NAG and HB with the standard choice of parameters or parameters that optimize the rate were found to be considerably less robust to noise in a number of inexact gradient settings (see e.g.  \cite{devolder2013exactness,aybat2018robust,aybat2019universally,mohammadi2020robustness,van2021speed}), resulting in a performance that can be worse than that of GD. This raises the natural question of whether one can identify alternative parameters for accelerated methods such as NAG or HB that can lead to a more robust behavior to inexactness while retaining their fast convergence rates as much as possible.\looseness=-1%and whether achieving robustness can be achieved without deteriorating the accelerated convergence rates much. 
%identify alternative parameters where momentum methods can be fast and robust in some gradient error settings. %and if so what are some fundamental limitations 
%This motivates the development of methods that can perform well \emph{universally}, i.e. in both exact and inexact gradient settings where convergence rates and robustness are two performance metrics of interest. %For stochasti%There are currently a%While there are ways of quan
%While for stochastic noise, robustness can be measured in terms of the average suboptimality 
%That being said, there is need to rigourously define what is meant by robustness in the literature. %That being said, % where there is further need to introduce notions of robustness to . 

In this work, we consider the trade-offs between the robustness and the convergence rate when designing a first-order method. Such trade-offs have been studied previously for momentum methods when the gradient errors are random where universally optimal algorithms are known \cite{aybat2018robust,aybat2019universally,ac-sa,ghadimi-lan-part-one,fallah2022robust}, but the deterministic error setting is relatively less studied where existing notions of robustness such as suboptimality at the last iterate or at an averaged iterate do not precisely capture the cumulative deviations in the whole trajectory due to accumulation of errors. %Existing notions of robustness are related to suboptimality at the last iterate or at an averaged iterate \cn, but do not consider the whole trajectory of the iterates. 
For unbiased stochastic errors, robustness can be defined through the asymptotic expected suboptimality of the iterates (normalized by the variance of the noise), which is equivalent to the squared \emph{$H_2$ norm} of a transformed dynamical system corresponding to the iterations when the objective is a quadratic and errors are i.i.d. Gaussian \cite{aybat2018robust}. The $H_2$ norm is a well-known robustness measure used to design control systems that are robust to stochastic perturbations \cite{arzelier2010h2,zhou1996robust}. On the other hand, for \emph{deterministic} worst-case errors, the natural analog notion of robustness that is well-studied in robust control theory to design robust control systems is the \emph{$\ell_2$ gain} \cite{zhou1996robust,zhou1998essentials,van2016l2}. The $\ell_2$ gain is a measure of sensitivity to external errors, quantifying how much the $\ell_2$ norm of the error input sequence to a system is amplified (from input to the output) in the worst-case. However, to our knowledge, $\ell_2$ gains of first-order optimization algorithms (such as GD, HB, NAG, TMM) have never been characterized for quantifying and comparing their robustness to worst-case gradient errors and for designing the parameters of a first-order algorithm to achieve systematic trade-offs between the rate and robustness,
%for the design of first-order methods that can achieve a Pareto-optimal performance %(with respect to robustness and convergence rate as competing performance metrics) 
which will be the main focus of this work. 
%In this work, our main contribution is to calculate and estimate the induced $\ell_2$ norm for a general class of momentum methods and then choose the parameters (stepsize and momentum) to systematically trade-off the convergence rate and robustness in a computationally tractable framework. 
%In this work, our contribution is to comp induced $\ell_2$ norms on a general class of momentum methods for strongly convex minimization 
 %In this work, our contribution will be to estimate the induced $\ell_2$ gains 
%as a function of their parameters such as the stepsize and the momentum. 
 %where we look at the ratio of the $\ell_2$ norm of the output sequence of an algorithm to the $\ell_2$ norm of the input error sequence $\{w_k\}$. %. Induced $\ell_2$ gains is a measure of how much the input noise is amplified %which is the ratio of the $\ell_2$ norm of the output $z_k$ to the $\ell_2$ norm of the errors up to a constant term
%which is a measure of how much the $\ell_2$ norm of the output of an optimization algorithm $z_k$ 
%which will be the focus of our paper. %Main novelty of our work is to propose cumulative suboptimality over iterations as a robustness metric which measures the deviation of the trajectory from the optimal value. %We then estimate this metric  %then,  and then 
%and place it in a computationally tractable framework to trade robustness with convergence rates. % (the convergence rate. %The novelty of our approach in this work is to propose a new robustness metric  
For this purpose, we consider a general class of momentum methods (GMM) for minimizing a strongly convex smooth\footnote{Here, smoothness means the gradient of $f$ is globally Lipschitz, see e.g. \cite{zhang2022sapd+,aybat2018robust,aybat2019universally}.} objective $f$ that has three constant parameters (the stepsize $\alpha$, and two momentum parameters $\beta $ and $\nu$). GMM is a rich class in the sense that it recovers GD and the momentum methods such as HB, NAG and TMM as special cases depending on the choice of parameters \cite{lessard2016analysis,hu2017dissipativity,can2022entropic}. %GMM generates iterates $\{x_k\}$ 
In the inexact setting with additive deterministic errors, GMM generates iterates $\{x_k\}_{k\geq 0}$ starting from the initialization $x_0 = x_{-1}\in \mathbb{R}^d$, based on momentum averaging and inexact gradients where $w_k$ is the additive deterministic gradient error at step $k$ (see Sec. \ref{sec-prelim} for details of GMM updates). %We can view
We consider gradient error sequences $w_k \in \ell_2$, i.e. the errors are square-summable satisfying $\|w\|_{\ell_2}:=\sum_{k\geq 0}\|w_k\|^2 < \infty$. Given fixed parameters $(\alpha,\beta,\nu)$, if $\{z_k\}$ is the output sequence of GMM that is of interest, then we say that the $\ell_2$ gain (from input $\{w_k\}$ to output $\{z_k\}$) is finite if there exists $\gamma>0$ such that for all $w_k \in \ell_2$, it holds %that
\begin{equation} %\vspace{-0.4in}
\setlength{\abovedisplayskip}{5pt}
\sum_{k\geq 0}\|z_k\|^2 \leq \gamma^2  \sum_{k\geq 0 } \|w_k\|^2 + \overline{H}(x_0),
\setlength{\belowdisplayskip}{10pt}
\label{eq-cum-dist-squared}
%\vspace{-0.13in}
\end{equation}where $\overline{H}:\mathbb{R}^{d}\to \mathbb{R}$ is a function that only depends on the algorithm parameters $(\alpha,\beta,\nu)$ and the problem parameters (strong convexity constant $\mu$ and the Lipschitz constant $L$ of the gradient $\nabla f$) while being independent of the error sequence $\{w_k\}$. The $\ell_2$ gain is then defined as the infimum of such $\gamma$, which we will denote by $L_{2,*}$ in this paper. %satisfying %H(x_*)=0$. 
As an example, if we take $z_k = x_k - x_*$ where $x_*$ is the minimum of $f$, then the $\ell_2$ gain would allow us to measure how much the error effects the cumulative distance squared over the iterations up to a constant term coming from the initialization. However, for optimization purposes, suboptimality is often more relevant. Therefore, in this paper we take the approach of setting $z_k$ as a weighted distance to the optimum so that $\|z_k\|^2 = f(x_k) - f(x_*)$ (see Sec. \ref{subsec-l2-gain-def} for details) and we will obtain bounds of the form \eqref{eq-cum-dist-squared} that control the cumulative deviation from the optimal value along the iterations. In this sense, $\ell_2$ gain is a natural robustness measure that captures the effect of gradient errors on the whole trajectory of the algorithm. % unlike previous approaches in the literature. %That being said, our results would extend to z_k = x_
Here, we also require $\overline{H}(x_*)= 0$ so that our bounds are tight in the errorless case.\footnote{Since the optimum $x_*$ is a fixed point of the iterations, when $x_0=x_*$ and $w_k = 0~ \forall k$, then $\sum_{k\geq 0 }( f(x_k) -f(x_*) )=0$.} Optimizing the parameters of control systems numerically with respect to multi-objective criteria involving $\ell_2$ gains and the speed of convergence is popular
%To our knowledge, $\ell_2$ gains were never used to design parameters of momentum algorithms for achieving (approximate) Pareto-optimal rate and robustness performance despite the popularity of such approaches 
for designing control systems that are robust to deterministic errors \cite{gumussoy2009multiobjective, zhou1996robust}, however such approaches are not yet fully applied to the problem of designing the parameters of first-order methods. Our aim is to fill this gap in the context of strongly convex smooth minimization. Our main contributions are as follows:\looseness=-1

First, we focus on strongly convex quadratics. In this case, the gradient has linear growth and the dynamical system corresponding to GMM iterations is linear. For linear systems, it is known that  $L_{2,*}$ coincides with the $H_\infty$ norm which is a key metric for assessing robustness of linear systems (see e.g. \cite{zhou1998essentials,zhou1996robust,burke2006hifoo}). Existing standard off-the-shelf algorithms with provable global convergence guarantees for $H_\infty$ norm computation %can only estimate the $\ell_2$ gain numerically, and this 
require solving algebraic matrix Riccati equations or eigenvalue problems requiring $\mathcal{O}(d^3)$ operations at each step \cite{boyd1990regularity,hinrichsen1991stability} and for general systems $H_\infty$ norms are not explicitly known. However, exploiting the block diagonal structure of the GMM updates in the frequency domain, %Using this connection and the frequency domain representation of the $H_\infty$ norm, 
we provide a closed-form formula for the $H_\infty$ norm and (hence for $L_{2,*}$) for any parameter choice $(\alpha,\beta,\nu)$ in Thm. \ref{thm-h-inf}, when GMM is globally convergent (otherwise $L_{2,*}$ is infinite). Using this formula, we characterize robustness of TMM, HB, GD, NAG methods for any choice of parameters including standard parameters previously used in the exact gradient case. To our knowledge, our results are the first to explicitly estimate $L_{2,*}$ for momentum methods as a function of parameters. Furthermore, in Thm. \ref{thm-h-inf}, we show the lower bound $H_\infty = L_{2,*}\geq \frac{1}{\sqrt{2\mu}}$ where $\mu$ is the strong convexity constant and characterize all possible choice of GMM parameters that attain the lower bound. This set of parameters corresponds to the most robust performance (in the sense of $\ell_2$ gain), revealing a fundamental lower bound on what level of robustness is achievable and when it can be achieved. Also, our results characterize the fundamental trade-offs between the linear convergence rate (that measures the performance in the exact gradient setting) characterized by the spectral radius of the iteration matrix and the worst-case robustness (measured by $L_{2,*}$). In particular, the stepsize $\alpha = 2/(L+\mu)$ leads to the fastest rate for GD but is not the most robust, whereas the stepsize $\alpha = 1/L$ is the most robust, but is not the fastest. We show NAG can be at the best robustness level, while achieving an accelerated $\rho = 1 -\Theta(\frac{1}{\sqrt{\kappa}})$ rate in the exact gradient setting where $\kappa = L/\mu$ is the condition number. HB with standard parameters can admit a faster (smaller) rate than NAG,  but at the cost of worsened robustness. Motivated by these observations, we develop the \emph{robustly stable heavy-ball method} that can achive the best robustness level while being faster than NAG by up to a constant factor (Prop. \ref{prop-robust-hb}). We also propose the \emph{robustly stable gradient descent} that is the fastest version of GD while being at the best robustness level. To our knowledge, this is the first time the $\ell_2$ gain is %xplicitly computed as a function of the parameters of a first-order algorithm and is 
used to design the parameters of momentum methods to achieve Pareto-optimal robustness with respect to convergence rates, i.e. for achieving best robustness level possible for a given rate. In addition, we construct gradient error sequences that correspond to the worst-case performance explicitly using a frequency domain analysis (Prop. \ref{prop-worst-case-noise}). Such a noise depends on the parameters and therefore on the underlying algorithm and, to our knowledge, was not known previously for a momentum algorithm when subject to worst-case square-summable noise.\looseness=-1

%computed as a function of the parameters of a first-order optimization algorithms. 
%used as a robustness measure for designing the parameters of an optimization algorithm
%GD with the fastest linear rate 
%
% We can characterize robustness of TMM, HB, GD, NAG methods for any choice of parameters. 
%
%We find that GD with the fastest rate, is not optimally robust. We propose the RS-GD method that can achieve the fastest rate while being optimally robust. We find that NAG with standard choice of parameters can achieve an accelerated rate at the best robustness level. HB can achieve a faster rate but at a degraded robustness level. Therefore, we propose RS-HB method which is the most robust possible, while being slightly faster than NAG. 

%It is known that $H_\infty$ norm corresponds to the multiplicative inverse of the stability radius, which is another robustness measure for dynamical systems. 
Second, for quadratic objectives, we consider multiplicative gradient noise where the norm of the gradient error is a fraction $p \in (0,1)$ of the norm of the gradient and the gradient error is a linear function of the iterates $x_k$ with possibly complex entries. In Sec. \ref{subsec-multiplicative-noise}, we show that the $\ell_2$ gain of GMM is related to the maximum amount of multiplicative gradient noise that can be tolerated to avoid divergence, i.e. the size of the multiplicative noise that can destabilize GMM is inversely proportional to $L_{2,*}$. %i.e. designing the GMM parameters with a small $L_{2,*}$ also enables us to have better robustness to multiplicative noise. 
%In this setting, the norm of the gradient error is a percentage $p \in (0,1)$ of the norm of the gradient and the gradient error is a linear function of the iterates $x_k$.
The results are achieved by leveraging the connections between the $\ell_2$ gain and the \emph{complex stability radius}. The latter is an alternative robust stability measure for a linear system (see e.g. \cite{gurbuzbalaban2012theory,hinrichsen-pritchard}), and in our context it measures the size of the linear perturbations that are necessary to destabilize the GMM iterations, allowing the perturbations to have complex values. If the perturbations are restricted to be real-valued, then the resulting robustness measure is called the \emph{real stability radius} \cite{qiu1995formula}. If the robustness is not at the best level, i.e. if $L_{2,*} = H_\infty>\frac{1}{\sqrt{2\mu}}$, then we can construct a worst-case multiplicative noise sequence (with potentially complex-valued entries) that will make GMM divergent for the relative noise level $p_* =\frac{1}{L_{2,*} \sqrt{2\mu} } \in (0,1)$ and show that when the noise level $p < p_*$ then GMM will converge (Remark \ref{rem-mult-noise-gmm}). Furthermore, we show that for most of the common choices of parameters, the worst-case multiplicative noise can be constructed to have real entries. The latter result is achieved by showing that the complex stability radius is equal to the real stability radius in many cases for GMM (Thm. \ref{thm-real-hinf}, Coro. \ref{coro-real-stab-equals-complex-stab-cases}). To our knowledge, these are the first explicit characterizations of the complex and real stability radii for momentum methods. Our quadratic results for robustness also serve as lower bounds on the worst-case robustness we can expect for more general strongly convex functions (Coro. \ref{coro-real-hinf-equals-complex-sometimes}, Remark \ref{rem-lower-bound}).\looseness=-1 %(eqn. \eqref{ineq-lower bound for minimal l2 gain}).

Third, we consider strongly convex functions with Lipschitz gradient and provide explicit bounds on $L_{2,*}$ for GD and NAG, leveraging strong convexity and smoothness of the objective where we consider both distances to the optimum and suboptimality as Lyapunov functions (Prop. \ref{prop-hinf-gd-bound} and Prop. \ref{prop-hinf-agd-bound}) to get tighter bounds. These upper bounds are tight when the stepsize is sufficiently small in the sense they cannot be improved more than a small constant factor (Remark \ref{remark-tightness-nag}). For more general GMM methods admitting arbitrary parameters $(\alpha,\beta,\nu)$, we provide a scalable matrix inequality-based approach
% that generalizes our explicit analysis for GD and NAG. This approach
that allows us to generate an upper bound on the $\ell_2$ gain numerically, %n a computationally tractable manner 
provided that a small ($4\times 4$) matrix inequality (MI) which depends on the GMM parameters as well as several auxilliary variables holds (Thm. \ref{thm-hinfty-bound}). Given GMM parameters, this allows us to efficiently estimate the $\ell_2$ gain by a simple grid search over the auxilliary variables. 
%numerically by checking the feasibility of a small ($4\times 4$) matrix inequality (MI). 
Our MI approach can be viewed as a generalization of our explicit analysis for GD and NAG; because for GD and NAG we can construct the auxilliary variables for which the MI holds by manual computations (Remark \ref{rema-MI-explicit-bounds}) and our MI approach can recover the same explicit bounds we obtained for GD and NAG. %where we can construct the auxilliary variables by handmade computations for which the MI holds (Remark \ref{rema-MI-explicit-bounds}). % and the MI approach yields the same explicit bounds we obtained for GD and NAG. %because we can show that the MI is feasible by handmade computations and the MI approach yields the same explicit bounds we obtained for GD and NAG.
Here, as a proof technique, we first obtain finite-time bounds for the cumulative suboptimality of the form $$\sum_{k=0}^K f(x_k) - f(x_*)\leq \gamma^2 \sum_{k=0}^K \|w_k\|^2 +\overline{H}(x_0)$$ for some $\gamma>0$ in terms of the errors $\{w_k\}_{k=0}^{K}$ seen up to step $K$ and a term $\overline{H}(x_0)$ depending on the initialization, where the validity of the MI ensures a decay in the Lyapunov function, otherwise the $\ell_2$ gain may be infinite. % when feasibility of the MI allows us to obtain a decay in a generalized Lyapunov function that combines suboptimality with distance to the solution.% where the feasibility of the MI is required to obtain a decay in the Lyapunov function. %, where the feasibility of the MI is equivalent to obtaining a decay for our generalized Lyapunov function that combines function values with distance squared to the solution. 
Then, we let $K\to\infty$ to achieve bounds of the form \eqref{eq-cum-dist-squared} where $\gamma$ serves as an estimate of the $\ell_2$ gain. From the convexity of $f$, it follows that these results directly imply $f(\bar{x}_K) - f(x_*) \leq \frac{  \gamma^2 \sum_{k=0}^K \|w_k\|^2}{K+1} + \frac{\overline H(x_0)}{K+1}$ for the averaged iterates $\bar{x}_K := \frac{x_0 + ... + x_K}{K+1}$ and as such we obtain new non-asymptotic ergodic convergence rate results for inexact GMM methods %where the constant $\gamma_K$ is optimal up to a constant for large $K$ and for some parameter choices and 
%$H(x_0)=0$ for $x_0=x_*$
 (Coro. \ref{coro-ergodic-rates} and Coro. \ref{coro-ergodic-rates-gmm}). It can also be seen that the right-hand side of these performance bounds stay bounded when %for non-square summable errors as long as the squared $\ell_2$ norm of the errors over the horizon $K$ grows at most linearly in $K$, i.e. if
$\sum_{k=0}^K \|w_k\|^2/(K+1)$ is bounded as $K\to\infty$. Therefore, our approach leads to performance bounds for gradient errors $w_k$ that are bounded (but not necessarily square-summable) as well, where we find that $\ell_2$ gains can also serve as a robustness measure beyond square-summable gradient errors (Coro. \ref{remark-non-square-summable}). In addition, our results highlight the trade-offs between convergence rate and worst-case robustness. In light of these results, we discuss how our characterizations of robustness can be stated as %placed in 
a small-scale optimization problem for selecting the GMM parameters to systematically achieve a desired trade-off (Sec. \ref{subsec-trading-rate-robustness}). Finally, in Sec. \ref{sec-num-experiments}, we provide numerical experiments that demonstrate that we can design GMM parameters to achieve these trade-offs, illustrating our results.\looseness=-1
%\vspace{-0.12in}
\section{Related work.}
When the gradient error is relative, Gannot \cite{gannot2022frequency} obtained linear convergence rates for the inexact gradient descent for objectives with a sector-bounded non-linearity, and linear convergence rates for inexact TMM for strongly convex smooth functions based on a frequency-domain analysis. In the relative noise setting, \cite{de2020worst} obtained worst-case convergence analysis of inexact gradient and Newton methods using the semi-definite programming performance estimation (PEP) technique developed by Drori and Teboulle \cite{drori2014performance}. The authors also studied inexact gradient descent subject to exact line search \cite{de2017worst}. In another line of work, Friedlander and Schmidt \cite{friedlander2012hybrid} consider $\mu$-strongly convex functions $f$ that are $L$-smooth (i.e. the gradient $\nabla f$ is $L$-Lipschitz) with inexact gradient descent subject to additive gradient errors $w_k$ at step $k$ admitting a bound $\|w_k\|^2 \leq B_k$. The authors show that for any (sub) linearly decreasing sequence $B_k$ and the choice of $\alpha =1/L$,  the algorithm has a sublinear (linear) convergence rate. More specifically, the authors show with an asymptotic analysis that $f(x_k) - f(x_*) = \mathcal{O}(C_k)$ as $ k\to\infty$ where $x_*$ is the minimum with $C_k = \max\{B_k, (1-\frac{\mu}{L}+\varepsilon)^k \}$ for any $\varepsilon < \frac{\mu}{L}$ although universal constants are not explicitly given and the analysis requires $B_k$ to be monotonically decreasing. Schmidt \emph{et al.} \cite{schmidt2011convergence} obtain guarantees for inexact proximal gradient and accelerated proximal gradient methods for optimizing the sum of a smooth convex function and a non-smooth convex function; for strongly convex objectives their result shows that if the errors decay to zero linearly sufficiently fast, then the algorithms' convergence rates will be comparable to those in the exact gradient settings. Luo and Tseng \cite{luo1993error} obtained asymptotic linear convergence under a local error bound condition by choosing the stepsize accordingly to obtain monotonic decay in the distances to the solution. Bertsekas and Tsitsiklis \cite{bertsekas2000} show that inexact gradient methods converge with square-summable relative and additive errors satisfying some assumptions. Many others also studied inexact gradient descent methods in various settings (see e.g.  \cite{bertsekas2015convex,bertsekas2000,iag,bertsekas2011incremental,gurbuzbalaban2019convergence,gurbuzbalaban2015globally} and the references therein), although robustness of inexact momentum-based methods such as NAG, HB or TMM to deterministic absolute errors is relatively understudied. 

Among the existing work that studied momentum methods, Devolder's Ph.D. thesis \cite{devolder2013exactness} and the related publications \cite{devolder2014first,devolder2013first,devolder2013intermediate,devolder2011stochastic}) considered an oracle model for first-order information that can capture inexactness in the gradients as well as function evaluations. Devolder \cite{devolder2013exactness} considers the primal gradient method (PGM) with the common stepsize $\alpha = 1/L$ and shows that for a convex objective $f$, the suboptimality of the averaged iterates $y_k$ satisfies $f(y_k) - f_* \leq \mathcal{O}(1/k) + \delta$ where $\delta$ is the (persistent) oracle error encountered at ever step. A similar result is shown for the dual gradient method (DGM), however for fast (accelerated) gradient methods (FGM), the faster convergence rate  $\mathcal{O}(1/k^2)$ is accompanied by a worse error term $\mathcal{O}(k\delta)$ that grows with iterations and the overall bound $f(y_k) - f_* \leq \mathcal{O}(1/k^2) + \mathcal{O}(k\delta)$ can be worse than that of the gradient method. Furthermore, it is shown that for convex objectives the accumulation of errors is unavoidable \cite[Sec. 4.8]{devolder2013exactness} in the sense that it is not possible to achieve acceleration with $\mathcal{O}(1/k^2)$ rate and not suffer from robustness issues, i.e. from error accumulation. Devolder also develops intermediate gradient methods (IGM) that can interpolate between the performance of the gradient method and the accelerated gradient method under different parameter choices. For a power decay stepsize/parameter rule, their method admits the intermediate performance $f(y_k) - f_* \leq \mathcal{O}({\frac{1}{k^p}}) +  \mathcal{O}(k^{p-1}\delta)$ for $p\in [1,2]$ as $k\to\infty$. 
 %These results show that for persistent additive errors that do not fade away, we do not expect to have an improved performance for large $k$ in general with accelerated methods for convex objectives. 
For strongly convex smooth objectives, it is also shown that accelerated gradient methods with standard parameters are also less robust compared to gradient descent %but the noise accumulation is not observed 
while having a faster decay rate,
i.e. we have $f(y_k) - f(y_*) = \mathcal{O}(e^{-k \kappa}) + \mathcal{O}(\delta) $ for both PGM and DGM; while for FGM we have $f(y_k) - f(y_*) \leq \mathcal{O}(e^{-k \sqrt{\kappa}}) + \mathcal{O}(\sqrt{\kappa}\delta) $. However, these results hold under an oracle noise model which does not apply to additive gradient noise on unbounded domains \cite[Remark 4.2]{devolder2013exactness} where the gradient errors are persistent  and non-square-summable. %Therefore, these results are not applicable to our setting. %where we have a different noise setup. %Another difference is that in our setup, the errors are not persistent but of decaying nature due to square-summability.%Also, our results are of ``positive nature" in the sense that they show that acceleration is possible without suffering from robustness
%Another difference is that our results show that acceleration and robustness can be achieved simultaneously when the errors are square-summable, therefore not persistent. 
In another line of work, d'Aspremont shows that optimal complexity of Nesterov’s method is preserved with proper averaging of the iterates, when the gradient admits uniformly bounded errors (by a parameter $\delta>0$) i.e. $f(y_k) - f(y_*) \leq \mathcal{O}(1/k^2) + \mathcal{O}(\delta)$. However, this result requires the constraint set to be convex and bounded and considers errors that are persistent \cite[Thm. 2.2]{daspremont}. Therefore, existing results from the literature do not apply to our setting where we are primarily interested in the worst-case behavior of a generalized class of momentum algorithms on an unbounded domain when the gradient errors are not persistent but are square-summable. \looseness=-1%for unconstrained minimization problems %f
\textbf{Notation.} A function $f:\mathbb{R}^d\to\mathbb{R}$ is called $\mu$\emph{-strongly convex} if the function $x \mapsto f(x)-\frac{1}{2}\mu \|x\|^2$ is convex on $\mathbb{R}^d$ for some constant $\mu>0$. A continuously differentiable function $f:\mathbb{R}^d\to\mathbb{R}$ is \emph{$L$-smooth} if its gradient is $L$-Lipschitz, i.e. satisfies $\|\nabla f(x) - \nabla f(y) \|\leq L \|x-y\|$ for all $x,y \in \mathbb{R}^d$. Let $\Cml$ denote the set of all functions $f$ that are both $\mu$-strongly convex and $L$-smooth at the same time. Such functions satisfy the inequalities
%\vspace{-0.1in}
\beq
\frac{\mu}{2}\|x-y\|^2 \leq f(x) - f(y) - \nabla f(y)^T (y-x) \leq \frac{L}{2}\|x-y\|^2,
\label{ineq-strcvx-smooth}
%\vspace{-0.1in}
\eeq
\cite[Sec. 9.1.2.]{boyd2004convex}. Due to strong convexity, $f$ admits a unique global minimum on $\mathbb{R}^d$, which we will denote by $x_*$. Let $f_* := f(x_*)$ denote the minimum value of $f$. We assume $\mu < L$ throughout the paper; otherwise, $\mu=L$ and the class $\Cml$ is trivial. Let $I_d$ and $0_d$ denote the $d\times d$ identity and zero matrices respectively; we drop the subscript $d$ in some cases if it is clear from the context. We use $0_{d_1 \times d_2}$ to denote the $d_1\times d_2$ rectangular matrix with all zero entries. The spectral radius $\rho(A)$ of a square matrix $A$ is the largest modulus of the eigenvalues of $A$. A block diagonal matrix $D$ with $i$-th diagonal block $D_i$ will be denoted as {$\underset{i=1,..,d}{\Diag} \left[D_i\right]$}. Let $\mathbb {C} ^{d}:=\mathbb {C} \times \mathbb {C} \times \cdots \times \mathbb {C}$ denote the $d$-fold Cartesian product of the complex plane $\mathbb{C}$. Given matrix $A$, let $A_{ij}$ denote the entry on the $i$-th row and $j$-th column of $A$. We let $\mathbb{R}_{+}$ denote the set of non-negative reals. For functions $g:\mathbb{R}_{+} \to \mathbb{R}_{+}$ and $h:\mathbb{R}_{+} \to \mathbb{R}_{+}$, we say $g = \Theta(h(u))$ as $u\to 0$ if there exists positive constants $c_0, c_1, u_0 >0$ such that $c_0 h(u) \leq g(u) \leq c_1 h(u)$ for all $u$ with $u< u_0$. For a matrix $A$, $\|A\|$ denotes the 2-norm, i.e. the spectral norm of $A$. %For a square summable sequence $\{w_k \}_{k\geq 0}$ with $w_k \in \mathbb{R}^d$, its $\ell_2$ norm is defined as $\|w\|_{\ell_2(\mathbb{R}^d)} := \sum_{k\geq 0}\|w_k\|^2$ where $\| \cdot\|$ denotes the Euclidean norm. 
 %where the superscript $*$ denotes the Hermitian transpose. 
%A real-valued function $h: W \to \mathbb{R}$ defined on an interval $W$ is \emph{quasi-convex} if the lower level sets $L_a := \{x \in \mathbb{R} : f(x)\leq a\}$ are convex for every $a\in \mathbb{R}$. 
Let $\mbox{arccos}: [-1,1] \to [0,\pi]$ denote the inverse of the cosine function, i.e. $y=\mbox{arccos}(x)$ if  $x = \cos(y)$ and $y \in [0,\pi]$. For a complex vector $v\in \mathbb{C}^d$, $v^T$ denotes the transpose and $v^H$ denotes Hermitian transpose where $\|v\|:= \sqrt{v^H v}$. We use $A \otimes B$ to denote the Kronecker product of the matrices $A$ and $B$.\looseness=-1
%\vspace{-0.05in}
\section{Preliminaries.}\label{sec-prelim}
%\mtodo{For dimension, in some places we use $d$ sometimes $n$. make consistent}
Let $f\in\Cml$ be given. We consider the unconstrained optimization problem of minimizing $f$ on $\mathbb{R}^d$. We consider the following class of generalized momentum methods (GMM): %Generalized Momentum Method (GMM) \cn consists of the following iterations
%\begin{subequations} \label{RBMM}
\begin{eqnarray}
x_{k+1}=x_{k}-\alpha \nabla f( y_k ) +\beta(x_{k}- x_{k-1}), \quad
y_{k}&=&x_k+\nu (x_k- x_{k-1}), \label{RBMM}
\end{eqnarray}
starting from the initialization $x_0=x_{-1} \in \mathbb{R}^d$ which admit three parameters $\alpha, \beta$ and $\nu$: The parameter $\alpha>0$ is the stepsize, where $\beta,\nu\geq 0$ are the \emph{momentum parameters}. These methods were previously studied in the exact gradient setting (see e.g. \cite{lessard2016analysis,hu2017dissipativity}) and referred to as GMM in \cite{can2022entropic}. 
GMM generalizes a number of momentum-averaging based first-order algorithms. If $\nu=0$, GMM is equivalent to Polyak's HB method \cite{polyakintroduction}. If we choose $\beta=\nu$, it recovers the NAG method \cite{nesterov2018lectures}. On the other hand, when $\beta=\nu=0$, this method reduces to GD. TMM, another momentum method which admits faster convergence rate $\rho$ than NAG (up to a constant factor) in the absence of noise, corresponds to a particular choice of parameters \cite{scoy-tmm-ieee}. Other choices of parameters have also been useful for minimizing the risk and the tail probabilities associated with suboptimality when the gradients are subject to stochastic noise \cite{can2022entropic}. 

There has been a growing literature about the reformulation of optimization algorithms as dynamical systems, where tools from control theory can be leveraged to obtain convergence rates of existing algorithms in the exact gradient as well as stochastic gradient settings \cite{hu2017dissipativity, lessard2016analysis,aybat2018robust,can2022entropic}. %as well as to design the parameters of algorithms \cn. 
%These works are in the exact gradient setting. U
%nbiased stochastic noise has been studied in \cite{can2022entropic}, 
However, the robustness of GMM algorithms to worst-case deterministic noise is relatively less understood. In our work, we will also consider such a reformulation and apply tools from robust %Linear-Exponential-Quadratic (LEQ) 
control theory to design the parameters of GMM algorithms to obtain systematic trade-offs between their sensitivity to worst-case deterministic errors in the gradients (in terms of the $\ell_2$ gain) and their convergence rate. The setup is as follows: %We focus on additive errors where %n this work, we will first study the performance of GMM algorithms subject to deterministic worst-case additive gradient noise, i.e.
Instead of the actual gradient $\nabla f(y_k)$ at step $k\geq 0$, we  assume we have access to its inexact version 
%\beq 
$\widetilde \nabla f(y_k, w_k) = \nabla f(y_k) + w_k$,
%\label{eq-gradient-noise}
%\eeq 
where $w_k \in \mathbb{R}^d$ is the additive gradient error at step $k$. In this context, throughout this work, we will use the terms \emph{errors} and \emph{noise} interchangably. We start with reformulating GMM iterations \eqref{RBMM} subject to additive errors as a dynamical system:
%\begin{subequations}
%\label{Sys: RBMM}
\begin{align} 
&\xi_{k+1}=A \xi_{k}+Bu_k, \quad y_k = C \xi_k, \quad u_k = \widetilde \nabla f(y_k, w_k)  = \nabla f(y_k) + w_k, \label{Sys: RBMM}
\end{align}
%\end{subequations}
where $\xi_k := \begin{bmatrix}  x_k^T  &
                           x_{k-1}^T 
          \end{bmatrix}^T$ is the \emph{state vector} which contains the last two iterates $x_k, x_{k-1}$ at time $k$ and $A, B$ and $C$ are system matrices defined as $A = \tilde{A}\otimes I_d$, $B = \tilde{B}\otimes I_d$ and $C = \tilde{C}\otimes I_d$ with
%\beq \xi_k := \begin{bmatrix}  x_k  \\
%                           x_{k-1} 
%          \end{bmatrix}, \quad A = \tilde{A}\otimes I_d,  \quad B = \tilde{B}\otimes I_d, \quad C = \tilde{C}\otimes I_d, 
%\label{def-ABC}
%\eeq
%with
\begin{align}\label{def: system mat for TMM}
\tilde{A}:=\begin{bmatrix} 
(1+\beta) & -\beta  \\ 
1 & 0
\end{bmatrix} , \quad \tilde{B}:=\begin{bmatrix} 
-\alpha \\ 
0
\end{bmatrix}, \quad \tilde{C}:=\begin{bmatrix}
(1+\nu) & -\nu 
\end{bmatrix}.
\end{align} 
%where $A, B$ and $C$ are system matrices with appropriate sizes,
%and
%\begin{equation}\xi_k := \begin{bmatrix}  x_k  \\
%                           x_{k-1} 
%          \end{bmatrix}
%\label{def-ksi-k}
%\end{equation}
%is the \emph{state vector} which contains the iterate at step $k$ and the previous iterate. 
%The vector $u_k$ is the noisy gradient evaluated at $y_k$. 
%\vspace{-0.1in}
We next discuss the notion of $\ell_2$ gain for this GMM system and its relevance to optimization.
\subsection{$L_2$ gain as a robustness metric.}\label{subsec-l2-gain-def} %The $H_\infty$ norm of a dynamical system is a measure that quantifies its robustness to worst-case (square-summable) noise. 
%Before defining the $H_\infty$ norm of the dynamical system corresponding to GMM, 
%Before defining the $\ell_2$ gain of the system corresponding to (noisy) inexact GMM iterations formally,
We first introduce the following assumption on the gradient noise, which says that it is deterministic and square-summable. 
%In the following, we consider
%the following assumption 
%$H_\infty$ norm of the dynamical system corresponding to GMM iterations is a way to quantify its robustness to worst-case square summable deterministic noise. In the following, we consider GMM iterations when the gradients are corrupted by a sequence of square-summable deterministic noise $\{{w}_k \}_k$. %instead of the stochastic noise $\{w_k\}$ considered previously:
%{\color{red}
%\begin{subequations}
%\label{Sys:GMM-deterministic-noise}
%\begin{eqnarray} 
%\hxi_{k+1}&=&A \hxi_{k}+B \nabla f(C\hxi_k) +  B  \hat{w}_k,\\
%\hz_k &=& F(\hxi_k) %\sqrt{f(\hat{x}_k) - f(x_*)} %F(\hxi_k)
%\end{eqnarray}
%\end{subequations}
%where 
%\beq \hat{\xi}_k := \begin{bmatrix}
%    \hx_k  \\
%    \hx_{k-1}
%\end{bmatrix}
%\eeq
%}
%{\color{red}We note that the only difference between these two systems \eqref{Sys:stoc-RBMM} and \eqref{Sys: TMM-deterministic-noise} is that the noise is deterministic in the latter whereas it is stochastic in the former. }
%We make the following assumption on the deterministic noise which says that it is square-summable. %, later in Section \cn, we will also consider more general deterministic noise that does not satisfy this condition.
 \begin{assum}\label{assump-deter-noise-real} Consider the noise vector ${w}_k$ that represents the gradient noise in \eqref{Sys: RBMM}. For every $k\geq 0$, ${w}_k \in \mathbb{R}^d$, ${w}_k$ is deterministic and $\|w\|_{\ell_2(\mathbb{R}^d)} = (\sum_{k=0}^\infty \|w_k\|^2) ^{1/2} < \infty$.
 \end{assum}
 
%We assume that the sequence $\hat{d}_k$ is in $\ell_2$, i.e. it is squaree summable with $\|d\|_{\ell_2} := (\sum_{k=0}^\infty d_k^2)^{1/2}$ being finite. Then, we
%This deterministic noise structure is clearly very different than the randomized (Gaussian) noise we considered previously in Assumption \ref{assump-gaussian-noise}. We will not require both assumptions to hold simultaneously. 
Notice that we can rewrite the noisy GMM iterations \eqref{Sys: RBMM} as
%\begin{subequations}
%\label{Sys:stoc-RBMM}
%\begin{eqnarray} 
%%\xi_{k+1}&=&A \xi_{k}+Bu_k +  B  \Sigma_{\pi} w_k,\\
%%y_k &=& C \xi_k, \\ 
%%u_k &=& \nabla f(y_k),
%\xi_{k+1}&=&A \xi_{k}+B \nabla f(C\xi_k) +  B  w_k,\\
%z_k &=& F(\xi_k),
%\end{eqnarray}
%\end{subequations}
%Consider the noisy GMM dynamics
%\begin{subequations}
%\vspace{-0.05in}
\begin{eqnarray} 
\label{Sys:stoc-RBMM}
%\xi_{k+1}&=&A \xi_{k}+Bu_k +  B  \Sigma_{\pi} w_k,\\
%y_k &=& C \xi_k, \\ 
%u_k &=& \nabla f(y_k),
\xi_{k+1}&=&A \xi_{k}+B \nabla f(C\xi_k) +  B  w_k,
%\vspace{-0.15in}
\end{eqnarray}
%\end{subequations} 
%the matrices $A,B$ and $C$ are as in \eqref{def: system mat for TMM} 
starting from an initialization ${\xi_0} = \begin{bmatrix} x_0^T & x_{-1}^T \end{bmatrix}^T $.  For simplicity of the presentation, throughout the paper we take $x_{-1} = x_0$, but our results would extend to an arbitrary choice of $x_{-1}$ in a straightforward manner. Without noise, i.e. when $w_k = 0$ for all $k$, and when GMM is convergent, the fixed point of the iterates \eqref{Sys:stoc-RBMM} is
%\beq
 $\xi_* = \begin{bmatrix} x_*^T & x_*^T \end{bmatrix}^T$.
%\label{def-xi-star}
%\eeq
%We will assume in this subsection that the noise sequence $\{ w_k \}_{k\geq 0}$ obeys Assump. \ref{assump-deter-noise-real} where the entries of the noise vector is real-valued. Later in this section, when we introduce the $H_\infty$ norm,  we will allow the entries of the noise vector $w_k$ to be complex-valued. 
We will also consider an output sequence $z_k$ satisfying
\beq z_{k} &=& F(\xi_{k}), \quad \|z_{k}\|^2 = \| F(\xi_{k})\|^2 =  {f(x_{k})-f(x_*)},
\label{eq-output-choice} 
\eeq
for $k\geq 0$ where the output $z_k$ is defined through a map $F:\mathbb{R}^{2n} \to \mathbb{R}^p$ for some $n\geq 1$ satisfying the latter inequality, i.e. the squared norm of $z_k$ is equal to the suboptimality $f(x_k)-f(x_*)$ at step $k$.  Therefore, with this choice of the $F$ map (whose exact definition will be provided later in Section \ref{sec-quad}), the $z_k$ sequence can be viewed as an \emph{error signal}, and its $\ell_2$ norm squared will coincide with the cumulative suboptimality over the iterations. In fact, for any given $x_k $ and $f\in\Cml$, we can write $f(x_k) - f(x_*) = \frac{1}{2}(x_k-x_*)^T Q_k (x_k-x_*)$ for a matrix $Q_k$ which can be viewed as an average Hessian of $f$ (along the line segment from $x_*$ to $x_k$) that satisfies $ \mu I_d \preceq Q_k \preceq L I_d$ (see e.g. [Sec. 1.1.3]\cite{polyakintroduction}). Therefore, $z_k$ satisfying \eqref{eq-output-choice} has the property that $\|z_k\| = \frac{1}{\sqrt{2}}\|Q_k^{1/2}(x_k - x_*)\|$ and as such it can be interpreted as a \emph{weighted distance to the optimum}. Here, $z_k$ is introduced only for analysis purposes, and is not actually directly computed over iterations. It can also be seen from our proof techniques that our theory would naturally extend to other choices of $F$ that are affine in the state variable $\xi_k$, for instance if we choose ${z}_k = x_k - x_*$ then our framework applies with minor modifications.\looseness=-1 %the $\ell_2$ norm of the error signal $z_k$ will track the cumulative distance squared to the solution.%Here, as we are in the optimization algorithms setting; we choose the map $F$ in a particular way so that the output captures information about suboptimality. %More specifically,
%We consider a particular choice of. More specifically, in Section \cn, we will choose $F$ to satisfy %W
%WFe will specify the choice of $F$ later. %In this paper, we are interested in the suboptimality $f(x_k) - f(x_*)$ and therefore we will choose the output $z_k$
%It will be chosen in a particular way so that the output $z_k$ captures information about the suboptimality $f(x_k) - f(x_*)$ in the  
%\beq \|z_k\|^2 = \|F(\xi_k)\|^2 =  {f(x_k)-f(x_*)},\eeq similar to the literature. % \cite[Remark 2.5]{fleming1995risk}. 
%More specifically, the norm square of $z_k$ is equal to the suboptimality $f(x_k)-f(x_*)$ at step $k$. Therefore, the $z_k$ sequence can be viewed as an \emph{error signal}, whose $\ell_2$ norm will coincide with the cumulative suboptimality over the iterations.
%\footnote{It can be seen from our proof techniques that our theory would naturally extend to other choices of $F$, for instance if we choose ${z}_k = F({x}_k) = {x}_k - x_*$ then the $\ell_2$ norm of the error signal $z_k$ will track the cumulative distance squared to the solution.}
%\begin{defi}{\textbf{$L_2$ gain.}} 

We say that the system \eqref{Sys:stoc-RBMM} with input $\{w_k\}_{k\geq 0} \in \ell_2(\mathbb{R}^d)$ (obeying Assumption \ref{assump-deter-noise-real}) and output $z_k$ has $L_2$ gain $\leq \gamma$ (from $\mathbb{R}^d$ to $\mathbb{R}^d$) if for every initialization $\xi_0 \in \mathbb{R}^{2d}$, there exists a function $H: \mathbb{R}^{2d} \to \mathbb{R}_+$ satisfying $H(\xi_*) = 0$ such that
%\vspace{-0.05in}
\beq \sum_{k\geq 0}  \big(f(x_k) - f(x_*) \big)  = \sum_{k\geq 0}\|z_k\|^2 \leq \gamma^2 \sum_{k\geq 0}\|w_k\|^2 + H(\xi_0), \quad \forall w_k \in \ell_2(\mathbb{R}^d),%\vspace{-0.65in}
\label{def-l2-gain-general}
\eeq
see e.g. \cite{tran2017qualitative, lin1996h,van2016l2}). Since we take $x_0=x_{-1}$, note that this inequality is equivalent to \eqref{eq-cum-dist-squared} by setting $H(\xi_0) = \overline{H}(x_0)$. If the inequality \eqref{def-l2-gain-general} holds, we have clearly $z_k \in \ell_2(\mathbb{R}^d)$. Roughly speaking, this inequality says %a smaller value of $\gamma$ shows 
that the input noise is not amplified by a factor more than $\gamma$ if we would compare the $\ell_2$ norms of the input and the output, albeit with a constant factor arising from the initialization. The $\ell_2$ gain is formally defined as 
\beq %%\vspace{-0.15in} 
L_{2,*} : = \inf\{ \gamma \in \mathbb{R} ~:~ \gamma \mbox{ satisfies } \eqref{def-l2-gain-general}\},
\label{def-optimal-L2-gain}%%\vspace{-0.25in}
\eeq which corresponds to the (smallest) best choice of $\gamma$ that leads to the tightest error bounds \cite{van2016l2}. 
 %In this work, we will estimate provide lower and upper bounds for the clas
%smallest such $\gamma$ is called the minimal $L_2 gain$ and will be denoted by 
 %Therefore, we can view the output sequence $z_k$ as an error signal, as its norm square is equal to the suboptimality $f(x_k)-f(x_*)$ at step $k$.
%\mtodo{Fleming-James says this formulation for $H_\infty$ as a ratio is for linear systems, and for nonlinear systems, we may need to define it according to \eqref{def-hinf-linear}. Check.} 
 %{\color{red} Initial cond. may need to be zero; see Hinrichsen-Son discrete-time stability radius paper.}
%\mtodo{According to \cite{hinrichsen1991stability}, it seems that the sequence ${d}_k$ should be complex-valued and be square summable.}
%Therefore, the $H_\infty$ norm of a system
%As its name suggests, \ell_2 gai can be viewed as the $\ell_2$ gain of the system from the noise input to the output. 
It is a measure of how much the output signal will be amplified in the worst case due to input perturbations %, if the perturbations $w_k$ to the system has a finite energy, i.e. a finite $\ell_2$ norm.
%how much cumulative errors can be in the worst case allowing deterministic perturbations 
$\{{w}_k\}$ with a finite energy (i.e. with a finite $\ell_2$ norm). 
%Therefore, it can be interpreted as a measure of sensitivity of a dynamical system to worst-case deterministic perturbations of finite energy. %if the system is started from the optimum at time zero. When the iterations start from the optimum subject to deterministic perturbations, 
%For analyzing GMM algorithms, we will define the output of the system as the suboptimality $f(\hx_k) - f(x_*)$. As $k$ goes to infinity, the errors ${w}_k$ will fade away as they are square-summable and consequently we have eventually the suboptimality $f(\hx_k) - f(x_*)$ going to zero. In this case, the quantity \eqref{eq-hinf-quad} can be interpreted %as a measure of how fast this transition happens where \eqref{eq-hinf-quad} can be thought 
%as the worst-case sensitivity of the underlying optimization algorithm to the deterministic perturbations that are square-summable.
The smaller $L_{2,*}$ of the dynamical system representation of GMM algorithm is, the more robust GMM is with respect to square-summable worst-case perturbations ${w}_k \in \mathbb{R}^d$. The value of $L_{2,*}$ will depend on the choice of the objective $f$ and the parameters $(\alpha,\beta,\nu)$. In particular, for making the GMM algorithm more robust with respect to such square-summable worst-case deterministic noise, a possible approach is to choose the parameters ($\alpha, \beta,\nu$) to make the $L_{2,*}$ %$\gamma$ 
value smaller. % which will be approach taken in this paper. 
%Since we are interested in the suboptimality at iteration $k$, throughout this paper we will set the output ${z}_k$ satisfying
%\beq \|{z}_k\|^2 = {f({x}_k) - f(x_*)}
%\label{def-out-seq}
%\eeq
%so that the $\ell_2$ norm of the output sequence ${z}_k$ will quantify the suboptimality accumulated over the iterations. However, our theory would extend for instance if we choose ${z}_k = F({x}_k) = {x}_k - x_*$ or more generally if we take $F$ as an affine function as the iterates $x_k$.
This will be the approach taken in this paper and we call the $L_{2,*}$ value of the dynamical system representation of GMM with respect to the output satisfying \eqref{eq-output-choice}, the \emph{(worst-case) robustness} of GMM.

\subsection{Relevance of $\ell_2$ gain to optimization.}\label{subsec-relevance} 
The idea of $\ell_2$ gain is illustrated in Fig. \ref{fig: transient-main} where we fix a noise sequence $\{w_k\}$ and visualize how the suboptimality trajectory could evolve over iterations. % if the iterations were started at the optimum.
For simplicity of the illustration, we consider a strongly convex function $f$ initialized at the optimum $\xi_0=\xi_*$ so that $f(x_0)-f_*=0$ where $f_*:=f(x_*)$ is the optimal value and $H(\xi_0)=H(\xi_*)=0$. We assume the parameters are such that without noise, GMM iterations are globally (linearly) convergent (otherwise the $\ell_2$ gain may not be finite). As we see in Fig. \ref{fig: transient-main}, the suboptimality will be increased in the early iterations due to the (worst-case) adversarial structure of the noise, that being said we expect the suboptimality go to ze-
\\
\noindent
\begin{minipage}[t]{0.55\textwidth}
%% PUT TEXT HERE
ro eventually. This is because by the assumption, the noise is square-summable and consequently is ``fading away" (i.e. for every $ \varepsilon$, there exists $K$ such that $\|w_k\|\leq \varepsilon$ for $k\geq K$). The $\ell_2$ gain squared ($L_{2,*}^2$) will be given by the worst-case accumulated suboptimality $\sum_{k\geq 0} [f(x_k) - f(x_*)]$ (which is approximated by %related to 
the area under the curve in Fig. \ref{fig: transient-main} up to a first-order term in the stepsize) normalized by the $\ell_2$ norm squared of the noise input. In particular, for the same noise budget ($\ell_2$ norm of the gradient errors), if the quantity $L_{2,*}$ is smaller, then this means that GMM iterations will converge
\end{minipage}\hfill
\begin{minipage}[t]{0.40\textwidth}
  \centering\raisebox{\dimexpr \topskip-\height}{%
  \includegraphics[width=\textwidth]{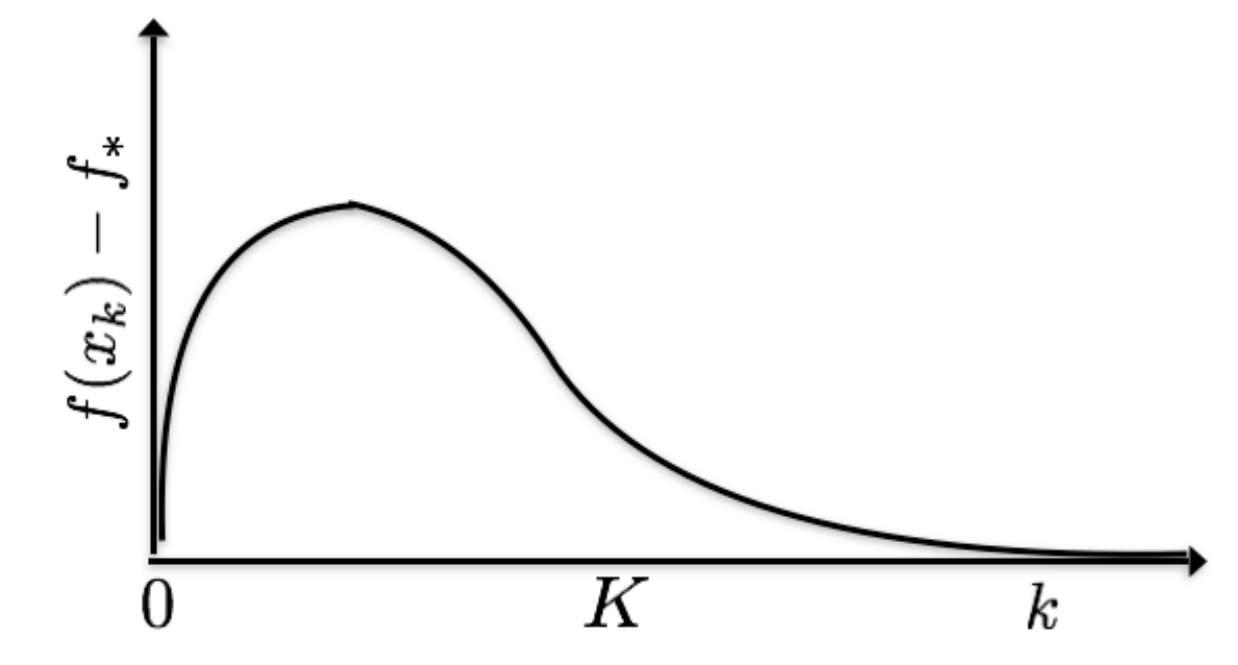}} 
  \captionof{figure}{Illustration of the $H_\infty$ norm \label{fig: transient-main}}
  \label{fig1}
\end{minipage}
%\begin{wrapfigure}{r}{0.3\textwidth}
 %   \centering
    %\FIGURE
   % \begin{figure}
    %\FIGURE
   % {\includegraphics[width=0.3\textwidth]{hinfinity-crop.pdf}}
    %{Illustration of the $H_\infty$ norm \label{fig: transient-main}}
    %{The suboptimality over iterations $k$.}
    %\end{figure}
%    \begin{wrapfigure}{r}{0.3\textwidth}
%     {\includegraphics[width=0.3\textwidth]{hinfinity-crop.pdf}}
%     {Illustration of the $H_\infty$ norm \label{test}}
%    \end{wrapfigure}
%     \caption{Illustration of the $H_\infty$ norm}
%\end{wrapfigure}
%the GMM algorithm is more robust to worst-case noise in the sense that 
to the optimum while having a smaller cumulative deviation ($\sum_{k\geq 0} [f(x_k) - f(x_*)]$) from the optimum along the way in the worst-case and hence will be more ``robust" to worst-case noise. If the iterations do not start from the optimum, we have a similar conclusion except that the cumulative suboptimality bounds will be shifted by a non-negative term $H(\xi_0)$ that depends on the initialization. 
%for the same noise level (in the $\ell_2$ sense).%, i.e. for any gradient error sequence with the same $\ell_2$ norm. 
%The $H_\infty$ norm captures a s similar effect while allowing the complex-valued noise structure. 
Also, as discussed in Sections \ref{subsec-multiplicative-noise} and \ref{subsec-bounds-on-l2-gain}, the multiplicative inverse of $L_{2,*}$ is related to the stability of GMM in the sense that the $\ell_2$ gain is related to how much multiplicative noise an optimization algorithm can tolerate before being divergent. These reasons make the $\ell_2$ gain fundamental to understanding the stability and robustness of optimization algorithms. %More generally, when the initial point $x_0$ is not necessarily the optimum, we will show that for GMM algorithms we can obtain guarantees of the form 
%\beq 
%\sum_{k=0}^\infty \big( \| {z}_k\|^2 - \gamma^2\|{w}_k\|^2 \big) \leq F(x_0) 
%\quad \forall {w}\in\ell_2, \label{def-hinf-linear-general}
%\eeq
%with $F(x_0) = c[f(x_0) - f(x_*)]$,  for any initial point $x_0$ and for some constants $\gamma, c>0$ that are independent from the initial point $x_0$.  For non-linear systems that are asymptotically stable (such as the GMM algorithm admitting with parameters leading to linear convergence in the abscence of noise), the smallest $\gamma$ that satisfies \eqref{def-hinf-linear-general} for some finite function $F(x_0)$ is called the $H_\infty$ norm (see e.g. \cite{tran2017qualitative, lin1996h}). 
%Next, we will consider stochastic noise and introduce the risk-sensitivity index as a performance metric. 
%\begin{figure}
%    \centering
%    \hfill
%    \quad\quad\quad\quad\quad\quad
%    \includegraphics[width=0.4\linewidth]{hinfinity-crop.pdf}
%  %  \includegraphics[width=0.4\linewidth]{figure-transient-2-crop.pdf}
%    \hfill
%    \caption{Illustration of the $H_\infty$ norm}
%    \label{fig: transient-main}
%\end{figure}

%In the next section, for strongly convex quadratic objectives, we show that the $\ell_2$ gain can be characterized with an explicit formula as a function of the parameters $(\alpha,\beta,\nu)$ and discuss the insights that can be obtained for the robustness of momentum algorithms. %Later in Section \ref{sec-strongly-convex}, we will study the case when $f\in \Cml$.

\section{Main results for quadratic objectives.}\label{sec-quad} In this section, we consider the special case when the objective is a strongly convex quadratic function of the form
    \beq\label{eq-quad}
        f(x):= \frac{1}{2}x^T Q x + p^T x + r,
    \eeq    
where $Q\in \mathbb{R}^{d\times d}$ is a positive-definite matrix, $p\in \mathbb{R}^d$ is a column vector and $r\in\mathbb{R}$. %First, we will discuss how the output $z_k$ %with system output $z_k$ can be written as
%\begin{align} 
%$\xi_{k+1}=A \xi_{k}+B(\nabla f(C\xi_k) + w_k)$ and $z_k = F(\xi_k)$,
%\label{system-TMM-with-output}
%\end{align}
%and from \eqref{eq-hinf-quad}, if the iterations start at the optimum (i.e. $\xi_0 = \xi_*$), we have
%\beq %\| G_{{w},{z}} \|_\infty 
%L_{2,*}
%\leq \gamma \iff \sum_{k=0}^\infty \big( f(x_k) - f(x_*) -  \gamma^2\|{w}_k\|^2 \big) \leq 0, \quad \forall {w}\in\ell_2(\mathbb{C}^d).
%\label{def-hinf-linear}
%\eeq
%\begin{subequations}
%\begin{eqnarray}
%x_{k+1}&=&x_{k}-\alpha [\nabla f(y_k)+ {w}_k] +\beta(x_{k}- x_{k-1}), \label{RBMM: Iter1-worst-case}\\
%y_{k}&=&x_k+\nu (x_k- x_{k-1}).  \label{RBMM: Iter2-worst-case}%\\
%\hat{z}_k &=& \sqrt{f(\hat{x}_k) - f(x_*)}\label{RBMM: Iter2-worst-case-out}
%\end{eqnarray}
%\end{subequations}
%where the output $\hat{z}_k$ is chosen such that the $\ell_2$ norm of the output sequence $\hat{z}_k$ will quantify the suboptimality accumulated over the iterations.
%We first introduce the ``centered" iterates
%\beq \xi_k^c:= \xi_k - \xi_*, ~ \xi_* := \begin{bmatrix} x_*^T\\ x_*^T\end{bmatrix},
%\label{def-centered-iter}
%\eeq
%where
%Note that the superscript $``c"$ is to highlight that these iterates are shifted to be ``centered" around the optimum; i.e. if GMM converges to the optimum, by definition $\xi_k^c\to 0$.   
%where the output $z_k$ is chosen according to \eqref{eq-output-choice}. 
First, we discuss how the output $z_k$ of the GMM system \eqref{Sys: RBMM} can be chosen to satisfy \eqref{eq-output-choice}. For this purpose, we consider the eigenvalue decomposition $Q = U \Lambda U^T$ of $Q$ where $\lambda$ is a diagonal matrix containing eigenvalues of $Q$ in increasing order, i.e. $Q_{ii} = \lambda_i$
where $\mu = \lambda_1 \leq \lambda_2 \leq \dots \leq \lambda_d = L$ are the eigenvalues of $Q$. Since 
$  f(x_k) - f(x_*) = \frac{1}{2} (x_k - x_*)^T Q(x_k - x_*) =  \frac{1}{2} (x_k - x_*)^T (U\Lambda U^T) (x_k - x_*), $
we have $\nabla f(C\xi_k) = Q(C\xi_k - x_*) = QC \xi^c_k $ where % $\xi_k^c$ defined as
\beq \xi_k^c:= \xi_k - \xi_* \quad \mbox{with} \quad \xi_* := \begin{bmatrix} x_*^T & x_*^T \end{bmatrix}^T,
\label{def-centered-iter}
\eeq
is the centered iterate. Note that the superscript $``c"$ is to highlight that these iterates are ``centered" around the optimum; i.e. if GMM converges to the optimum, by definition $\xi_k^c\to 0$. If we take
%In our discussion, we will choose $T$ as 
	\beq
	z_k = F(\xi_k) = T \xi_k^c, \quad T =\left[ \frac{1}{\sqrt{2}} \Lambda^{1/2} U^T \quad  0_d \right],
	\label{def-T}
	\eeq
then \eqref{eq-output-choice} holds as desired. In light of \eqref{def-T}, we can rewrite the GMM system \eqref{Sys: RBMM} with this output as	
%so that we have 
%  $$ \| \hz_k\|^2 = \|T(\hxi_k)\|^2 =  \frac{1}{2} (\hx_k - x_*)^T Q (\hx_k - x_*) = f(\hx_k) - f(x_*)$$
%where the last equality follows from considering the Taylor expansion of the function $f$ around its minimum $x_*$. 
%Since we are interested in the suboptimality at iteration $k$, we will set the output $\hat{z}_k$ as 
%$$ \hat{z}_k = \sqrt{f(\hat{x}_k) - f(x_*)}$$
%so that the $\ell_2$ norm of the output sequence $\hat{z}_k$ will quantify the suboptimality accumulated over the iterations.
%If we choose the output to be a linear function of the state, i.e. if $F(\xi_k) = T\tilde{\xi}_k$ for some matrix $T$, then these iterations can be rewritten as
\begin{eqnarray} 
{\xi}^c_{k+1} &=&A_Q {\xi}^c_{k} +  B  {w}_k, \quad z_k = T {\xi}^c_k,
\label{Sys:quad-deterministic-noise}
\end{eqnarray}
where $B=\tilde B\otimes I_d$ with $\tilde{B}$ as in  \eqref{def: system mat for TMM}, $T$ is defined by \eqref{def-T} and
\beq %{\hat{\xi}^c}_{k} :=\hat\xi_k - \xi_*, ~ \xi_* := \begin{bmatrix} x_*^T\\ x_*^T\end{bmatrix}, ~
A_Q:=\begin{bmatrix} 
(1+\beta)I_d -\alpha(1+\nu)Q & -\beta I_d +\alpha\nu Q \\
I_d & 0_d 
\end{bmatrix}.
\label{def-AQ}
\eeq
%Here, the superscript $c$ in $\hxic$ is to emphasize that the iterations $\hxic$ are centered, i.e. if TMM converges to the optimum, $\hxic\to 0$. 
%$\xi_* := [x_*^T x_*^T]$
For referring to this system, we will use $(A_Q, B, T)$ as a shorthand following the literature \cite[Sec. 5.3]{hinrichsen-pritchard}. In the next section, we will compute the $\ell_2$ gain corresponding to this system explicitly. %to characterize the worst-case robustness of GMM algorithms. %Also, as we discuss later in Section \cn, $H_\infty$ norm is relevant to computing the risk-sensitivity index; it will lead to a characterization of the threshold level $\bar{\theta}$ (the risk-sensitivity index will be finite only when $\theta < \bar{\theta}$).

We note that without any noise (i.e. when $w_k = 0$), it is well-known from the theory of iterative methods that the linear system \eqref{Sys:quad-deterministic-noise} converges to the optimum with arbitrary initialization if and only if $\rho(A_Q)<1$ and in this case, $\rho(A_Q)$ determines the  linear convergence rate, i.e the suboptimality decays exponentially $f(x_k) - f(x_*) \leq C_k \rho^{2k} (f(x_0) - f(x_*)) $ at rate $\rho^2=\rho^2(A_Q)$ where the pre-factor constant $C_k$ can have at most polynomial growth \cite{varga1999matrix}. In fact, for GMM methods, $C_k$ can grow at most quadratically with $k$ and a formula for the convergence rate $\rho(A_Q)$ is known \cite{can2022entropic}.\looseness=-1
 %Therefore, the $H_\infty$ norm of the linear dynamical system $(A_Q, B, T)$ given by \eqref{Sys:quad-deterministic-noise} can be interpreted as
%\beq \| G_{\hat{w},\hat{z}}\|_\infty^2 = \sup_{ \hat{d}~:~ \|\hat{d}\|_{\ell_2} <\infty, \hat d\neq 0}\frac{\sum_{k\geq 0 } f(\hx_k) - f(x_*)}{\sum_{k\geq 0} \|\hat{d}_k\|^2}.
%\label{eq-hinf-quad}
%\eeq
%starting with initialization $\hx_0 = x_*$. 
%the amount of worst-case cumulative deviation from optimality, i.e. 
%This can be interpreted as the amount of cumulative error that can be caused by injecting a worst-case perturbation to the gradients with a unit  
%whether the non-negative error signal we get $f(x_k) - f(x_*)$ is summable when we consider perturbations with a finite energy (i.e. finite $\ell_2$ norm).
\subsection{The $\ell_2$ gain of GMM for quadratic objectives.}
%Under this assumption, we will compute the $H_\infty$ norm of the GMM algorithm explicitly. 
%In general, 
When $f$ is a quadratic and GMM is globally convergent, then the dynamical system corresponding to GMM is stable and linear with system matrices $(A_Q, B,T)$. In this case, it is known that $L_{2,*}$ is equal to the \emph{$H_\infty$ norm} \cite{zhou1998essentials} 
of this system defined in the frequency domain according to the formula:
%which is defined in the frequency domain for a stable system with system matrices $(A_Q, B,T)$ as%; this is not the case for general non-linear systems.
 \beq\label{eq-h-inf-tf}% \|G\|_\infty 
H_\infty := \max_{ \omega \in [0,2\pi]} \sqrt{\lambda_{\max} \left[G(e^{i\omega})G(e^{i\omega})^H\right]} = \max_{z\in\mathbb{C}: \|z\|=1} \|G(z)\|,
	\eeq
\cite{zhou1996robust} %\footnote{With slight abuse of notation, we use $G(z)$ to denote the transfer function of the input-to-output map $G_{w, z}$ with $z$ being a complex scalar in this context.} 
where $\lambda_{\max}(\cdot)$ denotes the largest eigenvalue in this notation and 
 \beq G(z) := T ( z I - A_Q)^{-1}B, %\quad \mbox{for} \quad z\in \mathbb{C}\setminus \mbox{Spec}(A_Q),
 \label{def-transfer-matrix}
\eeq %
is called the \emph{transfer matrix} of the discrete-time system $(A_Q, B, T)$. As the name $H_\infty$ norm suggests, \eqref{eq-h-inf-tf} defines a norm in a properly defined space of transfer matrices \cite{zhou1996robust}. In addition, $H_\infty$ norms have some desirable properties such as invariance under linear transformations and are frequently used as a metric to design control systems that are robust to external disturbances \cite{hinrichsen-pritchard,zhou1998essentials}. Furthermore, $H_\infty$ norms provide
%studying robustness with respect to (perturbations) noise which can be complex-valued has also attracted a lot of interest in the robust $H_\infty$ control literature as it provides
 a rich set of information about the robustness of a system with respect to multiplicative noise \cite{hinrichsen1991stability,hinrichsen-pritchard}, which we will elaborate on in Sec. \ref{subsec-multiplicative-noise}. %In this section, we will compute the $\ell_2$ gain
%$\mbox{Spec}(A_Q) = \cup_{i=1}^d \{\lambda_i \}$ is the set of eigenvalues of $A_Q$  (see e.g. \cite{zhou1996robust}). 
%We note that $G(z)$ is well-defined for every $z\in \mathbb{C}$ except when $z$ is an eigenvalue of $A_Q$ (in which case the inverse of the matrix $zI - A_Q$ in \eqref{def-transfer-matrix} does not exist). In particular, if $\rho(A_Q)<1$, then $G(z)$ is well-defined and is continuous over the unit circle; consequently the supremum in \eqref{eq-h-inf-tf} can be replaced with a maximum.

Next, we will characterize the $\ell_2$ gain, based on calculating the $H_\infty$ norm. %through the formula \eqref{eq-h-inf-tf}. 
There are standard numerical methods for computing the $H_\infty$ norm of a linear dynamical system \cite{hinrichsen-pritchard} based on computing solutions to the discrete-time matrix Riccati equations or solving a sequence of eigenvalue problems \cite{boyd1990regularity,bruinsma1990fast,hinrichsen1991stability}, \cite[Sec. 5]{hinrichsen-pritchard}. % In particular, it is well-known that for linear systems, $H_\infty$ norm is related to the %solutions of parametrized
%certain matrix equations. More specifically, for a given parameter $\gamma>0$, consider 
%he
%discrete-time matrix Riccati equation (DARE) equation%for the system $(A_Q, B , T)$:
%\beq  X = A_Q^T X \left ( I + \frac{1}{\gamma^2} B B^T X\right) A_Q + T^T  T 
%\label{eq-disc-riccati}
%%\eeq
%\beq  X = A_Q^T X A_Q + A_Q^T X B (\gamma^2 I_{d} - B^T X B)^{-1}B^T X A_Q + T^T  T.  \label{dare}
%\eeq
%where $\gamma>0$ is a parameter. This equation admits a positive semi-definition solution $X \succeq 0 $  if and only if $H_\infty < \gamma$ \cite[Corollary 2.1, Prop. 3.7]{hinrichsen1991stability}. Therefore, a standard approach to compute the $H_\infty$ norm would be to solve the matrix Riccati equation numerically or certify that a solution does not exist for a given value of $\gamma$. Then, one can do bisection search over $\gamma$ to estimate the $H_\infty$ norm \cite{boyd1990regularity,bruinsma1990fast}. 
%The existence of a $2d\times 2d$ solution $X$ can be checked numerically by
%which require solving a $4d\times 4d$ symplectic eigenvalue problem \cite{arnold1984generalized}, 
 However, such approaches %this
%Another well-known approach would be to solve a sequence of symplectic eigenvalue problems whose solutions are related to matrix Riccati equations where each eigenvalue problem
would require typically $\mathcal{O}(d^3)$ operations per step combined with a bisection search and this can be numerically expensive for large $d$. There are also scalable approaches to large dimensions with cheaper steps such as \cite{GGO,mitchell2016hybrid,benner2014structured}, however these approaches are not guaranteed to converge globally, the convergence guarantees are of local nature. % or need to solve full eigenvalue problems to certify that an approximate solution is found. % (for example instead of $\lambda_{\max}$ in \eqref{eq-h-inf-tf}, they may potentially return another eigenvalue). Furthermore, numerical approaches would (clearly) not provide explicit solutions. %the guarantees are local under some assumptions such  %This is computationally expensive for $d$ large. %, furthermore additional bisection steps over $\gamma$ would be needed. 
%By exploiting the special structure of $A_Q, B$, and $T$ matrices, if we block diagonalize the $A_Q$ matrix with an orthogonal transformation to result in $2\times 2$ blocks, it can be shown that instead of solving the matrix equation \eqref{dare} where $X$ is $2d$ dimensional, it would suffice to solve DARE equations in dimension $2$. However, this approach would not directly lead to an explicit formula for the $H_\infty$ norm as a bisection search over $\gamma$ would still be needed.
 Instead of a numerical approach, we take an alternative path and use the formula \eqref{eq-h-inf-tf} directly to show that an explicit formula % which will lead to an explicit formula for the $H_\infty$ norm. More specifically, we consider the following frequency domain characterization of the $H_\infty$ norm:
% In the next result, we provide an explicit formula
for the $H_\infty$ norm of GMM can be obtained due to the special structure of the GMM system. %which equals $L_{2,*}$.
 Furthermore, we can characterize the parameters that achieve the best possible robustness level. %To our knowledge, this result is the first explicit characterization of the $L_{2,*}$ (and $H_\infty$ norm) of a momentum-based optimization algorithm in the literature, shedding further light into the robustness of GMM algorithms subject to worst-case square-summable deterministic (errors) noise. 
 The proof is based on the representation \eqref{eq-h-inf-tf} where we show that 
the transfer matrix $G(z)$ of GMM admits a special structure as the product of a diagonal matrix with an orthonormal matrix. Our proof also exploits certain quasi-convexity properties of the optimization objective in \eqref{eq-h-inf-tf} to show that $H_\infty$ depends only on the smallest and largest eigenvalues of $Q$ (which are $\mu$ and $L$) but not on the interior eigenvalues of $Q$.\looseness=-1%Note that in the special case %The details can be found in the appendix.
\begin{theo}\label{thm-h-inf}
Assume that the parameters $\alpha>0, \beta\geq 0, \nu\geq 0$ are such that $\rho(A_Q)<1$ and $f$ is a quadratic function of the form \eqref{eq-quad}. The worst-case robustness of the GMM algorithm %measured in terms of the $H_\infty$ norm of the corresponding dynamical system \eqref{Sys:quad-deterministic-noise} 
is
\begin{small}
	\begin{equation}\label{hinfty-quad}
	%\| G\|_\infty  
	L_{2,*} = H_\infty =  \frac{\alpha}{\sqrt{2}} \max_{\lambda \in \{ \mu, L\} }  \frac{ \sqrt{\lambda}}{r_\lambda},
\end{equation}
\end{small}
where 
\begin{equation} \label{def-r-lambda}
\begin{normalsize}
r_\lambda: = \begin{cases} | 1-c_\lambda | \sqrt{1 - \frac{b_\lambda^2}{4c_\lambda}} & \mbox{if   } c_\lambda>0 \mbox{ and } \frac{|b_\lambda| (1+c_\lambda)}{4c_\lambda} < 1,\nonumber \\
	\big| |1+c_\lambda| - |b_\lambda| \big| & \mbox{otherwise},
	\end{cases}
\end{normalsize}
\end{equation}	 
with 	
%\beq 
$b_{\lambda} :=\alpha\lambda (1+\nu)-(1+\beta)$ and $c_\lambda := \beta - \alpha \lambda \nu$.
%\label{def-b-c-lambda}
%\eeq
Furthermore, we have %the following lower bound %on the worst-case robustness level
$  L_{2,*} = H_\infty \geq \frac{1}{\sqrt{2\mu}}$  and the best robustness level $L_{2,*} = H_\infty = \frac{1}{\sqrt{2\mu}}$ is achieved if and only if $(\alpha,\beta,\nu) \in\mathcal{S}_1 \cap \mathcal{S}_2 $ where $$\mathcal{S}_1 := \big\{ (\alpha,\beta,\nu) : c_\mu \leq 0 \mbox{ or } |b_\mu (1+c_\mu)| \geq 4|c_\mu | \big\}, ~\mathcal{S}_2 := \big\{ (\alpha,\beta,\nu) : \alpha \sqrt{L \mu } \leq r_L\big\}.$$
%\mg{As a consequence, we have $\sum_{k\geq 0}\|z_k\|^2 \leq \gamma^2 \sum_{k\geq} \|w_k\|^2$ }
\end{theo} 
\begin{proof} 
%\proof {\textbf{Proof}.} 
The proof is given in Appendix \ref{app-proof-thm-hinf}. %\myqed
%\endproof
\end{proof}
In light of Thm. \ref{thm-h-inf}, we next discuss and compare the robustness of existing methods including GD, NAG, HB and TMM as a function of their parameters, starting with the following corollary. %for GD.
\begin{coro}[The $\ell_2$ gain of GD]\label{coro-robust-gd} Consider the setting of Theorem \ref{thm-h-inf}. For GD, we have $\beta = \nu = 0$ which results in $c_{\lambda} = 0$ and $b_{\lambda} = \alpha \lambda - 1$. Then, %the worst-case robustness becomes
	\beq  L_{2,*} = H_\infty = \begin{cases}  \frac{1}{\sqrt{2\mu}} & \mbox{if} \quad 0<\alpha\leq \frac{2}{L + \sqrt{L\mu}}, \\
\frac{\alpha\sqrt{L}}{\sqrt{2}(2-\alpha L)} & \mbox{if} \quad \frac{2}{L + \sqrt{L\mu}} < \alpha < \frac{2}{L}.
	\end{cases}
	\label{eq-hinf-gd}
	\eeq
	In particular, $L_{2,*} = H_\infty\to \infty$ as $\alpha \to \frac{2}{L}$.
\end{coro}

\begin{proof}This is a direct consequence of the formula \eqref{hinfty-quad}. %\myqed
\end{proof}
	\begin{figure}[ht!]
  \centering
  \hspace{-8pt}
    \begin{minipage}[t]{0.45\textwidth}
    \includegraphics[width=1\textwidth]{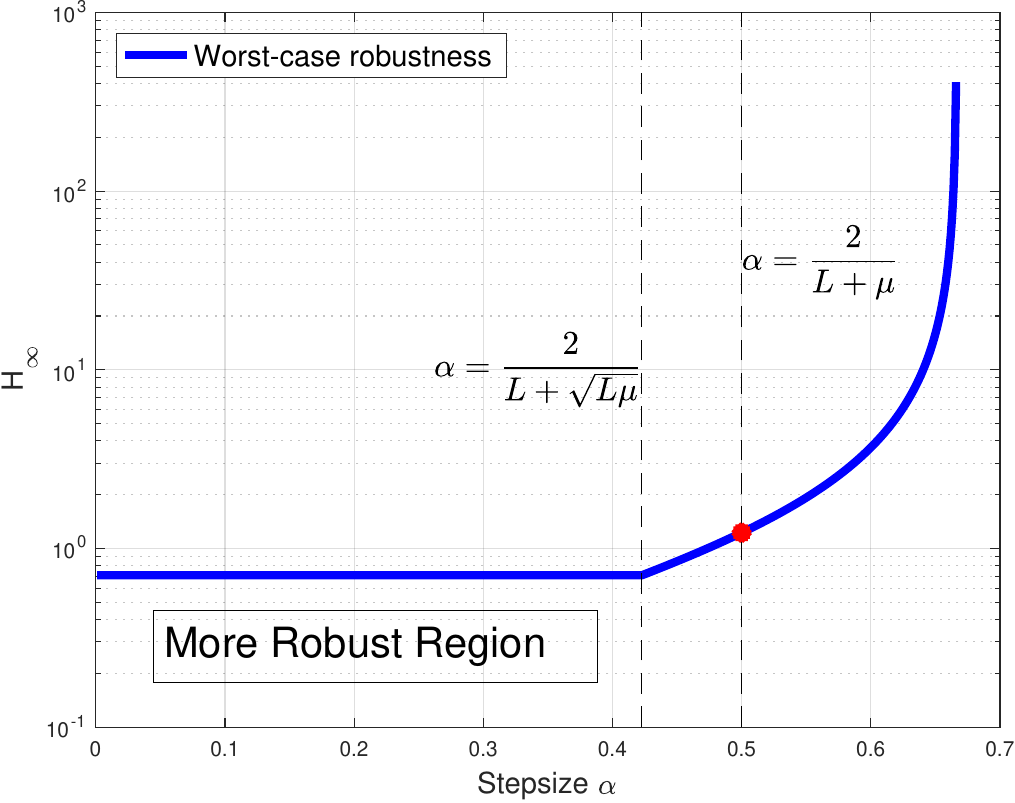}
  \end{minipage}
  \hfill
  \begin{minipage}[t]{0.45\textwidth}
  \hspace{-8pt}
       \includegraphics[width=1\textwidth]{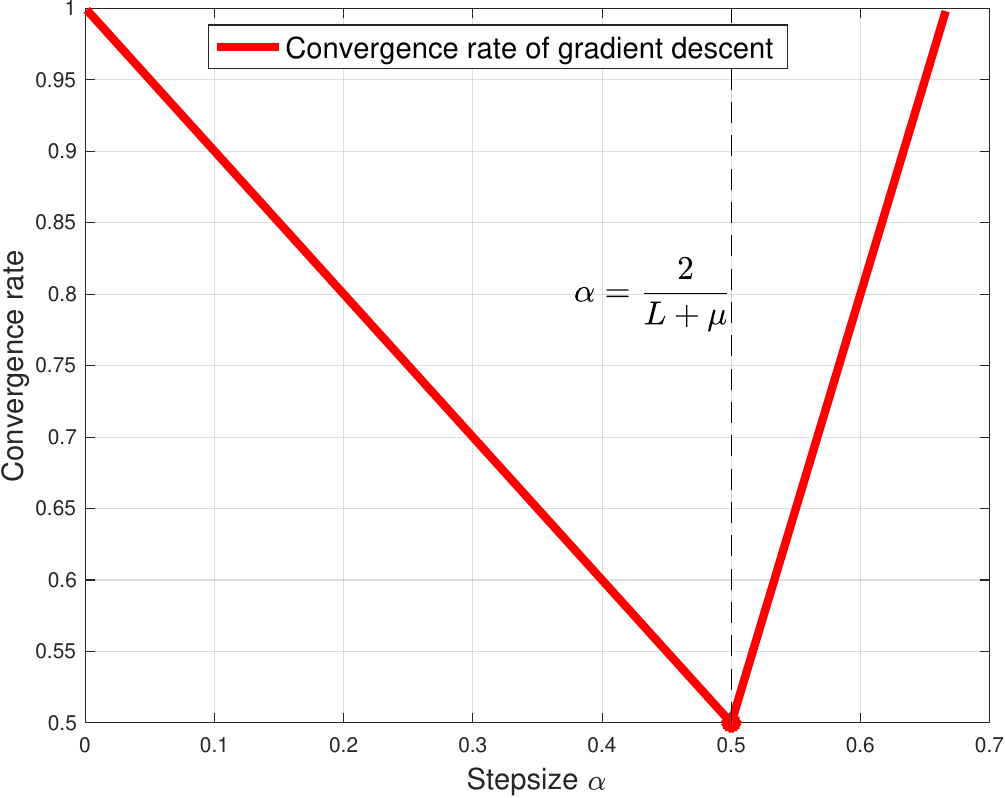}
  \end{minipage}
  \caption{\label{fig-gd} Robustness and convergence rate of GD as a function of stepsize.}
\end{figure}

\textbf{Robustness of GD.} Cor. \ref{coro-robust-gd} is illustrated in Fig. \ref{fig-gd} where on the left panel, we plot the worst-case robustness as a function of the stepsize, wheras on the right panel of Fig. \ref{fig-gd}, we plot the convergence rate $\rho_{GD} = \rho(A_Q)=\max( |1-\alpha\mu |, | 1- \alpha L |)$ of GD as a function of the stepsize for $\mu=1$, $L=3$. %where $\beta=\nu=0$.
 We observe that the best robustness is obtained when stepsize $\alpha \leq \alpha_c := \frac{2}{L + \sqrt{L\mu}}$ in which case $L_{2,*} =H_\infty = \frac{1}{\sqrt{2\mu}}$ is the smallest possible for GD. In particular, taking $\alpha=\alpha_c$ leads to the fastest rate possible $\rho_{GD} = 1- \frac{2}{\kappa + \sqrt{\kappa}}$ while retaining the best robustness level  $L_{2,*}=H_\infty = \frac{1}{\sqrt{2\mu}}$. We call this method \emph{Robustly Stable Gradient Descent} (RS-GD), to distinguish it from the other common parameter choices for GD. We also see that the common choice of stepsize $\alpha = \frac{1}{L}$ is strictly less than the critical level $\alpha_c$, leading to the best robustness level. %On the right panel of Fig. \ref{fig-gd}, we plot the convergence rate $\rho_{GD} = \rho(A_Q)=\max( |1-\alpha\mu |, | 1- \alpha L |)$ of GD where $\beta=\nu=0$. 
The red asterisk in both plots illustrates the stepsize choice $\alpha_* = \frac{2}{L+\mu}$ for which the fastest rate $\rho_* = \frac{\kappa-1}{\kappa+1}$ is attained. For this stepsize choice, the robustness becomes worse (than the best possible) at a level of $L_{2,*}=H_\infty = \frac{\sqrt{\kappa}}{\sqrt{2\mu}}$ where $\sqrt{\kappa} = \sqrt{L/\mu}>1$. These results show the trade-offs between the convergence rate and worst-case robustness. In particular, for $\alpha \leq \frac{2}{L+\mu}$, a smaller stepsize is accompanied by better robustness at the expense of a slower convergence rate. As we discuss next, similar trade-offs exist  %natural question we investigate next is whether a similar behavior exists 
for other GMM methods as well including NAG, HB and TMM.
%Then, it follows after a straightforward computation that $\|G\|_\infty $ is non-decreasing as a function of the stepsize $\alpha$. In particular, $\|G\|_\infty = \frac{1}{\sqrt{2\mu}}$ is a constant independent of the stepsize as long as $\alpha \in (0, \bar{\alpha}]$ with $\bar{\alpha} := \frac{2}{L+\sqrt{Lm}}$, whereas larger stepsizes than $\bar{\alpha}$ leads to a worst-case robustness with a higher $H_\infty$ norm given by $\| G\|_\infty =  \frac{\alpha\sqrt{L}}{\sqrt{2}(2-\alpha L)}$. In particular, we see that standard choice $\alpha_* = \frac{2}{L+\mu}$ that leads to the fastest convergence (i.e. smallest value of $\rho(A_Q)$) satisfies $\alpha_* > \bar{\alpha}$ unless $L=\mu$ and therefore is less robust with respect to a smaller stepsize $\alpha$ that satisfies $\alpha \leq \bar{\alpha}$. 

\begin{figure}[h!]
  \centering
  \begin{minipage}[t]{1.9\textwidth}
    \hspace{-11.8mm}
    \includegraphics[scale=0.36]{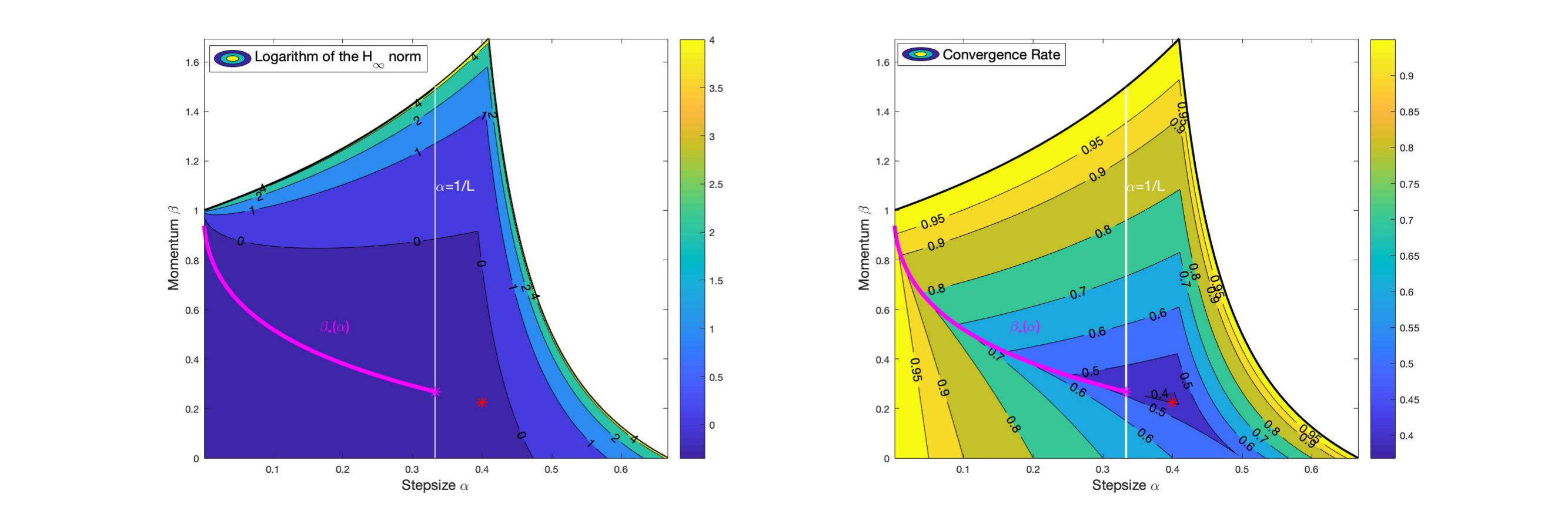}
     \vspace{-3mm}
        \caption[caption]{\label{fig-nag}(Left panel) Worst-case robustness in terms of the $H_\infty$ norm of the\\ NAG method where we report $\log(H_\infty)$ as a function of $\alpha,\beta$. (Right panel) Con-\\vergence rate of the NAG method. For both panels, we have $L=3, \mu=1$, $\nu=0$.\\ Red asterisk {\color{red}*} indicates the parameters leading to the fastest rate for NAG.}
  \end{minipage}
  \hfill
  \vspace{2mm}
  \begin{minipage}[t]{1.9\textwidth}
    \hspace{-9mm}
    \includegraphics[scale=0.39]{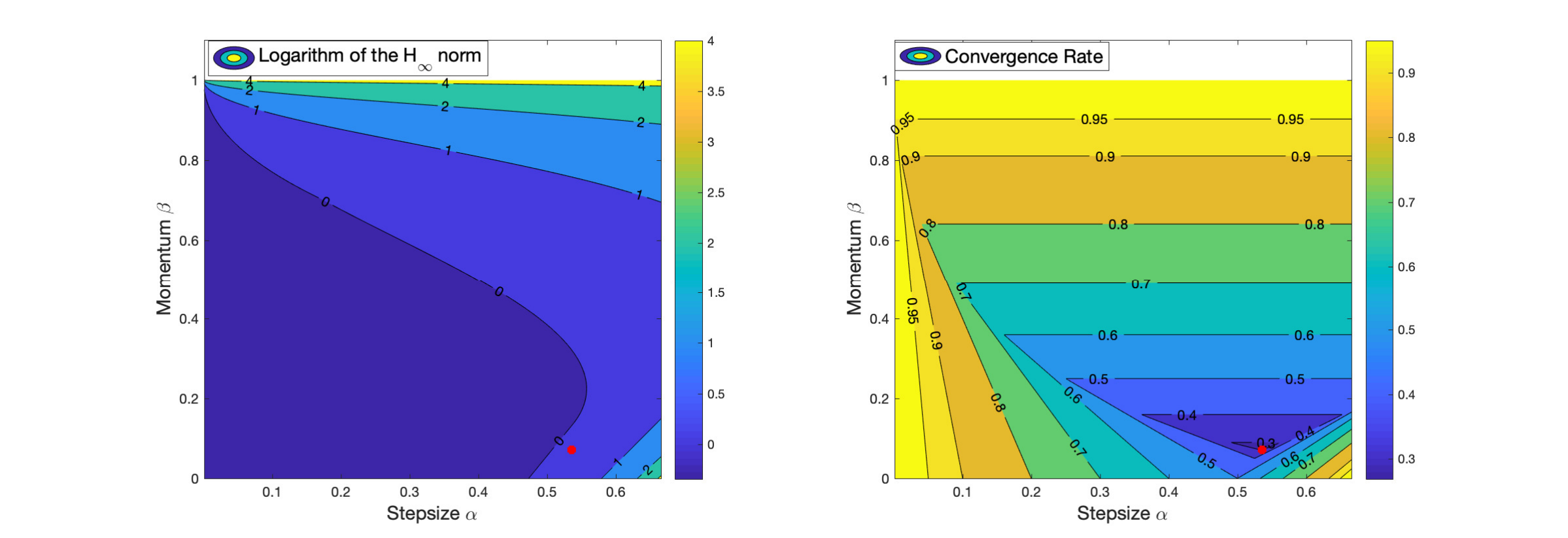}
    \hfill
    \vspace{-3mm}
    \caption[caption]{\label{fig-hb}(Left panel) Worst-case robustness of the HB method where we re-\\port $\log(H_\infty)$ as a function of $\alpha,\beta$. (Right panel) Convergence rate of the HB\\ method. For both panels, we have $L=3, \mu=1$, $\nu=1$. Red asterisk {\color{red}*} indicates \\the parameters leading to the fastest rate for HB. }
  \end{minipage}
\end{figure}

\textbf{Robustness of NAG.} On the left panel of Fig. \ref{fig-nag}, we compute the worst-case robustness of NAG based on the formula \eqref{hinfty-quad} and visualize the robustness level (as a contour plot) as the parameters are varied. %, we consider all possible parameters such that GMM is globally convergent (with $\rho(A_Q)<1$).
%as the parameters $\alpha$ and $\beta$ are varied. 
%visualize the worst-case robustness ($H_\infty$ norm) of the NAG method b
The original NAG algorithm proposed by Nesterov in the exact gradient setting takes $\alpha=1/L$ and $\beta = \nu = \frac{1-1/\sqrt{\kappa}}{1+1/\sqrt{\kappa}}$ \cite{nesterov2018lectures} but here with a slight abuse of notation, we allow any choice of parameters $\alpha>0, \beta\geq 0$ within NAG as long as the algorithm stays linearly convergent, i.e. satisfies $\rho(A_Q)<1$. On the right panel of Fig. \ref{fig-nag}, we plot the convergence rate $\rho(A_Q)$. The red asterisk indicates the choice of parameters ($\alpha = \frac{4}{3L+\mu}$, $\beta = \frac{\sqrt{3\kappa+1}-2}{\sqrt{3\kappa+1}+2})$ that leads to the fastest convergence rate (smallest $\rho(A_Q)$) for quadratic objectives. We also consider the parameter choices $\beta_*(\alpha) = \frac{1-\sqrt{\alpha \mu}}{1+\sqrt{\alpha\mu}}$ and $\alpha\in (0,1/L]$ (marked with the magenta curve) which have often been used, in fact this choice of $\beta$ leads to the fastest rate $\rho(A_Q)$ for a given $\alpha \in (0,1/L]$\cite[Lemma 2.1]{aybat2019universally}. % (this can also be observed from the right panel of Figure \ref{fig-nag})
In particular, when $\alpha = 1/L$, $\beta_*(1/L)) = \frac{1-1/\sqrt{\kappa}}{1+1/\sqrt{\kappa}}$, we recover the standard choice of parameters in the original NAG algorithm mentioned above \cite{nesterov2018lectures}. These standard
 parameters as well as the red asterisk %(corresponding to the fastest NAG rate possible) %$(\alpha,\beta) = (1/L, \beta_*(1/L))$ for NAG %are the standard choices 
%for the NAG method 
lead to an accelerated rate (i.e. the rate is $(1-\Theta(\frac{1}{\sqrt{\kappa}}))$ instead of the $1-\Theta(\frac{1}{\kappa})$ rate of GD) \cite{nesterov2018lectures,aybat2018robust} and at the same time
%Based on formula \eqref{hinfty-quad}, the red asterisk, the magenta curve (and hence the standard parameters  $\alpha = 1/L,\beta_*(1/L))$) all
lie in the (darkest blue) region which corresponds to the (best) lowest possible $H_\infty$ norm (characterized as the set $\mathcal{S}_1\cap \mathcal{S}_2$ in Thm. \ref{thm-h-inf}) and satisfy $L_{2,*} = H_\infty = \frac{1}{\sqrt{2\mu}}$. %we summarize these results in Table \ref{table}. 
These results show that the accelerated rates can be obtained for NAG at the best robustness level possible. The magenta curve also lies in the most robust region, where the rate deteriorates (gets larger) as the stepsize gets smaller. 
%The standard parameters and the red asterisk achieve acceleration with the best robustnes level possible. 
%The magenta curve corresponding to parameterized momentum $\beta_*(\alpha)$ is also in the most robust region, however the convergence rate is slower for $\alpha < 1/L$. 
Also, when $\alpha$ and $\beta$ are small enough, NAG will achieve this best robustness level. That being said, the boundary of the (most stable) darkest blue region is not a line segment but instead a curve, which means that when the parameters are close to the boundary of the darkest blue region but are outside of it, increasing the stepsize may potentially lead to an improved robustness. On the other hand, for a given fixed stepsize $\alpha>0$, when $\beta$ is larger than a threshold, robustness will get worse. These results are summarized in Table \ref{table}.\looseness=-1%where we observe that NAG can achieve an accelerated rate with standard parameters $\alpha = 1/L,\beta_*(1/L))$ while being as robust as possible. 

\textbf{Robustness of HB and TMM.} In Fig. \ref{fig-hb}, we provide analogous results for the HB method. On the left panel of Fig. \ref{fig-hb}, we plot the $H_\infty$ norm as parameters are varied, and on the right panel we plot the convergence rate $\rho(A_Q)$. We observe that for given fixed $\beta \in (0,1)$, the best possible robustness $L_{2,*} = H_\infty = \frac{1}{\sqrt{2\mu}}$ is only achievable when the stepsize $\alpha$ is small enough. In Fig. \ref{fig-hb}, the red asterisk marks the standard choice of  parameters in HB that leads to the fastest convergence rate for quadratic objectives. We found that the red asterisk is not in the region of best robustness level (marked with the darkest blue) but instead has an elevated $H_\infty$ norm at a level $L_{2,*}=H_\infty = \frac{\sqrt{\kappa}}{\sqrt{2\mu}}$. This shows that the robustness of HB can get arbitrarily worse as $\kappa$ increases.
We summarize all these findings about HB in Table \ref{table} where we also compared with the rate and robustness of TMM using the suggested stepsize from \cite{scoy-tmm-ieee}. We find that TMM's rate is slower than HB but TMM admits a relatively better robustness $H_\infty = \frac{2-\frac{1}{\sqrt{\kappa}}}{ \sqrt{2\mu}}$ than HB, although this robustness is worse than what NAG can achieve at the same rate. This begs the question of whether we can improve the robustness and stability of HB by modifying its parameters, which we discuss next. 
%We find that while %We find that 
%NAG can achieve the accelerated rate $1-\Theta(\frac{1}{\sqrt{\kappa}})$ %(which improves the standard $1-\Theta(\frac{1}{\kappa})$ rate of gradient descent for ill conditioned problems), 
%while being as robust as possible, % (with the smallest possible $L_{2,*} = H_\infty = \frac{1}{\sqrt{2\mu}}$ value). On the other hand, 
%HB with standard parameters and TMM methods are less robust (with a larger $H_\infty$ norm than $\frac{1}{\sqrt{2\mu}}$) while achieving the accelerated rates. In particular, HB is the fastest method in terms of rate and as such faster than NAG, but its robustness $H_\infty = \frac{\sqrt{\kappa}}{\sqrt{2\mu}}$ can get arbitrarily worse as the condition number of the problem increases. These results show that NAG can be both fast (admitting an accelerated rate of the form $1 - \mathcal{O}(1/\sqrt{\kappa})$) and robust, whereas HB and TMM methods with standard choice of parameters can be less robust to worst-case error. 

\textbf{Robustly stable heavy-ball (RS-HB).} There exists an explicit formula which provides the linear convergence rate $\rho=\rho(A_Q)$ of GMM methods \cite{can2022entropic, gitman2019understanding} including that of HB. Therefore, based on our explicit robustness characterization from Thm. \ref{thm-h-inf}, we can actually grid the parameter space and estimate numerically the fastest rate that can be obtained while having the best robustness level $1/\sqrt{2\mu}$ and the corresponding parameters. However, such an approach does not lead to explicit choice of parameters. In the next result, we approximate these parameters by hand, i.e., we 
propose a new modified set of parameters for HB in an explicit fashion which
%For example, if we choose
%$\alpha = a^2(\kappa)/L$ and $\beta = 1-\sqrt{\alpha \mu}$ for the HB method, we show in Lemma \ref{lemma-robust-hb} after some non-trivial calculations that the resulting method 
achieves a faster accelerated rate than NAG with standard parameters (up to a constant factor) while being at the best robustness possible. % as summarized in Table \ref{table}.
%the rate $\rho = 1 - \frac{\sqrt{2}}{\sqrt{\kappa}} + \mathcal{O}(\frac{1}{\kappa \sqrt{\kappa}})$ with $H_\infty = 1/\sqrt{2\mu}$. 
We call this method the \emph{robustly stable heavy-ball method} (RS-HB). %Basically, unlike the heavy-ball method with standard parameters, RS-HB has the best robustness level $\frac{1}{\sqrt{2\mu}}$ while admiting an accelerated rate that is faster than that of NAG for $\kappa$ large. 
To our knowledge, RS-HB method is the fastest method for quadratics (with the smallest rate $\rho$ when $\kappa \geq 32$) while having the best possible robustness level $L_{2,*}=H_\infty = \frac{1}{\sqrt{2\mu}}$. The proof relies on constraining the parameter choice to stay in the set $ \mathcal{S}_1 \cap \mathcal{S}_2$ (defined in Thm. \ref{thm-h-inf}) while improving the rate of NAG (by improving the constant that scales the $\sqrt{\kappa}$ term in the rate). The proof consists of tedious (but straightforward) computations and is included in Appendix \ref{sec-appendix-online-companion}.%\looseness=-1 %for the sake of completeness. 
\begin{prop}[RS-HB method]\label{prop-robust-hb}Consider the HB method with parameters 
$\alpha = a^2(\kappa)/L$, $\nu=0$ and $\beta = \left(1 - \frac{a(\kappa)}{\sqrt{\kappa}}\right)^2$ for minimizing a quadratic $f\in \Cml$ where $a(\kappa)=\frac{\sqrt{\kappa} (\sqrt{2\kappa -1 } -1 ) }{\kappa - 1}$ if $\kappa\geq  32$ and  $a(\kappa) =1$ otherwise,
%\beq a(\kappa) := \begin{cases}
%1 & \mbox{if} \quad \kappa < 32, \\
%\frac{\sqrt{\kappa} (\sqrt{2\kappa -1 } -1 ) }{\kappa - 1}
%& \mbox{if} \quad \kappa \geq 32,
%\end{cases} 
%\label{eq-a-kappa}
%\eeq
with $\kappa = L/\mu$. Then, this method, which we refer to as RS-HB, admits the rate $\rho = 1 - \frac{a(\kappa)}{\sqrt{\kappa}}$ and the best robustness level $L_{2,*}=H_\infty = \frac{1}{\sqrt{2\mu}}$ where $a(\kappa)$
%$ A_\kappa:= \frac{\sqrt{\frac{4}{\kappa} + 8(1- \frac{1}{\kappa}) } - \frac{2}{\sqrt{\kappa}} }{2(1- \frac{1}{\kappa})}
%$$
 is a non-decreasing function of $\kappa$ with $a(\kappa)> \frac{5}{4}$ for $\kappa \geq 32$ and $a(\kappa) = \sqrt{2} - \frac{1}{\sqrt{\kappa}} + \mathcal{O}(\frac{1}{\kappa})$.
% as $\kappa \to \infty$. 
 Therefore, $\rho = 1-\frac{\sqrt{2}}{\sqrt{\kappa}} + \mathcal{O}(\frac{1}{\kappa})$. %as $\kappa\to\infty$.
\end{prop}

\begin{figure}[h!]
  \centering
    \includegraphics[height=2.4in, width=1.7in, angle=90]{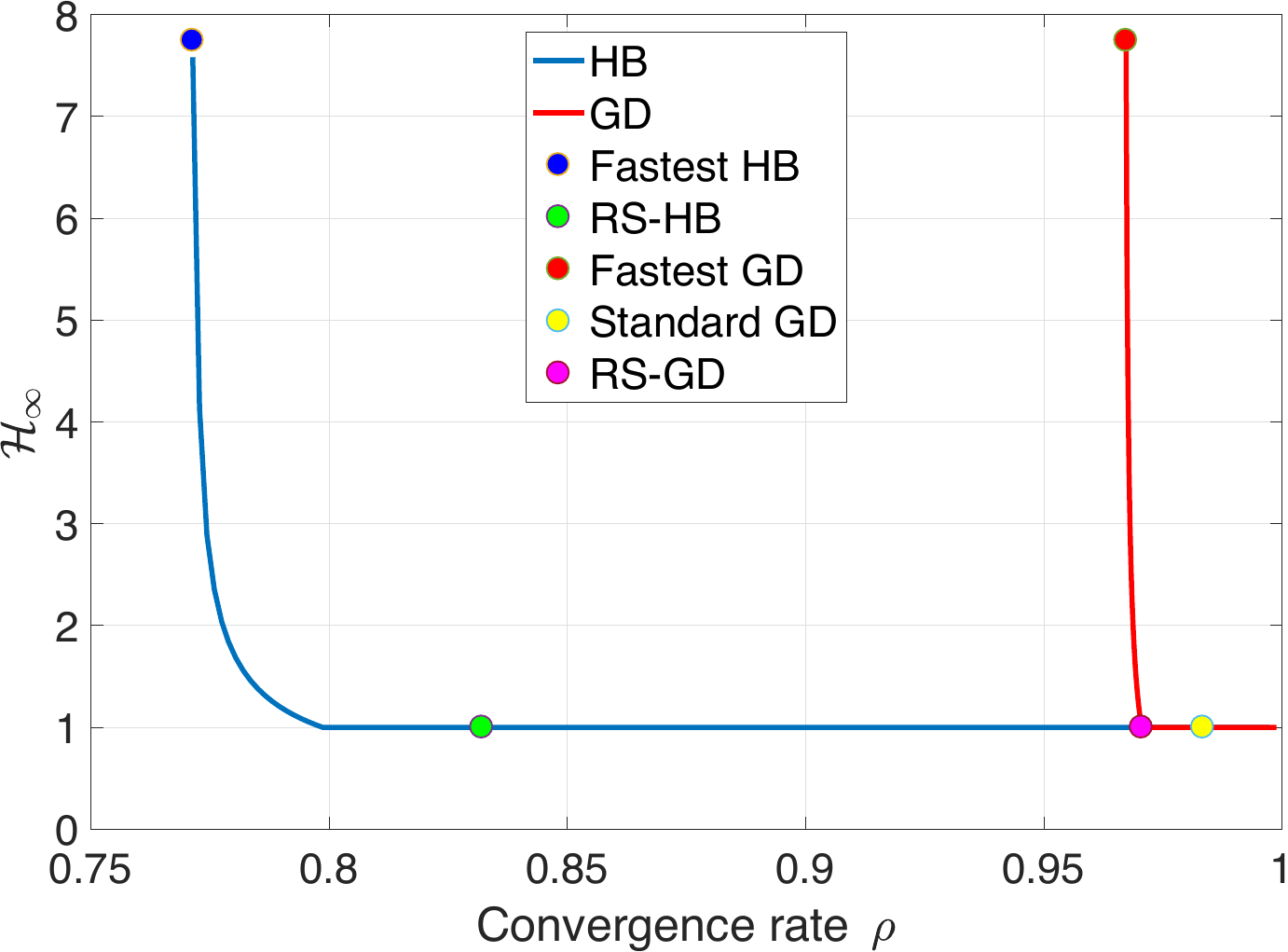}  
    \includegraphics[height=2.4in, width=1.7in, angle=90]{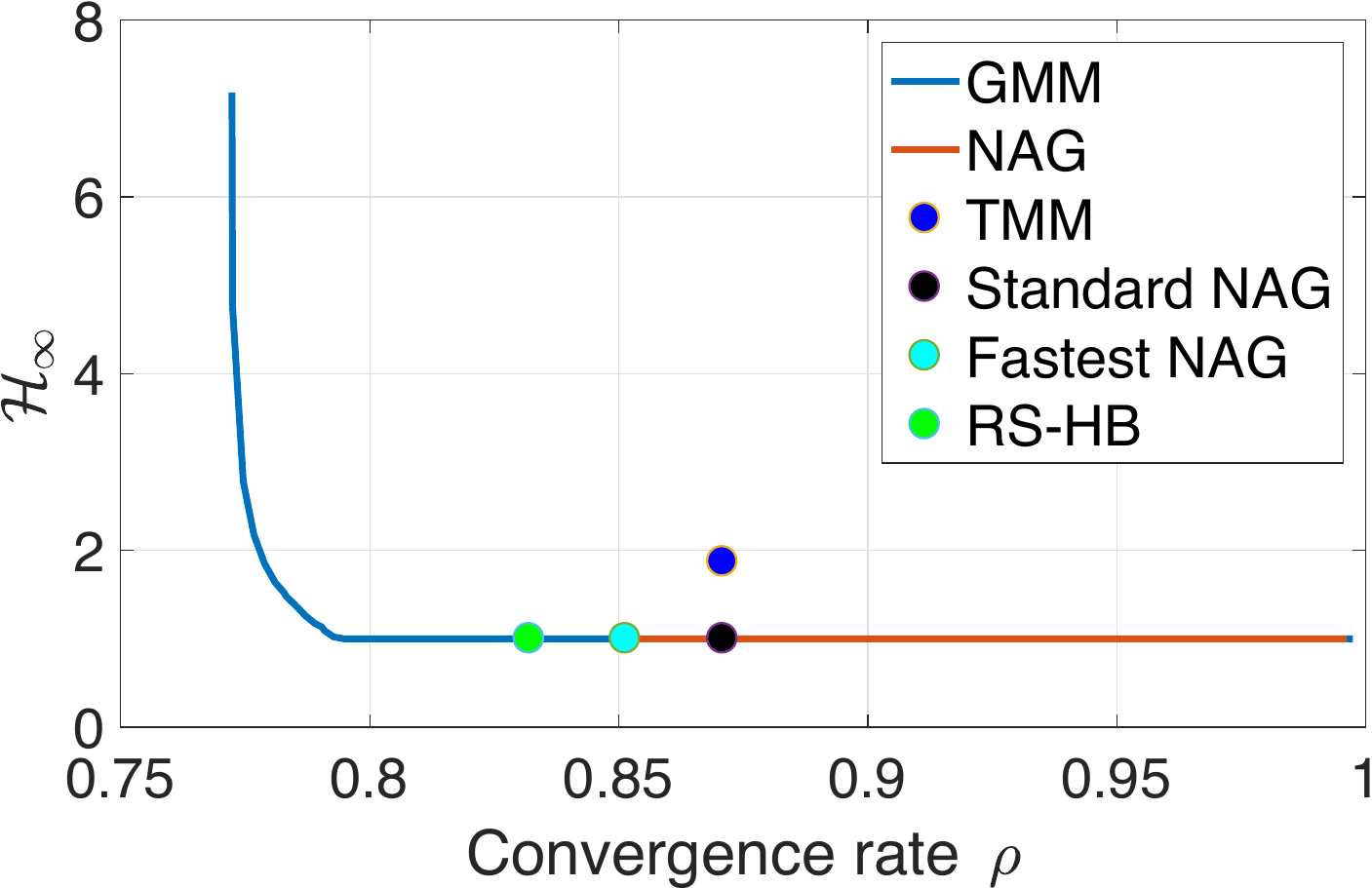}
  \caption{\label{fig-pareto-bad-cond}Pareto optimal boundary, HB, GMM, AG vs GD with standard and robust parameters $L=30, \mu=0.5$}
\end{figure}

In Fig. \ref{fig-pareto-bad-cond}, we illustrate the robustness versus convergence rate performance of the RS-HB, RS-GD, TMM, NAG and GD as well as GMM. We consider a quadratic $f$ with $\mu=\frac{1}{2}$,  $L=30$ and $\kappa = 60$ in which case the best robustness achievable is $L_{2,*} = H_\infty = 1$ according to Thm. \ref{thm-h-inf}. We display the Pareto-optimal boundary, i.e. the $x$-axis is the convergence rate and the $y$-axis displays the smallest $H_\infty$ norm possible while achieving this rate. This Pareto-optimal curve is obtained by a grid search over the parameters (except for GD where the Pareto-optimal boundary admits a closed-form expression based on Coro. \ref{coro-robust-gd}). We observe that the fastest rate in GD and the fastest rate in HB are accompanied by the largest $H_\infty$ norm ($\ell_2$ gain). %(which is $\sqrt{\kappa}/\sqrt{2\mu}$, see Table \ref{table}). 
The RS-GD and RS-HB methods we propose above achieve the best robustness level while achieving a faster rate than standard GD (with stepsize $\alpha=1/L$) and fastest NAG, lying on the Pareto optimal curve. In the right panel of Fig. \ref{fig-pareto-bad-cond}, TMM achieves an accelerated rate, but is not on the Pareto-optimal curve due to worsened robustness levels. %In Figure \ref{fig-pareto-bad-cond}, we repeat the same experiments in the same setup except that we take $L=30$ which corresponds to $\kappa=60$. 
%We note RS-GD and RS-HB have the best robustness level that is independent of $\kappa$, with a rate improvement
%The results are similar, the improvement of RS-GD and RS-HB 
%compared to standard GD (with stepsize $1/L$) and standard HB becomes more pronounced as $\kappa$ increases (in which case standard GD and HB are less robust). 
In Table \ref{table}, we summarize the $H_\infty$ norm (which equals $L_{2,*}$) corresponding to various choice of parameters for all the methods GD, HB, TMM, NAG and RS-HB. The last row of this table %relates to the other robustness metric $L_{2,*}$ and 
relates to worst-case error sequence construction and will be discussed in Section \ref{sec-worst-case-noise}. %\ref{subsec-bounds-on-l2-gain}.

{
\begin{table}[]
%\small
%\hspace{-5mm}
\centering
%\hspace{1in}
\begin{tabular}{|l|l|l|l|l| l| }
\hline
\textbf{Alg.}     & \textbf{Parameters}           & \begin{tabular}[c]{@{}l@{}}\textbf{Comments About Parameters} \end{tabular} &  \begin{tabular}[c]{@{}l@{}}\textbf{Conv.}\\ \textbf{Rate ($\rho$)} \end{tabular}                        & \begin{tabular}[c]{@{}l@{}} \textbf{\mbox{$L_{2,*}$}}\\\textbf{$(H_\infty)$} \end{tabular}  & \begin{tabular}[c]{@{}l@{}} \textbf{Info.} \end{tabular}\\ \hline
GD   & \begin{tabular}[c]{@{}l@{}} $\alpha=\frac{1}{L}$\\$\beta=\nu=0$ \end{tabular}
      & popular choice                                                    & $ 1 - \frac{1}{\kappa}$ &       $\frac{1}{\sqrt{2\mu}}$  &   \begin{tabular}[c]{@{}l@{}} $\omega_* = 0$\\$\lambda_*=\mu$ \end{tabular}
       \\ \hline
GD  & \begin{tabular}[c]{@{}l@{}}$\alpha = \frac{2}{L+\mu}$\\ $\nu=\beta=0$ \end{tabular}   & \begin{tabular}[c]{@{}l@{}}Fastest rate without noise\end{tabular} &      $ 1 - \frac{2}{\kappa+1}$                                                                                                                                               &    $\frac{\sqrt{\kappa}}{\sqrt{2\mu}}$    & \begin{tabular}[c]{@{}l@{}} $\omega_* = \pi $\\$\lambda_*=L$ \end{tabular}   \\ \hline
 \begin{tabular}[c]{@{}l@{}}RS-\\ GD \end{tabular}   & \begin{tabular}[c]{@{}l@{}}$\alpha = \frac{2}{L+\sqrt{L\mu}}$\\ $\nu=\beta=0$ \end{tabular}   & \begin{tabular}[c]{@{}l@{}}Fastest rate  while achieving \\ the best $H_\infty$ (Corollary \ref{coro-robust-gd}) \end{tabular} &      $ 1 - \frac{2}{\kappa + \sqrt{\kappa}}$                            &    $\frac{1}{\sqrt{2\mu}}$   &\begin{tabular}[c]{@{}l@{}} $\omega_* = 0$\\$\lambda_*=\mu$ \end{tabular}                                                                                                                     \\ \hline
NAG   &\begin{tabular}[c]{@{}l@{}}$\alpha=1/L, \nu=\beta$  \\ $ \beta = \frac{1-1/\sqrt{\kappa}}{1+1/\sqrt{\kappa}}$\end{tabular}     &                                                                 popular choice &     $1 - \frac{1}{\sqrt{\kappa}}$                         &                                                                 $\frac{1}{\sqrt{2\mu}}$   & \begin{tabular}[c]{@{}l@{}} $\omega_* = 0$\\$\lambda_*=\mu$ \end{tabular}                                                                                                                                                                                                                                       \\ \hline 
NAG   &   \begin{tabular}[c]{@{}l@{}} $\alpha = \frac{4}{3L+\mu},  \nu=\beta$ \\  $ \beta = \frac{\sqrt{3\kappa+1}-2}{\sqrt{3\kappa+1}+2}$\end{tabular}                         &                                                                   \begin{tabular}[c]{@{}l@{}}Fastest rate of  NAG for quadratics \\ without noise\end{tabular}  &     $ 1 - \frac{2}{\sqrt{3\kappa+1}}$ &    $\frac{1}{\sqrt{2\mu}}$   &  \begin{tabular}[c]{@{}l@{}} $\omega_* = 0$\\$\lambda_*=\mu$ \end{tabular}                                                                                                                                                                                                                                                                                                                                                                                                                                                   \\ \hline
NAG   &  \begin{tabular}[c]{@{}l@{}}$\alpha \in (0,\frac{1}{L}], \nu=\beta$  \\  $\beta = \frac{1- \sqrt{\alpha \mu}}{1+\sqrt{\alpha\mu}}$\end{tabular}                        &                                                                    \begin{tabular}[c]{@{}l@{}} Given $\alpha$ fixed, $\beta$ optimizes $\rho$\end{tabular}  &     $ 1 - \sqrt{\alpha\mu}$ &    $\frac{1}{\sqrt{2\mu}}$   & \begin{tabular}[c]{@{}l@{}} $\omega_* = 0$\\$\lambda_*=\mu$ \end{tabular}                                                                                                                                                                                                                                                                                                                                \\ \hline
TMM &  \begin{tabular}[c]{@{}l@{}}$\alpha=\frac{1+\rho}{L}, \beta = \frac{\rho^2}{2-\rho}$ \\ $\nu = \frac{\rho^2}{(1+\rho)(2-\rho)}$  \end{tabular}                          &                                                                  \begin{tabular}[c]{@{}l@{}}Proposed in \cite{scoy-tmm-ieee} \end{tabular}  &                                     $  1- \frac{1}{\sqrt{\kappa}} $                                                                                              & $\frac{2-\frac{1}{\sqrt{\kappa}}}{\sqrt{2\mu}}$  & \begin{tabular}[c]{@{}l@{}} $\omega_* = \pi $\\$\lambda_*=L$ \end{tabular}                           
   \\ \hline
\begin{tabular}[c]{@{}l@{}}HB \end{tabular} &   \begin{tabular}[c]{@{}l@{}}$\alpha=  \frac{4}{(\sqrt{L} + \sqrt{\mu})^2}$  \\ $\beta = (\frac{\sqrt{\kappa}-1}{\sqrt{\kappa}+1})^2,\nu=0$                   \end{tabular}       &                                                                  \begin{tabular}[c]{@{}l@{}}Fastest rate without noise\end{tabular}  &                                     $  1 -\frac{2}{\sqrt{\kappa}+1} $   & $\frac{\sqrt{\kappa}}{\sqrt{2\mu}}$   &\begin{tabular}[c]{@{}l@{}} $\omega_* = \pi $\\$\lambda_*=L$ \end{tabular}                                                                                                                                                                                                                                            
   \\ \hline   
   \begin{tabular}[c]{@{}l@{}}RS-\\ HB \end{tabular} &   \begin{tabular}[c]{@{}l@{}}$\alpha= \frac{a^2(\kappa)}{L}$  \\ $\beta = (1 - \frac{a(\kappa)}{\sqrt{\kappa}})^2$ \\ 
  % $a(\kappa)$ as in \eqref{eq-a-kappa}
    \end{tabular}       &                                                                  \begin{tabular}[c]{@{}l@{}}Proposed in this paper (Prop. \ref{prop-robust-hb})\end{tabular}  &                                     \begin{tabular}[c]{@{}l@{}}   $ 1 -\frac{\sqrt{2}}{\sqrt{\kappa}}$ \\ $~+ \mathcal{O}(\frac{1}{\kappa \sqrt{\kappa}})$ \\ as $\kappa\to\infty$ \end{tabular} 
   & $\frac{1}{\sqrt{2\mu}}$  &\begin{tabular}[c]{@{}l@{}} $\omega_* = 0$\\$\lambda_*=\mu$ \end{tabular}                                                                                                                                                                                                                                          
   \\ \hline
\end{tabular}
%\hspace{0.5in}
\caption{\label{table}Convergence rate, the $\ell_2$ gain ($H_\infty$ norm) and further information (about the quantities $\omega_*, \lambda_*$ arising in Prop. \ref{prop-worst-case-noise}) corresponding to different choice of GMM parameters when the objective $f\in\Cml$ is a quadratic.}
\end{table}
}
%\mtodo{Discuss figures about Pareto optimal boundary.}
%\begin{figure}[h!]
%  \centering
%    \includegraphics[height=2.4in, width=1.7in,angle=90]{hb_vs_gd_L3_muhalf.pdf}  
%    \includegraphics[height=2.4in, width=1.7in,angle=90]{gmm_vs_ag_L3_muhalf.pdf}
%  \caption{\label{fig-pareto-good-cond}Pareto-optimal boundary for HB, GMM, NAG, GD, RS-HB ($L=3, \mu=0.5$)}
%\end{figure}
%\mtodo{Change the names in the figures to RS-HB and RS-GD and to NAG to be consistent}
%\begin{figure}[h!]
%  \centering
%    \includegraphics[width=2.3in,height=1.7in]{risk_gmm_vs_ag_L3_muhalf.pdf}  
%    \includegraphics[width=2.3in,height=1.7in]{risk_hb_vs_gd_L3_muhalf.pdf}
%  \caption{Pareto optimal boundary for risk, HB, GMM, AG vs GD with standard and robust parameters $L=30, \mu=0.5$}
%\end{figure}
%\begin{figure}[h!]
%  \centering
%    \includegraphics[width=2.5in, height=1.7in]{risk_gmm_vs_ag_L30_muhalf.pdf}  
%    \includegraphics[width=2.3in, height=1.7in]{risk_hb_vs_gd_L30_muhalf.pdf}
%  \caption{Pareto optimal boundary for risk, HB, GMM, AG vs GD with standard and robust parameters $L=30, \mu=0.5$}
%\end{figure}
%\mtodo{Explain the new columns: freq. info. column in the text and the $L_{2,*}$ column.}
\subsection{The $\ell_2$ gain and worst-case multiplicative noise.}\label{subsec-multiplicative-noise}
In this section, our purpose is to show that the $\ell_2$ gain is related to worst-case multiplicative gradient noise that can be tolerated to avoid divergence. Assume that the gradient noise is a linear function of the output, i.e. ${w}_k = \Delta {z}_k$ for some matrix $\Delta \in\mathbb{C}^{d\times d}$. In this case, GMM iterations
 \eqref{Sys:quad-deterministic-noise} become
\begin{subequations}
\label{Sys:quad-deterministic-noise-2}
\begin{eqnarray} 
{\xi}^c_{k+1} &=& ( A_Q + B \Delta T) {\xi}^c_{k}, 
\end{eqnarray}
\end{subequations} 
where ${\xi}_k^c$ is as in \eqref{def-centered-iter} and $T$ is as in \eqref{def-T}. Here, we have $\nabla f(x_k) = Q({x}_k - x_*)$ where the eigenvalue decomposition $Q = U\Lambda U^T$ holds so that ${w}_k = \Delta z_k = \Delta T {\xi}_k^c = \frac{1}{\sqrt{2}} \Delta \Lambda^{1/2}U^T ({x}_k - x_*) $ and therefore,
\beq\|{w}_k \| &\leq& \frac{1}{\sqrt{2}}\| \Delta\Lambda^{-1/2}U^T\|  \| U \Lambda U^T (x_k - x_*)\| 
=\frac{1}{\sqrt{2}} \|\Delta \Lambda^{-1/2}\|  \|\nabla f(x_k)\| \nonumber \\ &\leq& \frac{\|\Delta\|}{\sqrt{2\mu}} \|\nabla f(x_k)\|, \label{ineq-multiplicative-noise}
\eeq
%$$\|\hat{w}_k \| \leq \frac{\|\Delta\|}{\mu} \|\nabla f(x_k)\|$$
where we used $U^T U = I$ and $\Lambda \succeq \mu I$. % $Q = U\Lambda U^T$ and
Here, we observe that the size of the gradient error ${w}_k$ is controlled by the size of the gradients; this setting is known as the \emph{multiplicative noise} or the \emph{relative gradient noise} when the relative error $p:= \frac{\|\Delta\|}{\sqrt{2\mu}} < 1$ \cite{polyakintroduction}. %In other words, instead of observing the actual gradient, GMM has access to the noisy gradient $\tilde \nabla f(x_k) = \nabla f(x_k) + w_k =  (Q + \Delta) (x_k-x_*)$ where the term $\Delta\in\mathbb{C}^{n\times n}$ is the source of the gradient error where the noise satisfies \eqref{ineq-multiplicative-noise}. %In particular, if $\|\Delta \| < \frac{1}{\sqrt{2\mu}}$ then
%$\frac{\|\Delta\|}{\sqrt{2\mu}} < 1$ and 
From the standard theory of iterative methods, the (asymptotic) linear convergence rate of \eqref{Sys:quad-deterministic-noise-2} will be determined by the quantity $\rho(A_Q + B\Delta T)$. In particular, for any given $\varepsilon \geq 0$, we can introduce 
\begin{equation}\rho_\varepsilon^{\mathbb{C}}(A_Q, B, T) := \max_{\Delta \in \mathbb{C}^{d\times d}: \|\Delta\| \leq \varepsilon} \rho(A_Q + B\Delta T)\label{def-pseudospec-radius}
\end{equation}
which is the worst-case spectral radius (convergence rate) over all possible (potentially complex) choices of the $\Delta$ matrix such that $\|\Delta\|\leq \varepsilon$. In fact, this quantity is known as the \emph{spectral value set radius} of $A_Q$ \cite{GGO}. We can interpret $\rho_\varepsilon(A_Q, B, T)$ as the worst-case (asymptotic) linear convergence rate of the iterations \eqref{Sys:quad-deterministic-noise-2} under the constraint $\|\Delta\|\leq \varepsilon$. In particular, for $\varepsilon$ small enough and $\rho(A_Q)<1$, we have $\rho_\varepsilon(A_Q, B, T)<1$. In this case, GMM is globally linearly convergent. However, it can be seen that when the noise level $\varepsilon>0$ is large enough, GMM diverges for some choices of $\|\Delta\|=\varepsilon$, i.e.
$\rho_\varepsilon^{\mathbb{C}}(A_Q, B, T) \geq 1$. It is known that $H_\infty$ norm (which equals $L_{2,*}$ for a linear system) is related to the smallest such $\varepsilon>0$, i.e. %admits the following characterization
%	\beq H_\infty = \frac{1}{\varepsilon_*},
	%\label{eq-hinfty}
%	\eeq
%where 
\beq \quad L_{2,*}=H_\infty = \frac{1}{\varepsilon_*}
\quad \mbox{where} \quad 
\varepsilon_* =r_{\mathbb{C}}(A_Q, B, T) := \inf \{ \varepsilon> 0 : \rho_\varepsilon^{\mathbb{C}}(A_Q, B, T) \geq 1\},
\label{def-complex-stab-rad}
\eeq
is the smallest $\varepsilon$ such that $\rho_\varepsilon^{\mathbb{C}}(A_Q, B, T) \geq1$ (see e.g. \cite{GGO} and the references therein). %, i.e. it is a solution of the equation
%$$\rho_{\varepsilon_*} = 1.$$ (see e.g. \cn). 
In other words, $L_{2,*}$ is given by the multiplicative inverse of the norm of the smallest perturbation $\varepsilon_*$ such that the eigenvalues of the perturbed matrix $A_Q + B\Delta T$ intersect the unit circle. This quantity  $\varepsilon_*$ is known as the \emph{complex stability radius} \cite{hinrichsen-pritchard}. 

If we minimize the $\ell_2$ gain, it will be harder to destabilize the matrix $A_Q$. The reason is, choosing the parameters of GMM algorithm to achieve a smaller $L_{2,*}$ (or equivalently a smaller $H_\infty$ norm) would imply that it would take a larger perturbation matrix $\Delta$ to destabilize GMM, i.e. to make GMM iterations no longer globally convergent. This shows that the $\ell_2$ gain can be also viewed as a robustness measure with respect to multiplicative gradient noise of the form $w_k = \Delta z_k$.
\begin{figure}[ht!]
  \centering
    \includegraphics[width=0.48\linewidth, height=0.45\linewidth]{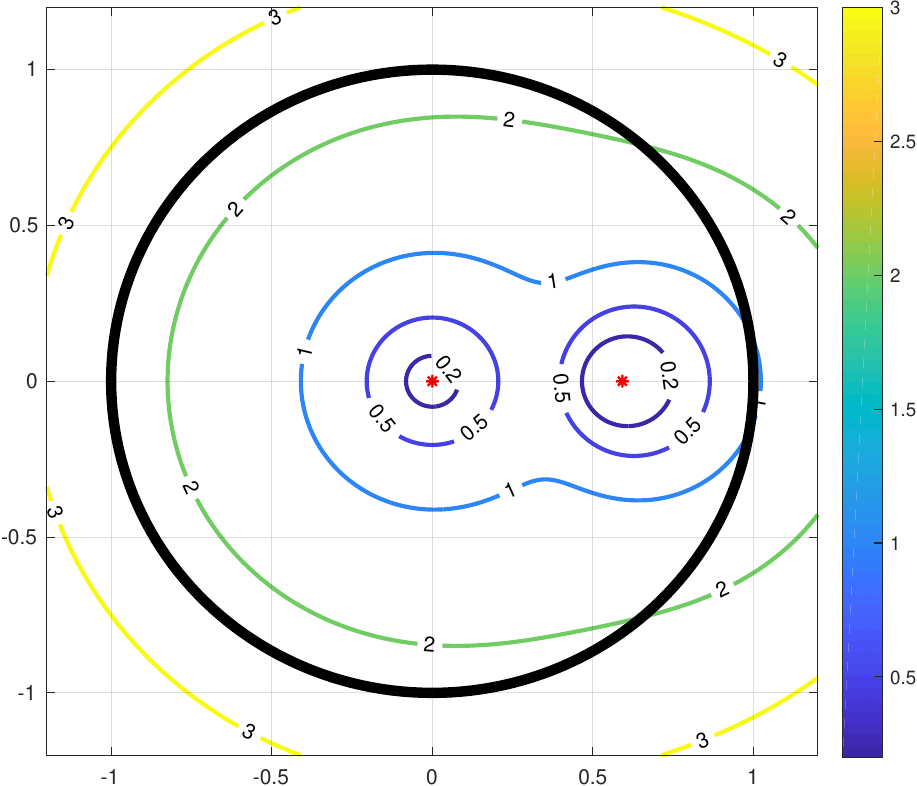}
      \includegraphics[width=0.44\linewidth, height=0.45\linewidth]{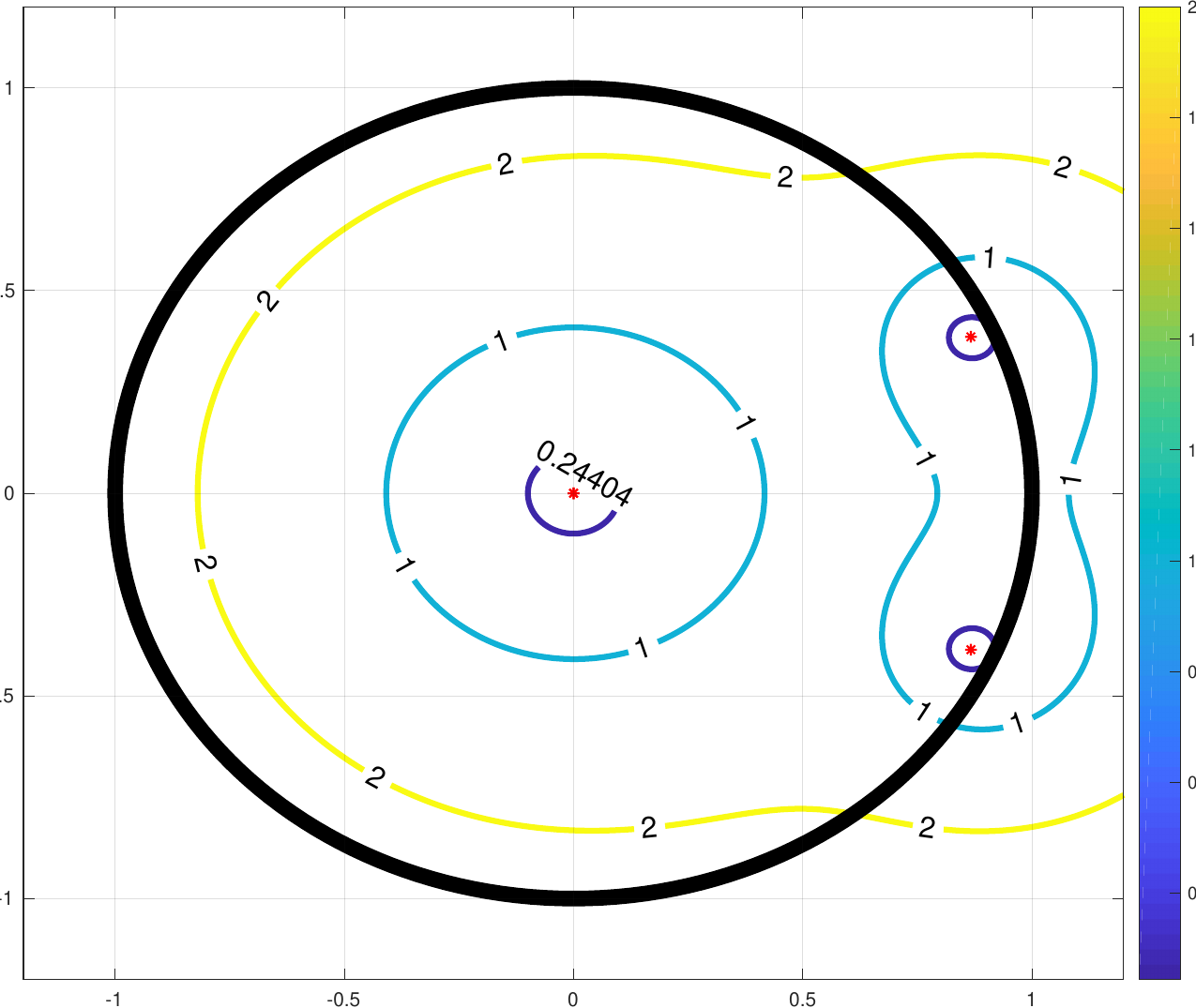}%{nag_spectral_sets_L3_muhalf_highbeta_v2-cropped.pdf}
        \caption[caption]{\label{fig-spec-val-set-nag} Boundary of the spectral value sets of NAG for $L=3$, $\mu=1/2$ and $\alpha = 1/L$. On the left panel, we take $\beta = \nu = \frac{1-\sqrt{\alpha\mu}}{1+\sqrt{\alpha\mu}}\approx 0.4202$ whereas on the right panel, we take $\beta = \nu = \frac{0.9}{1-\alpha \mu} \approx 1.08$.}
\end{figure}
 To illustrate this, in Fig. \ref{fig-spec-val-set-nag}, we plot the boundary of the sets 
$ \Lambda_\varepsilon:= \cup_{\Delta \in C^{n\times n} : \|\Delta\| \leq \varepsilon}\mbox{Spec} (A_Q + B\Delta T)$
for different values of $\varepsilon$ (tagging the boundary curves with the corresponding $\varepsilon$ value in the plots) as $\varepsilon$ is increased where $\mbox{Spec}(\cdot)$ denotes the spectrum (the set of eigenvalues) of a matrix. The sets $\Lambda_\varepsilon$ are called spectral value sets which reduce to the spectrum of $A_Q$ as $\varepsilon \to 0$ by the continuity of these sets %, we refer to
 \cite{trefethen1999spectra,hinrichsen1993spectral}. %for basic properties of spectral value sets. 
Basically, these sets relate to the worst-case asymptotic convergence rate of \eqref{ineq-multiplicative-noise}, noting that $\rho_\varepsilon^{\mathbb{C}}(A_Q, B, T) $ is defined as the largest modulus of the points lying in the set $\Lambda_\varepsilon$. %\cite{trefethen1999spectra,hinrichsen1993spectral}.
 We also display the unit circle in black. The left panel of Fig. \ref{fig-spec-val-set-nag} is for the NAG algorithm with $L=3, \mu = \frac{1}{2}$ with standard choice of parameters $\beta = \nu = \frac{1-\sqrt{\alpha\mu}}{1+\sqrt{\alpha\mu}}\approx 0.4202$ and $\alpha= \frac{1}{L}=\frac{1}{3}$. The spectral value sets $\Lambda_\varepsilon$ gets larger as $\varepsilon$ increases and hits the unit circle for $\varepsilon_* = 1$ (in which case we have $\rho_{\varepsilon_*}(A_Q,B,T)=1$. By \eqref{def-complex-stab-rad} and Thm. \ref{thm-h-inf}, we have $L_{2,*}=H_\infty = 1/\varepsilon_*= 1$. %\footnote{This would also follow from the $H_\infty$ formula for NAG given in Table \ref{table}.} 
On the right panel of Fig. \ref{fig-spec-val-set-nag}, we plot the boundary of the spectral value set for different $\varepsilon$ values for the same example except that $\beta$ is chosen larger to be $\beta = \frac{0.9}{1-\alpha \mu} \approx 1.08$. The red asterisks in both left and right panels show the eigenvalues of $A_Q$. We observe that on the right panel, $\beta$ is larger and the eigenvalues of $A_Q$ get closer to the unit circle compared to the figure on the left panel. 
%when $\beta$ is increased, 
%the eigenvalues of $A_Q$ become complex and gets closer to the unit circle. 
In this case, $\varepsilon_* \approx 0.2404$ and the robustness $L_{2,*}=H_\infty \approx \frac{1}{0.2404} \approx 4.16$ is worsened. For this particular example, we see that increasing $\beta$ leads to a decreased robustness to noise. In Appendix \ref{sec-appendix-online-companion}, we also provided additional plots that illustrate the spectral value sets of HB and TMM. 

\begin{rema}[Multiplicative noise tolerance of GMM]\label{rem-mult-noise-gmm} For the best robustness level $L_{2,*}=H_* = \frac{1}{\sqrt{2\mu}}$ (when the parameters lie in the set $\mathcal{S}_1\cap \mathcal{S}_2$ given in Thm. \ref{thm-h-inf}), noting the relationship \eqref{def-complex-stab-rad}, GMM will be convergent for any $\|\Delta \| < \varepsilon_* = \frac{1}{L_{2,*}}=\sqrt{2\mu}$. In other words, GMM will be convergent for any multiplicative noise level $p = \frac{\|\Delta\|}{\sqrt{2\mu}}<1$. However, if the robustness is not at the best level (i.e. if $L_{2,*}=H_\infty > \frac{1}{\sqrt{2\mu}}$), then there is a matrix $\Delta_*\in\mathbb{C}^{d\times d}$ with $\|\Delta_*\| = \varepsilon_* = \frac{1}{L_{2,*}}< \sqrt{2\mu}$ in which case the noise $w_k = \Delta_* z_k$ will be multiplicative satisfying $\|w_k\| \leq \delta_* \|\nabla f(x_k)\|$ with $p_* = \frac{\|\Delta_*\|}{\sqrt{2\mu}} = \frac{1}{L_{2,*} \sqrt{2\mu}}  <1$ and GMM with this multiplicative noise will be divergent. In this case, our results show that GMM cannot tolerate a multiplicative noise with parameter $p_*$ or more. Note that we can explicitly construct the noise matrix $\Delta_*$: It suffices to choose $\Delta_* = \varepsilon_* u_* v_*^H$ where $u_*$ and $v_*$ are the right and left singular vectors of the transfer matrix $G(e^{i\omega_*})$ corresponding to the largest singular value $\varepsilon_*$ and $\omega_*$ is a maximizer of $\|G(e^{i\omega})\|$ \cite{GGO,hinrichsen-pritchard}. We provide an explicit formula for $u_*, v_*$ in \eqref{def-ustar-vstar} of Appendix \ref{app-proof-of-worst-case-noise}. While this construction leads to an optimal matrix $\Delta_*$ that can be complex and noise vector $w_k$, we will show in Thm. \ref{thm-real-hinf} and Coro. \ref{coro-real-hinf-equals-complex-sometimes} that $\Delta_*$ and noise vector $w_k$ can be chosen to be real-valued in many interesting cases.% of practical interest. 
\end{rema}
\subsection{\label{sec-worst-case-noise}Construction of an almost worst-case noise sequence.} 
By the definition \eqref{def-optimal-L2-gain} of the $\ell_2$ gain as an infimum, with initialization $x_0 = x_{-1}=x_*$ in which case $H(\xi_0)=0$, given $h> 0$, there exists an input noise sequence $w^{(h)}:=\{w_k^{(h)}\}_{k\geq 0}$  for which the system output sequence $z^{(h)}:=\{z_k^{(h)} \}_{k\geq 0}$ has an $\ell_2$ gain $\gamma^{(h)}$ in the sense that
 % \beq
  % ${ (\sum_{k\geq 0}\|z_k^{(h)}\|^2)^{1/2} } / {(\sum_{k\geq 0}\|{w}_k^{(h)}\|^2)^{1/2} } = H^{(h)} \quad \mbox{with} \quad H^{(h)} \to H_\infty \mbox{ as } h \to 0$,
  $\gamma^{(h)} = \frac{\|z^{(h)}\|_{\ell_2}}{ \|w^{(h)}\|_{\ell_2}} \to L_{2,*}\mbox{ as } h \to 0$.
%  \label{def-almost-worst-case-noise}
 % \eeq 
%where $z_k^{(h)}$ is the output of the GMM system \eqref{Sys:quad-deterministic-noise} with input noise $w_k = w_k^{(h)}$ and 
 Therefore, for $h>0$ small, we can view $w^{(h)}$ as an ``almost worst-case sequence", noting that it leads to an approximately worst-case cumulative suboptimality, proportional to $(\gamma^{(h)})^2$. 
%without loss of generality we could take $\|w^{(h)}\|^2_{\ell_2(\mathbb{C}^d)} =\sum_{k\geq 0}\|{w}_k^{(h)}\|^2 = 1. $ In this case, if $x_k^{(h)}$ are the corresponding iterations, we have 
%$$ \frac{ (\sum_{k\geq 0}\|z_k^{(h)}\|^2)}{(\sum_{k\geq 0}\|{w}_k^{(h)}\|^2)} = \sum_{k\geq 0} f(x_k^{(h)}) - f(x_*) =  \left(H^{(h)}\right)^2 \to H_\infty^2.$$ 
A natural question is whether we can construct such an almost worst-case noise sequence  $\{w_k^{(h)}\}_{k\geq 0}$.\footnote{Note that for quadratic $f$, GMM system is linear and the $\ell_2$ gain $\gamma^{(h)}$ is invariant if the sequence $w^{(h)}$ is multiplied by a scalar due to linearity.  Hence, % (as this will also scale $z^{(h)}$ in the same way) and therefore 
the choice of $w^{(h)}$ achieving the performance $\gamma^{(h)}$ is not unique.}  In the following, we construct such a sequence based on a frequency domain analysis using the equivalence of the $\ell_2$ gain with the $H_\infty$ norm for linear systems. The proof constructs a sequence $\{w_k^{(h)}\}$ whose Fourier transforms admit the limit of Dirac's delta function (in the sense of the limits of generalized functions \cite{lighthill1958introduction}) as $h\to 0$ where the delta function is localized around a particular %worst-case 
frequency $\omega_*$ (for which $z_*=e^{i\omega_*}$ is a maximizer of the norm of the transfer function $\|G(z)\|$ on the unit circle). %previously calculated at the proof of Theorem \ref{thm-h-inf}. 
It turns out that inputs of the form with this frequency are amplified the most. The proof is given in Appendix \ref{app-proof-of-worst-case-noise}. \looseness=-1

\begin{prop}[Almost worst-case gradient errors]\label{prop-worst-case-noise} Let $f\in\Cml$ be a quadratic function of the form \eqref{eq-quad}. Consider GMM iterations given in \eqref{Sys: RBMM} with parameters such that $\rho(A_Q)<1$ and with initialization $x_0 = x_{-1} = x_*$. For given $h \in (0,1)$ and $k\geq 0$, we consider the noise sequence
\beq w_k^{(h)} =\sqrt{h (2-h)} (1-h)^k \cos(\omega_* k) u_* \in \mathbb{R}^d \quad \mbox{where} \quad u_* = \begin{cases} u_1 & \mbox{if } \lambda_* = \mu \\
                      u_d & \mbox{if } \lambda_* = L     
\end{cases},
\label{eq-worst-case-noise}
\eeq 
%with %\mtodo{take cos here, see Megretski and Zhou et al. (Essentials of Control book page 51), both $\omega_*$ and $-\omega_*$ are optimizers. these books say that $\ell_2$ gain is the same with $H_\infty$ norm }
 %with initialization $x_0 = x_{-1} = $ where $\begin{bmatrix} v^H & v^H \end{bmatrix}^H$ is a right eigenvector of $A_Q + B \Delta_* T$ corresponding to the eigenvalue $e^{i\omega_*}$ with $\Delta_*  = \frac{1}{H_\infty}u_* v_*^H$, 
%{\begin{small}$$
%u_* = \begin{cases} u_1 & \mbox{if } \lambda_* = \mu \\
%                      u_d & \mbox{if } \lambda_* = L     
%\end{cases}, ~
%$$ 
%\end{small}}
%\vspace{-0.2in}
{
\begin{small}
\begin{equation} \omega_*= \begin{cases} 
\mbox{arccos}(-\frac{b_{\lambda_*} (1+c_{\lambda_*})}{4 c_{\lambda_*}})  & \mbox{if}\quad c_{\lambda_*}>0 \mbox{ and } {|b_{\lambda_*} (1+c_{\lambda_*})|} < 4|c_{\lambda_*}| ,\nonumber \\
\pi & \mbox{if} \quad c_{\lambda_*} \leq 0 \mbox{ and }  (1-b_{\lambda_*} + c_{\lambda_*})^2 \leq (1+b_{\lambda_*} + c_{\lambda_*})^2, \\
\pi &  \mbox{if} \quad |b_{\lambda_*}(1+c_{\lambda_*})| \geq 4|c_{\lambda_*}| \mbox{ and }  (1-b_{\lambda_*} + c_{\lambda_*})^2 \leq (1+b_{\lambda_*} + c_{\lambda_*})^2, \\
0 & otherwise,
\end{cases} \label{eqn-freq-worst-case}
\end{equation}
\end{small}}where $\lambda_* \in \{\mu,L\}$ is the maximizer of the right-hand side of \eqref{hinfty-quad}, $u_1\in\mathbb{R}^d$ is a unit-norm eigenvector of the Hessian of $f$ corresponding to the eigenvalue $\mu$, and $u_d\in\mathbb{R}^d$ is a unit-norm eigenvector of the Hessian of $f$ corresponding to the eigenvalue $L$. Let $\{z_k^{(h)}\}_{k\geq0}$ be the output of GMM system with the input noise $ \{w_k^{(h)}\}_{k\geq 0}$ and $x^{(h)}:=\{x_k^{(h)}\}_{k\geq0}$ be the corresponding iterates. Then, we have $\sum_{k\geq 0}\|w_k^{(h)}\|^2<\infty $  for every $h \in  (0,1)$ and $\{w_k^{(h)}\}_{k\geq 0}$ is an almost worst-case sequence in the sense that 
	%\beq 
	$  \frac{ \sqrt{\sum_{k\geq 0}f(x_k^{(h)})-f(x_*)} }{\| w^{(h)}\|_{\ell_2}} = \frac{\| z^{(h)}\|_{\ell_2}}{\| w^{(h)}\|_{\ell_2}} \to L_{2,*}=H_\infty \mbox{ as } h \to 0$,
	%\label{eq-ratio-approx-hinf}
	%\eeq
where  $L_{2,*}$ is given by \eqref{hinfty-quad}.\looseness=-1
%For $\delta \in (0,1) $, consider the input noise sequence 
%$$ {w}_k^{(h)} =
%    \begin{cases}
%       \frac{1}{H_\infty + \delta} \frac{\sqrt{\mu}}{\sqrt{2}}    \frac{( R_{\omega_*}(\mu))^* }{\| R_{\omega}(\mu)\|}  u_1 u_1^T (x_k - x_*) & \mbox{if} \quad \lambda_* = \mu, \\
%         \frac{1}{H_\infty + \delta}\frac{\sqrt{L}}{\sqrt{2}}    \frac{( R_{\omega_*}(L))^* }{\| R_{\omega_*}(L)\|}  u_d u_d^T (x_k-x_*) & \mbox{if} \quad \lambda_* = L, 
%    \end{cases} \quad \mbox{for} \quad k\geq 0,
%    $$
%with    
%%to GMM, i.e. we take $w_k = w_k^{(h)}$ in \cn, with 
%%is a worst-case sequence attaining the supremum in \eqref{def-hinf-linear} with
%\begin{small}
%$$ R_{\omega_*}(\lambda_*):= \frac{1}{\sqrt{2}} \frac{-\alpha \sqrt{\lambda_*}}{e^{i\omega_*}-(1+\beta-\beta e^{-i\omega_*})+\alpha\lambda_i (1+\nu-\nu e^{-i\omega_*})} \quad \mbox{for} \quad \lambda_* \in \{\mu, L\},
%$$
%\end{small}
%where $H_\infty$ is given by \eqref{eq-hinf-quad}, 
\end{prop}
%\mtodo{$\delta$ is used in str. convex GMM bound proof as well  same symbol used twice.}
\begin{proof} 
%\proof  \textbf{Proof.} 
The proof is given in Appendix \ref{app-proof-of-worst-case-noise}. %\myqed
%\endproof
\end{proof}
\begin{rema}[Worst-case gradient errors depend on parameters]\label{remark-l2star-equals-Hinf}  %Prop. \ref{prop-worst-case-noise} shows that an  almost worst-case sequence can be chosen to have real entries; therefore it serves as a constructive proof of the fact that $L_{2,*}=H_\infty$. 
We observe from Prop. \ref{prop-worst-case-noise} that the worst-case noise \eqref{eq-worst-case-noise} depends on the parameters, for instance it is different for GD, HB or NAG. We will also use Prop. \ref{prop-worst-case-noise} in the numerical experiments section to compare the robustness of GD, HB, NAG methods when a worst-case noise sequence hits each algorithm. 
\end{rema}
%When $\omega_* = \pi$, $e^{i\omega_* k} = (-1)^k$ so that the noise sequence $w_k^{(h)}$ given in \eqref{eq-worst-case-noise} has an oscillatory patern whereas for $\omega_* = 0$,  $e^{i\omega_* k} = 1$ so that the same oscillatory behavior is not observed. Noting that the eigenvector $u_*$ can be chosen to have real entries, the noise sequence $w_k^{(h)}$ admits all real values for $\omega_* \in \{0,\pi\}$, but otherwise $e^{i\omega_* k}$ is complex-valued in general and  $w_k^{(h)}$ will have complex-valued entries. This result will also be useful for relating $L_{2,*}$ to $H_\infty$ as we discuss in the next section. We will also use Prop. \ref{prop-worst-case-noise} n the numerical experiments section to compare the robustness of GD, HB, NAG methods when a worst-case noise sequence hits each algorithm. %When an almost worst-case sequence $w_k^{(h)}$ can be constructed to have real entries for every $\delta$ small, then we have necessarily $L_{2,*}=H_\infty$. In the next section, we will discuss that this is indeed the case for many popular parameter choices for HB, GD, NAG. %For some functions, it will also% and for some functions in the class $\Cml$. 
%Next, we will estimate the minimal $\ell_2$ gain for linear systems, and show building on Prop. \ref{prop-worst-case-noise} that it is equal to the $H_\infty$ norm for many parameter choices and for some functions in the class $\Cml$. 
%\vspace{-0.1in}
\subsection{Estimating the real stability radius for quadratics.}\label{subsec-bounds-on-l2-gain} %In this section, we provide lower and upper bounds for the minimal $\ell_2$ gain $L_{2,*}$ for quadratic objectives that are tight in the sense that they can match for some choices of the parameters or for some particular choices of $f$ regardless of the parameters. 

Recall that from \eqref{def-complex-stab-rad}, the multiplicative inverse of the $H_\infty$ norm is equal to the complex stability radius. In the following, we define the real stability radius, which is the analogue of \eqref{def-complex-stab-rad} when the perturbation matrix $\Delta$ is restricted to be real and discuss its connections to $L_{2,*}$. The \emph{real stability radius} is defined as
\beq r_{\mathbb{R}}(A_Q, B, T) := \inf \{ \varepsilon : \rho_\varepsilon^{\mathbb{R}}(A_Q, B, T) \geq 1\},\eeq
with
\beq \quad \rho_\varepsilon^{\mathbb{R}}(A_Q, B, T) := \max_{\Delta \in \mathbb{R}^{n\times n}: \|\Delta\| \leq \varepsilon} \rho(A_Q + B\Delta T).
\label{eq-lim-perf}
\eeq 
%where  \begin{equation}\rho_\varepsilon^{\mathbb{R}}(A_Q, B, T) := \max_{\Delta \in \mathbb{R}^{n\times n}: \|\Delta\| \leq \varepsilon} \rho(A_Q + B\Delta T).
%\label{def-real-pseudospec-radius}
%\end{equation}
%$\rho_\varepsilon^{\mathbb{R}}(A_Q, B, T)\geq \rho_\varepsilon^{\mathbb{C}}(A_Q, B, T)$. 
Computing the real stability radius numerically is possible by solving a non-convex min-max optimization problem where subproblems require optimizing the second largest singular values of a $2d\times 2d$ matrix and this is typically computationally expensive in moderate to high dimensions (see \cite{qiu1995formula}). We note that by definition, $r_{\mathbb{R}}(A_Q, B, T)\geq r_{\mathbb{C}}(A_Q, B, T)$. The next result shows that for many choices of GMM parameters (including all the methods studied in Table \ref{table}), in fact the equality holds. This is relevant to optimization practice, because it implies that the $\Delta_*$ matrix and the worst-case multiplicative noise $w_k = \Delta_* z_k$ can be chosen to have all real entries (see also Remark \ref{rem-mult-noise-gmm}). For more general parameters, we also provide an upper bound on $r_{\mathbb{R}}(A_Q, B, T)$ in the following result by a frequency domain analysis approach.%where subproblems can be solved with a bisection search \cite{qiu1995formula} but this is computationally expensive unless the dimension is small. 
\begin{theo}\label{thm-real-hinf}Assume that $\rho(A_Q)<1$ and $f$ is a quadratic function of the form \eqref{eq-quad}. The worst-case robustness of the GMM algorithm measured in terms of $L_{2,*}$ of the corresponding dynamical system \eqref{Sys:quad-deterministic-noise} satisfies the following:
\begin{tightitemize}
\item [$(i)$] If $c_{\lambda_*} \leq 0   \quad \mbox{or} \quad {|b_{\lambda_*}(1 + 	c_{\lambda_*})|}\geq 4| c_{\lambda_*}|$, then 
$r_{\mathbb{R}}(A_Q, B, T)  = r_{\mathbb{C}}(A_Q, B, T) = H_\infty^{-1} = L_{2,*}^{-1}$
%$$L_{2,*} = H_\infty \quad \mbox{and} \quad 
%r_{\mathbb{R}}(A_Q, B, T) = L_{2,*}^{-1},$$
 where $H_\infty$ is given by \eqref{hinfty-quad},  $\lambda_* \in \{\mu, L\}$ is a maximizer of the optimization problem in \eqref{hinfty-quad} and $b_{\lambda_*}, c_{\lambda_*}$ are as in Thm. \ref{thm-h-inf}.%defined by \eqref{def-b-c-lambda}. 
\item [$(ii)$] Otherwise, i.e if $c_{\lambda_*}>0$ and $|b_{\lambda_*} (1 + 	c_{\lambda_*})| < 4|c_{\lambda_*}|$, we have 
%$$%\frac{\sqrt{\alpha}}{2}\min \bigg(\max_{\lambda \in \{\lambda_2, L\} } \frac{\sqrt{\lambda}}{r_\lambda}, \max_{\lambda \in \{\mu, \lambda_{d-1}\} } \frac{\sqrt{\lambda}}{r_\lambda} \bigg) 
%\frac{1}{\sqrt{2\mu}} \leq L_{2,*}^{\footnotesize\mbox{lb}}  \leq  L_{2,*} \leq H_\infty \quad \mbox{and} \quad H_\infty^{-1} \leq r_{\mathbb{R}}(A_Q, B, T) \leq (H_\infty^{\tiny\mbox{lb}} )^{-1}, $$
%$$H_\infty^{-1} = r_{\mathbb{C}}(A_Q, B, T) \leq r_{\mathbb{R}}(A_Q, B, T) \leq (H_\infty^{\tiny\mbox{lb}} )^{-1}$$
$ (H_\infty^{\tiny\mbox{lb}} )^{-1} \geq r_{\mathbb{R}}(A_Q, B, T) \geq  r_{\mathbb{C}}(A_Q, B, T) = H_\infty^{-1} = L_{2,*}^{-1}
$
where
\begin{footnotesize}
%$$ L_{2,*}^{\footnotesize\mbox{lb}} = \frac{\alpha}{\sqrt{2}} \max_{\lambda \in \{\mu, L\}}  \frac{\sqrt{\lambda}}{| |1+c_\lambda| - |b_\lambda| |},
%$$
 \beq  H_\infty^{\tiny\mbox{lb}} := \frac{\alpha}{\sqrt{2}} \max_{\omega \in [0,2\pi]}
\min \bigg( 
\max_{\lambda \in \{\lambda_2, L\}} 
\frac{ \sqrt{\lambda}}{\|e^{2i\omega} + b_{\lambda} e^{i\omega} + c_{\lambda} \|} ,  \max_{\lambda \in \{\mu, \lambda_{d-1}\}} 
\frac{\sqrt{\lambda}}{\|e^{2i\omega} + b_{\lambda} e^{i\omega} + c_{\lambda} \|}  
\bigg),
\label{def-hinfty-lb}
\eeq
 \end{footnotesize}
$\mu = \lambda_1 \leq \lambda_2 \leq \dots \leq \lambda_{d-1} \leq \lambda_d = L$ are the eigenvalues of $Q$, and $H_\infty, b_{\lambda}, c_{\lambda}$ are as in Thm. \ref{thm-h-inf}.%defined by \eqref{def-b-c-lambda}. 
\end{tightitemize} 
%\begin{small}
%	\begin{equation}
%	H_\infty^{\mathbb{R}} =  
%	\begin{cases}
%			H_\infty^{\mathbb{C}} & \mbox{if} \quad
%			c_{\lambda_*} \leq 0   \quad \mbox{or} \quad \frac{\|b_{\lambda_*}| (1 + 	c_{\lambda_*})}{4c_{\lambda_*}} \geq 1, \\
%		\end{cases}
%\end{equation}
%\end{small}
% Otherwise, in the remaining case when $c_{\lambda_*}>0$ and $\frac{|b_{\lambda_*}| (1 + 	c_{\lambda_*})}{4c_{\lambda_*}} < 1$, we have 
%with 
\end{theo}
\begin{proof}
%\proof {\textbf{Proof.}}
 The proof is given in Appendix \ref{app-proof-of-real-stab-rad}. %\myqed
%\endproof
\end{proof}
%\mtodo{Give examples of GMM parameters which lead to each of the two cases above.}
%\vspace{-0.2in}
In the last row of Table \ref{table}, we report the values of $\lambda_*$ and $\omega_*$ defined in Prop. \ref{prop-worst-case-noise} for all the algorithms and the corresponding parameter choices. It can be checked that we have either $c_{\lambda_*}\leq 0$ or  ${|b_{\lambda_*}(1 + 	c_{\lambda_*})}\geq 4| c_{\lambda_*}$ in every case for Table \ref{table}. Therefore we have the following corollary. %to part (i) of Theorem \ref{thm-real-hinf}.
\begin{coro}\label{coro-real-stab-equals-complex-stab-cases}
For any given $f\in\Cml$, for all the algorithms and for all the choice of parameters given in Table \ref{table}, it holds that  $r_{\mathbb{R}}(A_Q, B, T)  = r_{\mathbb{C}}(A_Q, B, T) = H_\infty^{-1} = L_{2,*}^{-1}.$%$L_{2,*} = H_\infty$.
\end{coro} 

%\tau_1(G(e^{i\omega})) \geq 
Part $(i)$ of Theorem \ref{thm-real-hinf} and its consequence Coro. \ref{coro-real-stab-equals-complex-stab-cases} show that $r_{\mathbb{R}}(A_Q, B, T)  = r_{\mathbb{C}}(A_Q, B, T)$ for parameter choices that satisfy certain inequalities. The next result shows that even if these inequalities are not satisfied (when we are in the setting of part $(ii)$ of Theorem \ref{thm-real-hinf}), we have also $r_{\mathbb{R}}(A_Q, B, T)  = r_{\mathbb{C}}(A_Q, B, T)=(H_\infty^{\tiny\mbox{lb}} )^{-1}$ provided that the eigenvalues $\mu$ and $L$ of the Hessian of $f$ have multiplicity two or more (this requires at least four eigenvalues and hence requires $d\geq 4$). This shows that the lower and upper bounds obtained in part (ii) of Theorem \ref{thm-real-hinf} for $r_{\mathbb{R}}(A_Q, B, T)$ are tight in the sense that the bounds match for some choices of $f$ regardless of the parameter choice as long as the GMM is globally convergent without errors, i.e. if $\rho(A_Q)<1$. 
%That being said, 
%%even if the smallest and largest eigenvalues of the Hessian of $f$ has multiplicity two or more, then $L_{2,*} = H_\infty$. 
%certifying the multiplicity of the eigenvalues $\mu$ and $L$ would require computing all the eigenvalues, which may be computationally expensive. Since our purpose is to obtain practical upper bounds to $L_{2,*}$ that  only requires the estimate of $\mu$ and $L$ without requiring the knowledge of all the eigenvalues of the Hessian or the multiplicity of the eigenvalues. The next result shows that the best upper bound to $L_{2,*}$ we hope to obtain with this approach is $H_\infty$ norm.
%The following result shows that for $d\geq 4$, in the worst case over the choice of function in the class $\Cml$, $H_\infty$ norm will be the same regardless of the fact that the noise is complex-valued or real-valued.
%When the multiplicity of the eigenvalues $\mu$ and $L$ are one, it turns out that there exists some parameter choices for which $L_{2,*}< H_\infty$ (see Appendix \cn). 
%Therefore, the metrics $L_{2,*}$ and $H_\infty$ are not totally equivalent. That being said, if the multiplicity of $\mu$ and $L$ are two or more, then we will have $L_{2,*}< H_\infty$ for every parameter choice satisfying $\rho(A_Q)<1$. 
\begin{coro}\label{coro-real-hinf-equals-complex-sometimes} Assume that $ f\in \Cml$ is a quadratic function of the form \eqref{eq-quad} with a Hessian matrix $Q$ and that $\rho(A_Q)<1$. For $d\geq 4$, if $\lambda_2 = \mu$ and $\lambda_{d-1} = L$, then 
	$r_{\mathbb{R}}(A_Q, B, T) = r_{\mathbb{C}}(A_Q, B, T) = (H_\infty^{\tiny\mbox{lb}})^{-1} =  H_{\infty}^{-1}$.	%Hence, for $d\geq 4$, t%here exists quadratic $f\in\Cml$ with a Hessian $Q$ for which $r_{\mathbb{R}}(A_Q, B, T)=r_{\mathbb{C}}(A_Q, B, T).$
%	$ \sup_{f\in \Cml} r_{\mathbb{R}}(A_Q, B, T)= \sup_{f\in \Cml} r_{\mathbb{C}}(A_Q, B, T).
%$
\end{coro}
\begin{proof}
%\proof {\textbf{Proof}.} 
%By assumption, $\lambda_2 = \mu$, and $\lambda_{d-1} = L$ holds.%\footnote{Since $\mu \neq L$, $d\geq 4$ is a necessary condition for this to hold.} 
In light of Thm. \ref{thm-real-hinf}, it suffices to show that $H_\infty^{\tiny\mbox{lb}} = H_\infty$. Plugging $\lambda_2=\mu$ and $\lambda_{d-1}=L$ into \eqref{def-hinfty-lb}, we obtain
	$ H_\infty^{\tiny\mbox{lb}}   = \frac{\alpha}{\sqrt{2}} \max_{\omega \in [0,2\pi]} \max_{\lambda \in \{\mu, L\}} 
\frac{ \sqrt{\lambda}}{\|e^{2i\omega} + b_{\lambda} e^{i\omega} + c_{\lambda} \|} = 
\frac{\alpha}{\sqrt{2}}  \max_{\lambda \in \{\mu, L\}} 
\frac{ \sqrt{\lambda}}{ r_\lambda }, 
$
by swapping the max operators. Then, from Thm. \ref{thm-h-inf}, we conclude that $H_\infty= H_\infty^{\tiny\mbox{lb}} $. This proves the desired result.\looseness=-1 %\myqed
\end{proof} 
%\endproof
%\vspace{-0.2in}
\begin{rema}\label{rem-lower-bound} %Since the class $\Cml$ includes strongly convex quadratic functions, 
Based on Coro. \ref{coro-real-hinf-equals-complex-sometimes} and Thm. \ref{thm-h-inf}, we can argue that we have %we have the lower bound %the lower bound 
\beq 
 \sup_{f \in \Cml} L_{2,*} 
\geq  
\footnotesize{ \sup_{ \begin{tabular}[c]{@{}l@{}} $f \in \Cml$,\\ $f \mbox{ is a quadratic}$ \end{tabular}} L_{2,*} } 
=
\footnotesize{ \sup_{ \begin{tabular}[c]{@{}l@{}} $f \in \Cml$,\\ $f \mbox{ is a quadratic}$ \end{tabular}} H_\infty } \geq \frac{1}{\sqrt{2\mu}}, 
\label{ineq-lower bound for l2 gain}
\eeq
so our results from Thm. \ref{thm-h-inf} can serve as lower bounds for more general $f\in\Cml$.
\end{rema}

Next, we obtain upper bounds for $L_{2,*}$ for the class $f\in\Cml$. To our knowledge, algorithms that can compute $L_{2,*}$ accurately with running time complexity that is independent of the dimension $d$ do not exist for general nonlinear systems. For GMM, we will obtain upper bounds for $L_{2,*}$ based on checking whether the parameters satisfy a $4\times 4$ matrix inequality, regardless of the dimension.\looseness=-1
%\vspace{-0.1in}
\section{Main results for strongly convex functions.}\label{sec-strongly-convex}
When the objective $f$ is strongly convex and smooth but is not a quadratic, then the dynamical system corresponding to the GMM iterations are non-linear (because the gradients are non-linear functions of the state) and the %GMM dynamics can no longer be characterized with constant system matrices and the 
$H_\infty$ representation \eqref{eq-h-inf-tf} is no longer valid. Instead, the $\ell_2$ gains for nonlinear systems are often characterized by obtaining numerical solutions to matrix inequalities (MIs). For example, 
%To our knowledge, algorithms that can compute $L_{2,*}$ exactly and efficiently (with a running time independent of the dimension) for general nonlinear do not exist. 
for some special nonlinear systems such as those with sector-bounded nonlinearities, using the distances as a Lyapunov function (i.e. taking a Lyapunov function of the form $V(\xi):=\xi_c^T P \xi_c$ for some $P\succ 0$), there exist linear matrix inequality-based approaches to compute an upper bound for the $\ell_2$ gain of a system \cite[Section 4.3]{turner2007mathematical}. However, such an approach based on the off-the-shelf solvers would require solving matrix inequalities involving $\mathcal{O}(d)\times \mathcal{O}(d)$ matrices which would be expensive in high dimensions and would not lead to explicit estimates. Secondly, 
%does not directly apply to our setting as the nonlinearity we have in the gradients has a different structure than sector-boundedness studied in \cite{turner2007mathematical}, in particular 
having a tight rate analysis for accelerated methods requires exploiting several inequalities between the gradient size and distance to the optimum obeyed specifically by strongly convex smooth functions that would not necessarily be satisfied for systems with sector-bounded nonlinearity (see e.g. \cite{aybat2018robust, gannot2022frequency,hu2017dissipativity}). %\cite{hu2017dissipativity}. 
%This is relevant, because ideally for a tight analysis we would like our robustness analysis to recover the existing convergence analysis of deterministic GMM when the noise input is set to zero. 
Third, accelerated methods are not always monotonic in distance, %\footnote{\label
%{footnote-nonmonotone}For instance, in the quadratic case without any gradient noise, with the parameter choice $\alpha \in (0,1/L]$ and $\beta = \frac{1 - \sqrt{\alpha\mu}}{1 + \sqrt{\alpha\mu}}$ and $P = I_{2d}$, the iteration matrix $A_Q $ admits a non-trivial Jordan block of size 2, therefore for some particular choices of the initialization, the Lyapunov function $V(P):= \xi_k^c P \xi_k^c $ will be proportional to $k^2\rho(A_Q)^{2k}$ (see \cite{can2022entropic,can2019accelerated}) and the latter quantity is not monotonically decreasing over $k$ for sufficiently small stepsize in which case $\rho(A_Q)$ is close to $1$. } 
so relying only on distances as a Lyapunov function does not lead to a tight rate analysis that can recover the fastest accelerated rates from the literature. To address these challenges and gain further insight into how $\ell_2$ gains can be estimated, we first consider GD and NAG algorithms and obtain explicit bounds for the $\ell_2$ gain by leveraging various inequalities for strongly convex and smooth functions. Our analysis shows that an approach based on more general Lyapunov functions that are weighted sums of distances and suboptimality yields tighter robustness bounds not only for NAG but also for GD. We then present a scalable MI approach in Sec. \ref{subsec-str-cvx-MI} for more general GMM parameters. %which, if feasible, provides a bound for the $\ell_2$ gain.%\looseness=-1 %This MI is motivated by our analysis for GD and NAG in the sense that it is designed to contain enough parameters and to be flexible enough to generalize our hand-made analysis for GD and NAG. %an MI is actually equivalent to inequalities. %MI can be thought of various inequalities %in the sense that fo
\subsection{Explicit bounds for $L_{2,*}$ under strong convexity.}\label{subsec-str-cvx-explicit}
Consider the noisy GMM updates \eqref{Sys: RBMM} for $f\in\Cml$. %We recall that the $\ell_2$ gain  $L_{2,*}$ of this nonlinear system is defined according to \eqref{def-l2-gain-general} and \eqref{def-optimal-L2-gain}. Therefore, for proving that $\ell_2$ gain $\leq \gamma$, it suffices to show that
% \beq \sum_{k\geq 0} \left[ f({x}_k) - f_* - \gamma^2 \|{w}_k\|^2 \right] \leq h_0 \left[ f(x_0) - f(x_*) \right], 
 %\label{def-hinfty-nonlinear}
 %\eeq
%starting from any initial point  $x_0 = x_{-1} \in \mathbb{R}^{d}$ where $h_0>0$ is a constant independent of $x_0$. Basically, this inequality corresponds to choosing $H(\xi_0) = h_0 \left[ f(x_0) - f(x_*) \right]$ in \eqref{def-l2-gain-general} where we recall that we take $\xi_0 = \begin{bmatrix} x_0^T & x_0^T\end{bmatrix}^T$. 
First, we will derive explicit bounds for the $\ell_2$ gain of GD. %when the objective $f$ is a strongly convex smooth function.% We start with 
%gradient descent.%, where
%\subsection{Explicit $L_{2,*}$ bounds for gradient descent}
%n the case of gradient descent, 
%iterations simplifies to
%\beq {x}_{k+1}=\Agd {x}_{k}+\Bgd \nabla f({x}_k) + \Bgd {w}_k, \quad \Agd= I_d, \quad \Bgd=-\alpha I_d.
%\label{def-GD-deter-noise}
%\eeq
%the following result which provides bounds on the worst-case robustness of gradient descent in terms of $\ell_2$ gain for strongly convex smooth objectives. 
\begin{prop}[Explicit $L_{2,*}$ bound for GD]\label{prop-hinf-gd-bound} Consider minimizing $f \in \Cml$ with GD subject to gradient errors satisfing Assump. \ref{assump-deter-noise-real} with constant stepsize $\alpha \in (0, 2/L)$. Then, %The worst-case robustness %in terms of $L_{2,*}$ 
%of GD method satisfies
$$\small L_{2,*} \leq  \begin{cases} 
	\frac{1}{\sqrt{2\mu}} & \mbox{if} \quad 0 < \alpha \leq \frac{1}{L}, \\
    \frac{1}{\sqrt{2\mu}}\frac{\alpha L}{2-\alpha L} & \mbox{if} \quad \frac{1}{L} <\alpha \leq \frac{2}{L+\sqrt{L\mu}},\\
      \frac{1}{\sqrt{2\mu}}  \sqrt{\kappa} & \mbox{if} \quad  \frac{2}{L+\sqrt{L\mu}} < \alpha \leq \frac{2}{L+\mu}, \\
\frac{\alpha\sqrt{L}}{\sqrt{2}(2-\alpha L)}  & \mbox{if} \quad  \frac{2}{L+\mu}<\alpha < \frac{2}{L}. \\
\end{cases}
$$
%where 
%where $\rhogd := \max(|1-\alpha \mu|, |1-\alpha L|) \in (0,1)$.
\end{prop} 
\begin{proof}
%\proof {\textbf{Proof}.} 
By \eqref{ineq-strcvx-smooth} due to $L$-smoothness, a sufficient condition for the $\ell_2$ gain $\leq \gamma$ is %\eqref{def-hinfty-nonlinear} to hold is 
	\beq \frac{L}{2}\|{x}_k - x_*\|^2 - \gamma^2 \|{w}_k\|^2 \leq h_0 \left[f(x_0) - f(x_*) \right] % \quad \mbox{with} \quad {x}_0 - x_* = 0 
	\quad \mbox{for every} \quad k\geq 0,
	\label{ineq-hinf-gd}
	\eeq
for some positive constant $h_0$.We first use distance squared to the optimum as a Lyapunov function, which is 	
%The Lyapunov function $L(x) = \|x\|^2$	is frequently used 
standard in the analysis of GD methods. From the GD iterations subject to errors,
\beqs \|{x}_{k+1} -x_*\|^2 &=& \|{x}_{k}-\alpha\nabla f({x}_k) - x_*\|^2 + \alpha^2 \|{w}_k\|^2 + 2  \langle {x}_{k}-\alpha\nabla f({x}_k) - x_*, \alpha {w}_k \rangle \\
	&\leq&  \|{x}_{k}-\alpha\nabla f({x}_k) - x_*\|^2 + \alpha^2 \|{w}_k\|^2 + 2\| {x}_{k}-\alpha\nabla f({x}_k) - x_*\| \| \alpha {w}_k \|.
\eeqs
%Here, we note that the last inequality is tight in the sense that it would be an inequality for some choices of the (worst-case) noise ${w}_k$.
On the other hand, from the study of exact GD methods, it is well-known that for $\alpha \in (0,2/L)$,
 $$  \|{x}_{k}-\alpha\nabla f({x}_k) - x_*\| \leq \rhogd\|{x}_k - x_*\|, \quad \rhogd := \max(|1-\alpha \mu|, |1-\alpha L|),$$
(see \cite{lessard2016analysis}). Then, using this equality and Young's inequality, we obtain
\beqs \|{x}_{k+1} -x_*\|^2 
	&\leq&  \rhogd^2 \|{x}_{k} -x_*\|^2  + \alpha^2 \|{w}_k\|^2 + c\rhogd^2 \|{x}_{k} -x_*\|^2  +  \frac{1}{c} \alpha^2 \| {w}_k \|^2, 
\eeqs
for any $c>0$. Summing both sides for $k\geq 0$, we obtain
	$$  \sum_{k\geq 0}\|{x}_{k} -x_*\|^2 \leq H_{\tiny\mbox{GD}}(c)\|w\|_{\ell_2} + \frac{\|x_0 - x_*\|^2}{1-\rhogd^2(1+c)} \quad \mbox{with} \quad
	\quad H_{\tiny\mbox{GD}}(c):= \frac{\alpha^2 (1+\frac{1}{c})}{1-\rhogd^2(1+c)}, 
	$$
provided that $1-\rhogd^2(1+c)>0$. In this case, \eqref{ineq-hinf-gd} holds for $\gamma^2 \geq  \frac{L}{2}\frac{\alpha^2 (1+1/c)}{1-\rhogd^2(1+c)}$. The choice of $c = \frac{1}{\rhogd}- 1$ minimizes the right-hand side of this inequality and shows  \eqref{ineq-hinf-gd} holds for
$\gamma^2 = \frac{L}{2}  \frac{\alpha^2}{(1-\rhogd)^2}$ and $h_0 = \frac{\kappa}{1 - \rhogd^2(1+c)} = \frac{\kappa}{1-\rhogd} $ where we used the inequality $\frac{\mu}{2}\|x-x_*\|^2 \leq f(x)- f(x_*)$ which is a consequence of \eqref{ineq-strcvx-smooth}.
By taking square root of both sides, we conclude that 
\beq L_{2,*} \leq  \frac{\sqrt{L}}{\sqrt{2}}  \frac{\alpha}{1-\rhogd}  = \begin{cases}
\frac{1}{\sqrt{2\mu}} \sqrt{\kappa} & \mbox{if} \quad 0 < \alpha \leq \frac{2}{L+\mu}, \\
\frac{\sqrt{L}}{\sqrt{2}}  \frac{\alpha}{2-\alpha L} & \mbox{if} \quad  \frac{2}{L+\mu}<\alpha < \frac{2}{L},  \\
\end{cases} 
\label{ineq-gd-hinfty-bd1}
\eeq
where we used the definition of $\rhogd$. When the stepsize is large and close to $2/L$, this bound is tight in the sense that it matches the $H_\infty$ norm \eqref{eq-hinf-gd} we computed explicitly in the quadratic case (see also Theorem \ref{thm-real-hinf}, Corollary \ref{coro-real-hinf-equals-complex-sometimes} and \eqref{ineq-lower bound for l2 gain}). However, when the stepsize is sufficiently small, this bound is not as tight. For obtaining better bounds that are more aligned with the quadratic case, next we derive an alternative bound based on suboptimality $f_{k}:=f(x_k) - f(x_*)$ as a Lyapunov function for stepsize $\alpha \in (0,2/L)$. By smoothness, arguing similarly as before, 
\beq f_{k+1} &\leq& f_k - \alpha \langle \nabla f(x_k), \nabla f(x_k) +  w_k \rangle + \frac{L}{2}\alpha^2 \|\nabla f(x_k) + w_k \|^2 \nonumber\\
&\leq& f_k - \alpha (1-\frac{\alpha L}{2})\|\nabla f(x_k)\|^2 + \frac{L}{2}\alpha^2 \| w_k \|^2 -\alpha (1-\alpha L) \langle \nabla f(x_k), w_k \rangle  \nonumber\\
&\leq& f_k - \left( \alpha (1-\frac{\alpha L}{2}) - \frac{\alpha |1-\alpha L|r}{2}  \right) \|\nabla f(x_k)\|^2 + \frac{L}{2}\alpha^2 \| w_k \|^2 + \frac{\alpha |1-\alpha L|}{2} \frac{\|w_k\|^2}{r}  \nonumber \\
&\leq &  \left(1 -  2\mu \alpha (1-\frac{\alpha L}{2}) + \mu \alpha |1-\alpha L|r  \right)f_k + 
\alpha \left( \frac{L}{2}\alpha  + \frac{|1-\alpha L|}{2r}\right) \|w_k\|^2, \nonumber
\eeq
for any scalar $r>0$ where we used the Polyak-Lojasiewicz (PL) inequality $\|\nabla f(x_k)\|^2 \geq 2\mu f_k$ due to strong convexity \cite[Section 9.1.2]{boyd2004convex}. Summing the inequality above for $f_{k+1}$ over $k$,%\looseness=-1%and simplifying the terms,
$$ \sum_{k=0}^\infty f_k \leq   \frac{ \frac{L}{2}\alpha  + \frac{|1-\alpha L|}{2r}}{2\mu  (1-\frac{\alpha L}{2}) - \mu  |1-\alpha L|r }\sum_{k=0}^\infty \|w_k\|^2  + \frac{f_0}{ \alpha (2\mu  (1-\frac{\alpha L}{2}) - \mu  |1-\alpha L|r )},
$$
provided that $r>0$ satisfies $2\mu  (1-\frac{\alpha L}{2}) - \mu  |1-\alpha L|r >0$. If $\alpha \in (0,1/L)$, we choose $r=1$ so that $L_{2,*}^2 \leq \frac{1}{2\mu}$ for $\alpha \in (0, \frac{1}{L}]$.
%\beq 
%L_{2,*}^2  \leq  \frac{ \frac{L}{2}\alpha  + \frac{|1-\alpha L|}{2}}{2\mu  (1-\frac{\alpha L}{2}) - \mu  |1-\alpha L| } = \frac{1}{2\mu} \quad \mbox{for} \quad \alpha \in (0, \frac{1}{L}].
%\label{ineq-gd-hinfty-bd2}
%\eeq
On the other hand, if $\alpha \in (\frac{1}{L}, \frac{2}{L})$, we choose $r = \frac{2-\alpha L}{\alpha L}$ which yields $L_{2,*}^2  \leq \frac{(\alpha L)^2}{2\mu (2-\alpha L)^2 }$.
%\beq L_{2,*}^2  \leq  \frac{ \frac{L}{2}\alpha  + \frac{|1-\alpha L| \alpha L}{2(2-\alpha L)}}{2\mu  (1-\frac{\alpha L}{2}) - \mu  |1-\alpha L| \frac{(2-\alpha L)}{\alpha L}} = \frac{(\alpha L)^2}{2\mu (2-\alpha L)^2 } \quad \mbox{for} \quad \alpha \in ( \frac{1}{L},  \frac{2}{L}). 
%% \frac{ \frac{L}{2}\alpha  + \frac{(\alpha L - 1)\alpha L}{2(2-\alpha L)}}{\mu  (2-\alpha L) - \mu (\alpha L - 1) (2-\alpha L)/(\alpha L) } = \frac{(\alpha L) (\frac{L}{2}\alpha  + \frac{(\alpha L - 1)\alpha L}{2(2-\alpha L)}) }{\mu  (2-\alpha L)  }  \quad \mbox{for} \quad \alpha \in (1/L, 2/L). 
%\label{ineq-gd-hinfty-bd3}
%\eeq
%Combining \eqref{ineq-gd-hinfty-bd1}, \eqref{ineq-gd-hinfty-bd2}, \eqref{ineq-gd-hinfty-bd3} and 
Taking the pointwise minimum of these upper bounds on $L_{2,*}$ with the previous bound \eqref{ineq-gd-hinfty-bd1}, we obtain the desired result. %\myqed
%\endproof
\end{proof}

 %then GMM algorithm is linearly convergent with rate $\rho \in (0,1)$. In this sense, the Lyapunov function \eqref{eq-lyap-GMM} is a more suitable Lyapunov function for GMM methods.
%We consider the Lyapunov function
%$$ V_P({\xi}_k) :=  f({x}_k) - f(x_*)  +({\xi}_k^c)^T P {\xi}_k^c$$ 
%where $P = \tilde{P}\otimes I_d$ for some symmetric positive semi-definite $2\times 2$ matrix $\tilde{P}$ and ${\xi}_k^c$ is the centered initial iterate (around optimum) given by \eqref{def-centered-iter}. Such Lyapunov functions were considered in deterministic analysis of triple momentum methods where where obtaining the fastest convergence rates requires non-trivial non-zero choice of $\tilde{P}$ obtained through grid search \cn. This would suggest that when the noise level is small enough, non-zero choices of $P$ can potentially lead to tigher bounds in our analysis. 
% \subsection{Explicit $L_{2,*}$ bounds for NAG}
%\vspace{-0.2in}
Comparing Cor. \ref{coro-robust-gd} and Prop. \ref{prop-hinf-gd-bound}, we observe Prop. \ref{prop-hinf-gd-bound} is tight in the sense it leads to the same bounds as in the quadratic case for $\alpha \leq \frac{1}{L}$ or for $\alpha\in [\frac{2}{L+\sqrt{L\mu}},\frac{2}{L})$. For NAG and more generally for GMM, one difficulty for obtaining tight $L_{2,*}$ bounds is that the squared distance to the optimal solution $\| {x}_k - x_*\|^2$ is not necessarily monotonically decreasing even without noise %(see \cite[Theorem 5]{can2019accelerated})%\footnote{This is the case even for quadratics as the iteration matrix may contain a non-trivial Jordan block (see \cite[Proof of Theorem 5]{can2019accelerated}).}, %(see also footnote \ref{footnote-nonmonotone}),% due to the polynomially increasing prefactor term in the error (see e.g. \cite[Proof of Theorem 5]{can2019accelerated}), 
and the Lyapunov function we used for GD does not give tight estimates. Therefore, we introduce 
%For the analysis of (deterministic) GMM algorithms without noise, the following (more general) Lyapunov function has been introduced in the literature: % Therefore, for deterministic analysis of GMM methods the following Lyapunov function is introHowever, it is known  that if we introduce a more general Lyapunov function, 
%  \beq  V_P(x_k) := f({x}_k) - f(x_*) + ({x}_k - x_*)^T P ({x}_{k} - x_*)
%  \label{eq-lyap-GMM}
%  \eeq
   \beq  V_{P,c_1}({\xi}_k) := c_1 \left( f({x}_k) - f(x_*) \right) + \begin{bmatrix} {x}_k - x_* \\ {x}_{k-1} - x_*  \end{bmatrix}^T P\begin{bmatrix} {x}_k - x_* \\ {x}_{k-1} - x_*  \end{bmatrix},
  \label{eq-lyap-tmm-general}
  \eeq
which is a more general Lyapunov function where $c_1\geq 0$, $P= \tilde{P}\otimes I_d$ and $\tilde{P}$ is a $2\times 2$ positive semi-definite matrix. This Lyapunov function has been beneficial in the study of exact NAG and GMM methods \cite{fazlyab2018analysis, can2022entropic,aybat2018robust,hu2017dissipativity}. In the special case $c_1 = 0$, we require $\tilde{P}_{11}>0$.\footnote{Otherwise if $\tilde{P}_{11}=0$, we have necessarily $\tilde{P}_{12}=\tilde{P}_{21} = 0$ by the positive semi-definiteness of $\tilde{P}$ and $V_{P,c_1}$ would depend only on $x_{k-1}$ but not on the last iterate $x_k$ at step $k$ which would lead to a loose analysis.} In other words, we require $c_1\geq 0$, $\tilde{P}\succeq 0$ and $c_1 +  \tilde{P}_{11}>0$. We first provide a lemma for GMM that characterizes how this Lyapunov function evolves over the iterations. Its proof requires leveraging various inequalities obeyed by smooth strongly convex functions. This lemma will be useful to estimate $L_{2,*}$
of NAG and GMM with general parameter choices.
\begin{lemm}\label{lemma-gmm-lyap-evolution}Consider GMM algorithm subject to gradient errors $\{w_k\}$ with parameters $(\alpha, \beta, \nu)$ for minimizing $f\in \Cml$ with $\alpha>0$ and $\beta,\nu \geq 0$. Let $P = \tilde{P}\otimes I_d \in\mathbb{R}^{2d\times 2d}$ with $\tilde{P}\succeq 0$ and $c_1 + \tilde{P}_{11}>0$. For given non-negative scalars $\rho_0, \rho_1, \rho_2,\rho_3 \in [0,1)$ and $a,b,c_0,c_1 \geq 0$, %\mg{$b\geq 0$}, $c_1\geq 0$, and $c_0\geq 0$, 
it holds that 
{
\begin{eqnarray} V_{P,c_1}({\xi}_{k+1})   &\leq& (\rho_0^2 \mg{+c_1\rho_1^2 + \rho_2^2}) V_{P,c_1}({\xi}_{k})  - \begin{bmatrix} 
 	{\xi}_k^c \\
 	\nabla f({y}_k) \\
 	{w}_k
 \end{bmatrix}^T
 \left(
 (  \tilde{M}_2
   +c_1 \tilde{M}_1 +c_0 \tilde{M_0}) \otimes I_d \right) 
 \begin{bmatrix} 
 	{\xi}_k^c \\
 	\nabla f({y}_k) \\
 	{w}_k
 \end{bmatrix} \nonumber \\
 && + \frac{c_1 L}{2}\alpha^2 \|{w}_k\|^2 + w_k^T B^T P  B w_k  \mg{ + a\|w_k\|^2 + bc_1 \|\nabla f(y_k)\|^2} + \mg{\rho_3^2 V_{P,c_1}({\xi}_{k-1}) },  \label{ineq-gmm-lyap-evolution}
\end{eqnarray}
}
for $k\geq 0$ with the convention $\xi_{-1}:=\xi_0$ where
{\small
\beqs
\mbox{$\tilde{M}_0$} = \begin{bmatrix}   2\mu L\tilde{C}^T\tilde{C} & -(\mu+L)\tilde{C}^T& 0_{2 \times 1}\\ -(\mu+L)\tilde{C} & 2 & 0 \\
0_{1 \times 2} &0 & 0 
\end{bmatrix},~ \tilde{M}_1 = 
\left[ \begin{array}{@{}c|c@{}}
  \huge{\tilde{X}_0 + \tilde{Z}}
      & \begin{matrix} \frac{L\alpha \delta}{2} \tilde{M}_3^T \\ \frac{\alpha (1-L\alpha)}{2} \end{matrix} \\
   \cmidrule[0.4pt]{1-2}
   \begin{matrix} \frac{L\alpha \delta}{2} \tilde{M}_3 & \frac{\alpha (1-L\alpha)}{2} \end{matrix}  & 0 \\
\end{array} \right],
%\mbox{$\tilde{M}_1$} = \mbox{$\begin{bmatrix}
%     \mbox{${\tilde{Y_1}}$} &  \mbox{${\tilde{Y_2}}$} &  \mbox{$-\frac{L\alpha (\beta - \nu)}{2} \tilde{M}_3^T$} \\
%     \mbox{${\tilde{Y_2}^T}$} & \mbox{${\tilde{Y_3}}$} &  \frac{\alpha (1-L\alpha)}{2}\\
%     - \frac{L\alpha (\beta - \nu)}{2} \tilde{M}_3 & \frac{\alpha (1-L\alpha)}{2} & \mbox{\normalsize$0$}
% \end{bmatrix}$},
\eeqs
}
 %  -\tilde{A}^T \tilde{P} \tilde{B}   -\tilde{B}^T\tilde{P}\tilde{A}
 %\vspace{-0.2in}
 \beqs
\mbox{$\tilde{M}_2$} =  \begin{bmatrix}	 -\tilde{A}^T \tilde{P} \tilde{A} + \rho_0^2  %\mg{\rho_1^2 + \rho_2^2}
 \tilde{P}   & - \tilde{A}^T \tilde{P} \tilde{B}  &  - \tilde{A}^T \tilde{P} \tilde{B}  \\
 						   - \tilde{B}^T \tilde{P}  \tilde{A}   &  - \tilde{B}^T \tilde{P} \tilde{B} \mg{+bc_1}   &   - \tilde{B}^T \tilde{P}  \tilde{B}   \\
 						  - \tilde{B}^T \tilde{P} \tilde{A}   & - \tilde{B}^T \tilde{P} \tilde{B}  & \mg{a} % \tilde{B}^T \tilde{P}  \tilde{B} 
 \end{bmatrix},%\\
 %\mg{\mbox{$\tilde{M}_3$} =  c_1 \begin{bmatrix}   -\mu \nu^2 (2L-\mu)  - 2\mu \tilde{P}_{11} & - \mu  \nu^2 (2L-\mu) -2\mu \tilde{P}_{12} & \mu \nu & 0\\ 
%- \mu  \nu^2 (2L-\mu) -2\mu \tilde{P}_{12}  &   -\mu \nu^2 (2L-\mu)  - 2\mu \tilde{P}_{22}  & -\mu v & 0 \\
 %\mu v & -\mu v & -1 & 0 \\
 %0 & 0 & 0 & 0
%\end{bmatrix}},
%\\
%\mg{\tilde{M}_5  = \begin{bmatrix} - \rho_1 \tilde{P}_{11} - \rho_2 &- \rho_1\tilde{P}_{12} & 0 & \tilde{A}^T \tilde{P} \tilde{B}  \\ 
%- \rho_1\tilde{P}_{12} &  -\rho_1 \tilde{P}_{22}-\rho_3  & 0 & 0 \\
% 0& 0 & 0 & 0 \\
 % \tilde{B}^T\tilde{P}\tilde{A} & 0 & 0 & 0
%\end{bmatrix}}
 \eeqs 
are $4\times 4$ symmetric matrices,  $\tilde{A},\tilde{B},\tilde{C}$ are defined by \eqref{def: system mat for TMM}, $\tilde{M}_3 := \begin{bmatrix} 1 & - 1\end{bmatrix}, \quad \delta:= \beta - \nu$, 
\begin{equation} \tilde{Z} = \begin{bmatrix} \rho_1^2 \tilde{P}_{11} +\frac{\mu}{2} \rho_2^2 &  \rho_1^2 \tilde{P}_{12} & 0\\
 							   				 \rho_1^2 \tilde{P}_{12} & \rho_1^2 \tilde{P}_{22} +\frac{\mu}{2}\rho_3^2 & 0 \\
 							   				 					0   &				0			& 0
 \end{bmatrix} \in \mathbb{R}^{3\times 3}, \quad \tilde{X}_0 = \tilde{X}_1 + \rho_0^2 \tilde{X}_2 + (1-\rho_0^2) \tilde{X}_3  \in \mathbb{R}^{3\times 3}, %\tilde{M}_3 := \begin{bmatrix} 1 & - 1\end{bmatrix}, \quad \delta= \beta - \nu
 \label{def-tilde-X}
\end{equation}
% \beq \begin{scriptsize} \quad\quad
% %~\tilde{X}_0 = \tilde{X}_1 + \rho_0^2 \tilde{X}_2 + (1-\rho_0^2) \tilde{X}_3, ~ 
% \tilde{X}_1 = \frac{1}{2}
%\begin{bmatrix}
%-L\delta^2 & L\delta^2 & -(1-\alpha L)\delta \\
%L\delta^2 & -L\delta^2 & (1-L\alpha)\delta \\
%-(1-L\alpha)\delta & (1-L\alpha)\delta & \alpha(2-L\alpha)
%\end{bmatrix},
%\end{scriptsize}
%\label{def-tilde-X}
%\eeq 
%\vspace{-0.2in}
\begin{align*}\begin{scriptsize} %\begin{small}
 \tilde{X}_1 = \frac{1}{2}
\begin{bmatrix}
-L\delta^2 & L\delta^2 & -(1-\alpha L)\delta \\
L\delta^2 & -L\delta^2 & (1-L\alpha)\delta \\
-(1-L\alpha)\delta & (1-L\alpha)\delta & \alpha(2-L\alpha)
\end{bmatrix}, ~
\tilde{X}_2 = \frac{1}{2}
\begin{bmatrix} 
	\nu^2\mu& -\nu^2\mu & -\nu \\
    -\nu^2\mu & \nu^2\mu & \nu  \\
    -\nu  & \nu  & 0
\end{bmatrix}, 
\end{scriptsize}  
\end{align*}
\beqs
\tilde{X}_3 = \frac{1}{2}
\begin{bmatrix} 
	(1+\nu)^2\mu & -\nu(1+\nu)\mu & -(1+\nu) \\
    -\nu(1+\nu)\mu & \nu^2\mu & \nu  \\
    -(1+\nu)  & \nu  & 0 
\end{bmatrix}. 
\eeqs

\end{lemm} 
\begin{proof} 
%\proof {\textbf{Proof}.} 
The proof is given in Appendix \ref{app-Lyap-eval-GMM}. %\myqed
%\endproof
\end{proof}
%\vspace{-0.2in}
For NAG, we have $\delta=\beta-\nu = 0$ in which case the matrices in Lem. \ref{lemma-gmm-lyap-evolution} simplifies. Leveraging this fact, and extending the analysis of NAG methods without noise to the worst-case noise setting, we obtain an explicit bound on the $\ell_2$ gain of NAG in the following result. We use the common choice of the momentum $\beta = \frac{1-\sqrt{\alpha \mu}}{1-\sqrt{\alpha \mu}}$ \cite{aybat2018robust,aybat2019universally}.
%The following result provides . %relying on matrix inequality techniques. 
 \begin{prop}[Explicit $L_{2,*}$ bound for NAG]\label{prop-hinf-agd-bound}Consider minimizing $f \in \Cml$ with NAG with constant stepsize $\alpha \in (0,1/L]$, $\beta = \frac{1-\sqrt{\alpha\mu}}{1+\sqrt{\alpha\mu}}$ and initialization $x_0 = x_{-1}\in\mathbb{R}^d$. The worst-case robustness of NAG method satisfies the upper bound
{\small
\beq  \hspace{-0.13in} {\quad \Large L_{2,*}} &\leq& \overline{L}_{NAG}(\alpha):= %\bigg( 
\sqrt{ \frac{4 \mg{(5-2\sqrt{\alpha\mu} + \alpha\mu)} }{ \mu (1+\sqrt{\alpha\mu})^2 } +
\frac{\sqrt{\alpha}( 1 + \alpha L)}{\sqrt{\mu}} 
 + \frac{8 \mg{\alpha^3} L^4\left(4+(1-\sqrt{\alpha\mu})^2\right) }{ \mu^2 (1+\sqrt{\alpha\mu})^2 }  }\\
% \bigg)^{1/2}
 \label{ineq-hinf-bound-agd-by-hand} 
 &=&\frac{2\sqrt{5}}{\sqrt{\mu}} + \mathcal{O}(\sqrt{\alpha}) \nonumber %\mbox{ as } \alpha\to 0. \nonumber
 \eeq
 }
\end{prop} 

\begin{proof} 
%\proof {\textbf{Proof}.} 
We introduce $V_k := \mathcal{V}_{P,c_1}(\xi_k)$ and let $P = \tilde{P}\otimes I_d$ with the choice
{\small
\beq \tilde{P} = \begin{bmatrix} 
\sqrt{\frac{1}{2\alpha}} \\
\sqrt{\frac{\mu}{2}} - \sqrt{\frac{1}{2\alpha}} 
\end{bmatrix} \begin{bmatrix} 
\sqrt{\frac{1}{2\alpha}} &
\sqrt{\frac{\mu}{2}} -\sqrt{\frac{1}{2\alpha}} )
\end{bmatrix} = \frac{1}{2\alpha} \begin{bmatrix} 
1 \\
-(1-\sqrt{\alpha\mu)}
\end{bmatrix} \begin{bmatrix} 
1&~
-(1-\sqrt{\alpha\mu})
\end{bmatrix}.
\label{def-lyap-agd}
\eeq}Recall that for NAG, we have $\delta=\beta-\nu=0$. The inequality  %, as well. Therefore, 
\eqref{ineq-gmm-lyap-evolution} with our choice of $\beta$ and $ \rho_1=\rho_2=\rho_3 = a=b =  0$ is equivalent to
\beqs 
%\small
V_{k+1}   &\leq& \rho_0^2 V_k  - \begin{bmatrix} 
 	{\xi}_k^c \\
 	\nabla f({y}_k) \\
 	{w}_k
 \end{bmatrix}^T
 \left[ \begin{array}{@{}c|c@{}}
  \huge{\mathcal{S}_{\rho_0}(\tilde{P}) \otimes I_d}
      & \begin{matrix}   -A^T P B \\ Y \end{matrix} \\
   \cmidrule[0.4pt]{1-2}
   \begin{matrix} - B^T P A  &Y \end{matrix}  & 0 \\
\end{array} \right]
 \begin{bmatrix} 
 	{\xi}_k^c \\
 	\nabla f({y}_k) \\
 	{w}_k
 \end{bmatrix} \\
 && \quad 
 + \frac{c_1 L}{2}\alpha^2 \|{w}_k\|^2 + w_k^T B^T P  B w_k,
 \eeqs
 where  $Y := (c_1\frac{\alpha (1-L\alpha)}{2}I_d - B^T P B) $ and the square matrix
 {\footnotesize
 \beq \footnotesize ~\qquad\qquad \mathcal{S}_{\rho_0}(\tilde{P}) : = \begin{bmatrix}   -\tilde{A}^T \tilde{P} \tilde{A} + \rho_0^2 \tilde{P}   + c_0 2mL\tilde{C}^T\tilde{C}   &- \tilde{A}^T \tilde{P} \tilde{B}     -c_0 (m+L)\tilde{C}^T \\ - \tilde{B}^T \tilde{P} \tilde{A}  -c_0(m+L)\tilde{C}   & -  \tilde{B}^T \tilde{P} \tilde{B} +  2c_0 \\
\end{bmatrix} + c_1 (\tilde X_0 + \tilde Z), 
\label{def-S-rho}
\eeq
}is the leading $3\times 3$ principal submatrix of $\tilde{M_2} + c_1\tilde{M}_1 + c_0 \tilde{M}_0$. If we also take $c_1 = 1$ and $c_0=0$, then it follows after some straightforward computations (similar to those in \cite[App. C]{aybat2019universally}) that $\mathcal{S}_{\rho_0}(\tilde{P}) \succeq 0$ with $\rho_0^2 = \rho^2_{\tiny \mbox{NAG}}:=1 - \sqrt{\alpha\mu}$. %\cite[Appendix C]{aybat2019universally}. %\cite[Thm. 2.3]{aybat2019universally}.  
Then, using $\tilde{P}_{11}=\frac{1}{2\alpha}$  and $B^T P B = \tilde{P}_{11}\alpha^2 I_d = \frac{\alpha}{2}I_d$ yields%we obtain %from \eqref{eq-evolution-nonquad-Lyap} and \eqref{eq-evolution-quad-Lyap} that
\beq
V_{k+1} &\leq& \rho^2_{\tiny \mbox{NAG}} V_k + \frac{L\alpha^2}{2}\|w_k\|^2 + \frac{\alpha}{2} \|w_k\|^2 + 2w_k^T B^T P A \xi_k^c + \mg{\alpha^2} L 
 \nabla f(y_k)^T w_k,\\
\label{ineq-lyap-decay-nag}
% \eeq
%We also have
%\beq 
 \| B^T P A \xi_k^c\|^2 &=& \alpha^2 \big\| (\tilde{P}_{11} (1+\beta) + \tilde P_{12}) (x_k - x_*) -\beta \tilde{P}_{11} (x_{k-1} - x_*) \big\|^2 =  \frac{1}{4(1+\sqrt{\alpha\mu})^2} I_k(\alpha), \label{def-I-k-alpha}
\eeq
with $I_k(\alpha) := \big\| ( 1 + \alpha \mu ) (x_k - x_*) - (1-\sqrt{\alpha\mu})(x_{k-1} - x_*) \big\|^2$. Moreoever, 
{\small
\beq %\hspace{-0.24in}
 I_k(\alpha) %&:=& \big\| ( 1 + \alpha \mu ) (x_k - x_*) - (1-\sqrt{\alpha\mu})(x_{k-1} - x_*) \big\|^2\\ %\nonumber \\
 %&=& 2\alpha (\xi_k^c)^T P (\xi_k^c) + (\alpha^2 \mu^2 + \mg{2\alpha\mu}) \|x_k - x_*\|^2 \nonumber \\
  %&& 
  %- 2\alpha \mu(1-\sqrt{\alpha\mu}) (x_k - x_*)^T (x_{k-1}-x_*) \nonumber \\
  &\leq & 2\alpha (\xi_k^c)^T P (\xi_k^c) + (\alpha^2 \mu^2 + \mg{2\alpha\mu} + \alpha \mu (1-\sqrt{\alpha\mu} )) \|x_k - x_*\|^2  +\alpha \mu(1-\sqrt{\alpha\mu}) \|x_{k-1}-x_*\|^2  \label{ineq-I-k-alpha}\\
 &\leq& 2\alpha V_k + (2 \alpha^2 \mu + \mg{4\alpha}) (f(x_k) - f(x_*)) %\nonumber \\
 %&&  
 + \alpha\mu (1-\sqrt{\alpha\mu}) (\|x_k - x_*\|^2 +\|x_k - x_*\|^2) \nonumber \\
 &\leq& %2\alpha (\mg{3}+\alpha\mu) V_k + 2 \alpha (1-\sqrt{\alpha\mu}) (V_k + V_{k-1}) \nonumber \\
% &=&
 2\alpha \mg{( 4-\sqrt{\alpha\mu} + \alpha\mu)} V_k +  2 \alpha (1-\sqrt{\alpha\mu}) V_{k-1}, \label{ineq-Ik-alpha-final}
 %4\alpha V_k +  2 \alpha (1-\sqrt{\alpha\mu}) V_{k-1},
\eeq
}where %in the first and the second inequality 
we used $\frac{\mu}{2}\|x_j-x_*\|^2 \leq f(x_j) - f(x_*) \leq V_j$ (which holds because of \eqref{ineq-strcvx-smooth} and $c_1 = 1$) for $j=k-1$ and $j=k$. Using the inequality $2ab \leq r a^2 + \frac{b^2}{r}$ for any scalars $a,b$ and for any $r>0$, Cauchy-Schwarz inequality, the inequalities \eqref{ineq-lyap-decay-nag}, \eqref{def-I-k-alpha} and \eqref{ineq-Ik-alpha-final}; we have for any scalars $s_1, s_2>0$, 
\beqs
V_{k+1}% &\leq& %\rho^2_{\tiny \mbox{NAG}} V_k + \frac{\alpha( 1 + \alpha L)}{2}\|w_k\|^2 + 2w_k^T B^T P A \xi_k^c + \mg{\alpha^2} %\sout{\alpha}  
%L 
 %\nabla f(y_k)^T w_k \\
 &\leq&  \rho^2_{\tiny \mbox{NAG}} V_k + \frac{\alpha( 1 + \alpha L)}{2}\|w_k\|^2  + \frac{ \mg{\alpha^2}  L}{2} 
( \frac{\| \nabla f(y_k)\|^2}{s_2}  + s_2 \|w_k\|^2)  + s_1 {\|w_k\|^2} \\
 && +   \frac{ 2\alpha \mg{( 4-\sqrt{\alpha\mu} + \alpha\mu)} V_k +  2 \alpha (1-\sqrt{\alpha\mu}) V_{k-1}}{4s_1(1+\sqrt{\alpha\mu})^2} \\
% & \leq & \rho^2_{\tiny \mbox{NAG}} V_k + \frac{\alpha( 1 + \alpha L)}{2}\|w_k\|^2  + \frac{ \mg{\alpha^2}  L}{2} 
%( \frac{\| \nabla f(y_k)\|^2}{s_2}  + s_2 \|w_k\|^2)  + s_1 {\|w_k\|^2} \\
 %&& +   \frac{ 2\alpha \mg{( 4-\sqrt{\alpha\mu} + \alpha\mu)} V_k +  2 \alpha (1-\sqrt{\alpha\mu}) V_{k-1}}{4s_1(1+\sqrt{\alpha\mu})^2}\\
 & \leq & \rho^2_{\tiny \mbox{NAG}} V_k + \frac{\alpha( 1 + \alpha L)}{2}\|w_k\|^2  + \frac{\mg{\alpha^2}  L}{2} 
( \frac{2L^2 (1+\beta)^2 \| x_k - x_*\|^2 + 2L^2 \beta^2 \|x_{k-1}-x_*\|^2}{s_2}  )  \\
&&+ \frac{\mg{\alpha^2}  L}{2}  s_2 \|w_k\|^2 + s_1 {\|w_k\|^2} 
 +   \frac{ 2\alpha \mg{( 4-\sqrt{\alpha\mu} + \alpha\mu)} V_k +  2 \alpha (1-\sqrt{\alpha\mu}) V_{k-1}}{4s_1(1+\sqrt{\alpha\mu})^2}\\ 
 & \leq & \rho^2_{\tiny \mbox{NAG}} V_k + \frac{\alpha( 1 + \alpha L)}{2}\|w_k\|^2  + { 2\mg{\alpha^2}  L^3}
( \frac{ (1+\beta)^2 V_k +  \beta^2 V_{k-1}}{s_2 \mu}  )  \\
&&+ \frac{\mg{\alpha^2}  L}{2}  s_2 \|w_k\|^2 + s_1 {\|w_k\|^2} 
 +   \frac{ 2\alpha \mg{( 4-\sqrt{\alpha\mu} + \alpha\mu)} V_k +  2 \alpha (1-\sqrt{\alpha\mu}) V_{k-1}}{4s_1(1+\sqrt{\alpha\mu})^2},
 \eeqs
with the convention that $V_{-1} := V_0$, where in the second and third inequalities, we used $L$-smoothness and strong convexity of $f$. We choose $s_i = \sqrt{\alpha} \hat{s}_i$ for $i=1,2$ for some $\hat{s}_i >0 $ that we will specify next. Using $1+\beta = 2/(1+\sqrt{\alpha\mu})$, for $k\geq 1$, %Letting $\hat{s}_2 = \frac{2}{\alpha L}\hat{s}_1$,
 \beq
  V_{k+1} &\leq& \Big(\rho^2_{\tiny \mbox{NAG}} +  \frac{8\mg{\alpha}\sqrt{\alpha} L^3}{\mu (1+\sqrt{\alpha\mu})^2 \hat{s}_2}  + \frac{\sqrt{\alpha}\mg{(4-\sqrt{\alpha\mu} + \alpha\mu)}}{\mg{2}\hat{s}_1 (1+\sqrt{\alpha\mu})^2 } \Big) V_k + \mg{\bigg(}\frac{\sqrt{\alpha} (1- \sqrt{\alpha\mu}) }{2\hat{s}_1  (1+\sqrt{\alpha\mu})^2 } \nonumber \\
  && \mg{+ \frac{2\alpha \sqrt{\alpha} L^3 \beta^2}{\hat{s}_2 \mu} \bigg)}
  V_{k-1} + \bigg( \frac{\alpha( 1 + \alpha L)}{2} + \sqrt{\alpha}\hat{s}_1 + \frac{\mg{\alpha^2} \sqrt{\alpha} L}{2}\hat{s}_2 \bigg)\|w_k\|^2. 
  \label{ineq-to-sum-side-by-side-nag}
 \eeq
We choose 
 \beq \hat{s}_1 = \frac{2\mg{(5-2\sqrt{\alpha\mu} + \alpha\mu)}}{\sqrt{\mu} (1+\sqrt{\alpha\mu})^2}, \quad  \hat{s}_2 =  \frac{8L^3 \mg{\alpha} }{\mu\sqrt{\mu} (1+\sqrt{\alpha\mu})^2}\left(4+(1-\sqrt{\alpha\mu})^2\right),
 \label{def-hat-r1-r2}
 \eeq
so that we have
{\small
\beq \frac{8\mg{\alpha}\sqrt{\alpha} L^3}{\mu (1+\sqrt{\alpha\mu})^2 \hat{s}_2}  +  \mg{\frac{2\alpha \sqrt{\alpha} L^3 \beta^2}{\hat{s}_2 \mu}}
= 
\frac{\sqrt{\alpha \mu}}{4} =   \frac{\sqrt{\alpha} \mg{(4-\sqrt{\alpha\mu} + \alpha \mu)}}{ \mg{2}\hat{s}_1 (1+\sqrt{\alpha\mu})^2 } + \frac{\sqrt{\alpha} (1- \sqrt{\alpha\mu}) }{2\hat{s}_1  (1+\sqrt{\alpha\mu})^2 }. \label{def-r1-r2}
\eeq
}
Summing the inequality \eqref{ineq-to-sum-side-by-side-nag} for $k=0, 1, 2,\dots, K$ with $V_{-1}:=V_0$, and rearranging terms
\beq   && \frac{V_{K+1}}{1 - \rho^2_{\tiny \mbox{NAG}} - \frac{\sqrt{\alpha\mu}}{4}  - \frac{\sqrt{\alpha\mu}}{4} } + \sum_{k=0}^K V_k \\&& \quad \leq \frac{1}{1 - \rho^2_{\tiny \mbox{NAG}} - \frac{\sqrt{\alpha\mu}}{4}  - \frac{\sqrt{\alpha\mu}}{4} } \big( \frac{\alpha( 1 + \alpha L)}{2} + \sqrt{\alpha}\hat{s}_1 + \frac{\mg{\alpha^2} \sqrt{\alpha} L}{2}\hat{s}_2 \big) \sum_{k=0}^K \|w_k\|^2 \nonumber \\
&&\quad \quad + \frac{1}{1 - \rho^2_{\tiny \mbox{NAG}} - \frac{\sqrt{\alpha\mu}}{4} - \frac{\sqrt{\alpha\mu}}{4} }   \big(1 +  \frac{\sqrt{\alpha} (1- \sqrt{\alpha\mu}) }{2\hat{s}_1  (1+\sqrt{\alpha\mu})^2 } + \frac{2\alpha \sqrt{\alpha} L^3 \beta^2}{\hat{s}_2 \mu}  \big) V_0 \nonumber \\ %\label{ineq-nag-cum-lyap-bound}
%&=& 
%\frac{2}{\sqrt{\alpha\mu}}  \bigg( \frac{\alpha( 1 + \alpha L)}{2} + \frac{2\sqrt{\alpha}\mg{(5-2\sqrt{\alpha\mu} + \alpha\mu)}}{\sqrt{\mu} (1+\sqrt{\alpha\mu})^2} %\frac{2\sqrt{\alpha}(3-\sqrt{\alpha\mu}) }{ \sqrt{\mu} (1+\sqrt{\alpha\mu})^2 }
%+\frac{4L^4 \mg{\alpha^3\sqrt{\alpha} \left(4+(1-\sqrt{\alpha\mu})^2\right)} }{\mu\sqrt{\mu} (1+\sqrt{\alpha\mu})^2} % + \frac{4 \alpha \sqrt{\alpha} L^4 }{ \mu \sqrt{\mu} (1+\sqrt{\alpha\mu})^2 }  
%\bigg) \sum_{k=0}^\infty \|w_k\|^2 
% \\
% &&+ \frac{2}{\sqrt{\alpha\mu}} 
% \bigg(1 +  \frac{\sqrt{\alpha\mu} (1- \sqrt{\alpha\mu}) }{4 (5-2\sqrt{\alpha\mu} +\alpha\mu) }  + \frac{\sqrt{\alpha\mu} (1-\sqrt{\alpha\mu})^2}{4(4+ (1-\sqrt{\alpha\mu})^2)}  \bigg) V_0 \\
&&\quad =   %\bigg( \frac{\sqrt{\alpha}( 1 + \alpha L)}{\sqrt{\mu}} + \frac{4 \mg{(5-2\sqrt{\alpha\mu} + \alpha\mu)} }{ \mu (1+\sqrt{\alpha\mu})^2 }
 %+ \frac{8 \mg{\alpha^3} L^4 \left(4+(1-\sqrt{\alpha\mu})^2\right) }{ \mu^2 (1+\sqrt{\alpha\mu})^2 }  \bigg) 
  \big(\overline{L}_{NAG}(\alpha)\big)^2 \sum_{k=0}^K \|w_k\|^2  + \overline{H}(x_0), \quad \mbox{where} \quad
  	 \label{ineq-nag-cum-lyap-bound}
\eeq
the last term $\overline{H}(x_0):= \frac{2}{\sqrt{\alpha\mu}} 
 	\bigg(1 +  \frac{\sqrt{\alpha\mu} (1- \sqrt{\alpha\mu}) }{4 \mg{(5-2\sqrt{\alpha\mu} + \alpha\mu)} } + 		\frac{\sqrt{\alpha\mu} (1-\sqrt{\alpha\mu})^2}{4(4+ (1-\sqrt{\alpha\mu})^2)} \bigg) 
 V_0$.
%{\footnotesize
%\beq \overline{H}(x_0):= \frac{2}{\sqrt{\alpha\mu}} 
% 	\bigg(1 +  \frac{\sqrt{\alpha\mu} (1- \sqrt{\alpha\mu}) }{4 \mg{(5-2\sqrt{\alpha\mu} + \alpha\mu)} } + 		\frac{\sqrt{\alpha\mu} (1-\sqrt{\alpha\mu})^2}{4(4+ (1-\sqrt{\alpha\mu})^2)} \bigg) 
% V_0
% \label{def-H-x0}
% \eeq
%} 
Then, letting $K\to \infty$, using $V_k\geq f(x_k)-f(x_*)$ and the definition of $L_{2,*}$, we conclude. %\myqed %\eqref{def-hinfty-nonlinear}.
\end{proof}
%\endproof
%\mtodo{Mention that the params in MI is there because we want an analysis that recovers hand-made GD and NAG analysis. Give the params for MI which implies the hand-made GD and NAG bounds. Say that finite time bounds $\sum_{k =0}^K f(x_k) -f(x_*)$ can also be obtained from our proof technique easily but here for simplicity we let $K\to \infty$.}
%\vspace{-0.2in}
\begin{coro}[$\mathcal{O}(1/k)$ ergodic rate for inexact NAG]\label{coro-ergodic-rates} In the setting of Prop. \ref{prop-hinf-agd-bound}, the ergodic average $\bar{x}_K = \frac{1}{K+1} \sum_{j=0}^K x_K$ satisfies
$ [f(\bar x_K) - f(x_*) ] \leq \frac{1}{K+1} \sum_{j=0}^K [f(x_j) - f(x_*)] \leq \big( \overline{L}_{NAG}(\alpha) \big)^2 \frac{\sum_{j=0}^{K} \|w_j\|^2}{K+1} + \frac{H(x_0)}{K+1} - \frac{2}{\sqrt{\alpha\mu}}\frac{f(x_{K+1}) -f(x_*) }{K+1} $
where $\overline{H}(x_0)$ is as in \eqref{ineq-nag-cum-lyap-bound}.
%with $H(x_0) = \frac{2}{\sqrt{\alpha\mu}} 
% \bigg(1 +  \frac{\sqrt{\alpha\mu} (1- \sqrt{\alpha\mu}) }{4 \mg{(5-2\sqrt{\alpha\mu} + \alpha\mu)} } + \frac{\sqrt{\alpha\mu} (1-\sqrt{\alpha\mu})^2}{4(4+ (1-\sqrt{\alpha\mu})^2)} \bigg) V_0$ for any $k\geq 0$. % and $w \in \ell_2(\mathbb{R}^d)$.
 %\sum_{k=1}^K f(x_k) - f(x_0)%By the convexity, ergodic iterate% $\bar{x}_K = \sum_{k=1}^K \frac{1}{K} x_k$.
%While \eqref{ineq-nag-cum-lyap-bound} has an asymptotic nature, it is also straightforward to obtain non-asymptotic bounds on the cumulative sums $\sum_{k=0}^K \[f(x_k) - f(x_*)\]$ in terms of the noise sequence $\{w_k\}_{k=0}^K$ from the recurrence \eqref{ineq-to-sum-side-by-side-nag}.
\end{coro}
%\proof {\textbf{Proof}.} 
\begin{proof}
This follows from the convexity of $f$, inequalities \eqref{ineq-to-sum-side-by-side-nag}, \eqref{ineq-nag-cum-lyap-bound}, $V_k = \mathcal{V}_{P,c_1}(\xi_k) \geq f(x_k)-f(x_*) 
\geq 0$ for $c_1 = 1$ combined with $1 - \rho^2_{\tiny \mbox{NAG}} - \frac{\sqrt{\alpha\mu}}{4}  - \frac{\sqrt{\alpha\mu}}{4} = \frac{\sqrt{\alpha\mu}}{2}$. %where we recall that $\rho^2_{\tiny \mbox{NAG}}=1 - \sqrt{\alpha\mu}$. 
%\myqed 
%\endproof
\end{proof}
%\vspace{-0.2in}
\begin{rema}[Tightness of NAG and GD  analysis]\label{remark-tightness-nag} By \eqref{ineq-lower bound for l2 gain}, Theorems \ref{thm-h-inf} and \ref{thm-real-hinf}, there are quadratic functions for which $L_{2,*}=H_\infty$ with $H_\infty = \frac{1}{\sqrt{2\mu}}$ (see Table \ref{table}). Therefore, for fixed parameters $(\alpha,\beta,\nu)$, we have $\sup_{f\in\Cml} L_{2,*}\geq \frac{1}{\sqrt{2\mu}}$ and the lower bound is attained for certain choice of parameters and quadratic $f$. Proposition \ref{prop-hinf-agd-bound} for NAG is tight in the small stepsize regime in the sense that when stepsize is small enough, it gives the bound $L_{2,*} \leq \frac{2\sqrt{10}}{\sqrt{2\mu}}$ which matches the lower bound up to a universal constant (of $2\sqrt{10}$). Similarly for GD, comparing Prop. \ref{prop-hinf-gd-bound} and Coro. \ref{coro-robust-gd}, we see that our analysis for $L_{2,*}$ is tight when $\alpha \leq 1/L$ or for $\alpha \in [\frac{2}{L+\sqrt{L\mu}}, \frac{2}{L})$.
\end{rema}
%\vspace{-0.1in}
\subsection{Bounding $L_{2,*}$ for general GMM parameters with matrix inequalities.}\label{subsec-str-cvx-MI}

%$V_{P,c_1}(\xi_k) = \tilde{P}_{11} \|x_k - x_*\|^2 $ vanishes only when $x_k = x_*$. For $c_1>0$, $\tilde{P}$ can be an arbitrary $2\times 2$ matrix as long as $\tilde{P}\succeq 0$.
  In the abscence of gradient errors, a linearly convergence rate $\rho$ can be certified for GMM methods if $\rho$ and $\tilde{P}$ satisfy a $3\times 3$ matrix inequality \cite{hu2017dissipativity}. Unfortunately, the convergence rate $\rho$ and the Lyapunov matrix $\tilde{P}$ pair is only known explicitly for some particular choice of the parameters $(\alpha,\beta,\nu)$ to our knowledge \cite{can2022entropic,van2021speed} and in the general case they can be computed numerically with a simple grid search. In the following, we will derive a new $4\times 4$ matrix inequality (MI) based on Lemma \ref{lemma-gmm-lyap-evolution}; if the parameters $\rho$ and $\tilde{P}$ satisfy this MI, then we will obtain an immediate upper bound on $L_{2,*}$. This approach recovers the bounds we derived for GD and NAG (in the sense that by choosing the parameters of this MI by hand, we can obtain the same explicit bounds we obtained for GD and NAG in Prop. \ref{prop-hinf-gd-bound} and Prof. \ref{prop-hinf-agd-bound}, see Remark \ref{rema-MI-explicit-bounds}) and generalizes it to more general parameter choice of GMM. In other words, our MI is motivated by our explicit analysis, and is designed to have enough parameters and flexibility to recover our explicit bounds.%\looseness=-1%After all, an MI can be thought of %After all, MIs are equivalent to inequalities. %Basically, we can view the MI approach as a way of bundling various useful inequalities satisfied by the iterates and gradients due to strong convexity and smoothness.
%In the next result, we provide an upper bound on the $H_\infty$ norm of TMM algorithm in terms of solutions of a small dimensional ($4\times 4$) matrix inequality.  
%By the positive-semi definiteness of $P$, we have  $V_P({\xi}_k^c) \geq  f({x}_k^c) - f_*$. Therefore, if 
%	\beq \sum_{k=0}^\infty V_P({\xi}_k) - \gamma^2 \|w_k\|^2 \leq 0
%	\label{ineq-suff-cond-nonlinear-hinf}
%	\eeq
%then \eqref{def-hinfty-nonlinear} is satisfied and the $H_\infty$ norm is smaller than $\gamma$. 
\begin{theo}[$L_{2,*}$ bound for GMM]\label{thm-hinfty-bound} Consider GMM algorithm subject to worst-case noise with parameters $(\alpha, \beta, \nu)$ for minimizing $f\in \Cml$ with $\alpha>0, \beta,\nu \geq 0$ and $x_0 = x_{-1}\in \mathbb{R}^d$. For given non-negative scalars $\rho_0, \rho_1, \rho_2,\rho_3 \in [0,1),  a\geq 0$, \mg{$b\geq 0$}, $c_1\geq 0$, and $c_0\geq 0$. Assume that $\tilde{P}\in\mathbb{R}^{2\times 2}$ is a positive semi-definite matrix that satisfies $c_1 + \tilde{P}_{11}>0$ and
%\beq  
 $\tilde{M}_4:=\tilde{M}_2 +c_1 \tilde{M_1} + c_0 \tilde{M}_0 \succeq 0$,
%\label{constraint-sdp}
%\eeq
where $\tilde{M}_0, \tilde{M}_1$ and $\tilde{M}_2$ are $4\times 4$ as in Lemma \ref{lemma-gmm-lyap-evolution}. %Then,  real symmetric matrices defined as
Then, the robustness of GMM satisfies%the worst-case robustness %in terms of $L_{2,*}$ 
%of the GMM algorithm admits %the bound %
%\beq h({a,\rho,c_1,c_0}) := && ~ \min_{\tilde{P}\in\mathbb{R}^{2\times 2}, \tilde{P}\succeq 0}  {\frac{\alpha^2}{1-\rho^2} (\frac{c_1L}{2}+ (1+a) \tilde{P}_{11})} \label{eq-sdp-tomin}\\
%&&\mbox{subject to }\tilde{M}_2(\tilde{P}) +c_1 \tilde{M_1} + c_0 \tilde{M}_0 \succeq 0. %\rho^2 \in [0,1),  a\geq 0 
%\eeq
% $$ 
% A=\begin{pNiceArray}{cc|cc}[baseline=line-3]
% A & B & 0 & 0 \\
% X & X & 0 & 0 \\
% \hline
% 0 & 0 & A & B \\
% 0 & 0 & D & D \\
% \end{pNiceArray}
% $$
\beq %\quad\quad 
L_{2,*} \leq \overline{L}_{GMM}(\alpha,\beta,\nu)%L_{2,*}^{ub}
:= \sqrt{
\frac{1 }{1-s} {  \tfrac{{
\alpha^2 ( \frac{c_1L}{2}+  \tilde{P}_{11}) \mg{+ a} } }{{ c_1 + \frac{2}{L} r(\tilde{P}) }}, 
} 
} %\mbox{  where  } \begin{small} r(\tilde P) := \begin{cases}
 %\tilde{P}_{11} - {\tilde{P}_{12}^2}/{\tilde{P}_{22}} & \mbox{if} \quad \tilde{P}_{22}\neq 0, \\
% \tilde{P}_{11} & \mbox{if} \quad \tilde{P}_{22}=  0,
  %\end{cases}\end{small}
\label{def-Hinfty-ub}
\eeq
provided\footnote{We use the convention that if $c_1 + \frac{2}{L} r(\tilde{P}) = 0$, then $\overline{L}_{GMM}(\alpha,\beta,\nu):=\infty$ in which case we obtain a trivial bound.} that $\begin{small} s:=\rho_0^2+c_1\rho_1^2 +\rho_2^2 + \rho_3^2  + \frac{4b(\nu^2 + (1+\nu)^2) L^2}{\mu} c_1 <1\end{small} $ where 
%$\begin{small} r(\tilde P) := \begin{cases}
% \tilde{P}_{11} - {\tilde{P}_{12}^2}/{\tilde{P}_{22}} & \mbox{if} \quad \tilde{P}_{22}\neq 0, \\
 %\tilde{P}_{11} & \mbox{if} \quad \tilde{P}_{22}=  0.
  %\end{cases}\end{small}$ 
  %\mg{ $\operatorname {sgn}(c_1):= 1$ if $c_1>0$ and $\operatorname {sgn}(c_1):= 0$ otherwise}, 
%and 
\beq
 r(\tilde P) = \begin{cases}
 \tilde{P}_{11} - {\tilde{P}_{12}^2}/{\tilde{P}_{22}} & \mbox{if} \quad \tilde{P}_{22}\neq 0, \\
 \tilde{P}_{11} & \mbox{if} \quad \tilde{P}_{22}=  0.
  \end{cases}
  \label{def-r-P}
\eeq
\end{theo} 
%\mg{Modify the proof as in our RAGM submission so that it cover GD as a special case. One should have a coefficient $c_0$ in front of $L/2$ in the formula for $c_{a,\rho}$ below and potentially recover the GD solution.}
%\begin{proof} 
%Using $y_k = \tilde{C}x_k$ and ?, we have
%\proof {\textbf{Proof}.} 
\begin{proof}
Using Lemma \ref{lemma-gmm-lyap-evolution}, the fact that $f\in\Cml$ and the MI constraint $M_4\succeq 0$,%$\tilde{M}_2 +c_1 \tilde{M_1} + c_0 \tilde{M}_0 \succeq 0$,
\begin{small}
\beq V_{P,c_1}({\xi}_{k+1})   
%&\leq& (\rho_0^2 \mg{+c_1\rho_1^2 + \rho_2^2}) V_{P,c_1}({\xi}_{k})  - \begin{bmatrix} 
% 	{\xi}_k^c \\
% 	\nabla f({y}_k) \\
% 	{w}_k
% \end{bmatrix}^T
% \left(
% (  \tilde{M}_2
%   +c_1 \tilde{M}_1 +c_0 \tilde{M_0}) \otimes I_d \right) 
% \begin{bmatrix} 
% 	{\xi}_k^c \\
% 	\nabla f({y}_k) \\
% 	{w}_k
% \end{bmatrix} \nonumber \\
% && + \frac{c_1 L}{2}\alpha^2 \|{w}_k\|^2 + w^T B^T P  B w  \mg{ + a\|w_k\|^2 + bc_1 \|\nabla f(y_k)\|^2} + \mg{\rho_3^2 V_{P,c_1}({\xi}_{k-1}) }  \label{ineq-evolution-Lyap-with-noise}\\
 &\leq& (\rho_0^2  \mg{+c_1\rho_1^2 + \rho_2^2}  )V_{P,c_1}({\xi}_{k}) + \frac{c_1 L}{2}\alpha^2 \|{w}_k\|^2 + w_k^T B^T P  B w_k \nonumber \\
 &&\quad \mg{+a\|w_k\|^2 + b c_1 L^2 \|y_k - x_*\|^2} +\mg{\rho_3^2 V_{P,c_1}({\xi}_{k-1}) }   \label{ineq-evolution-Lyap-with-noise}\\
 &\leq &  (\rho_0^2  \mg{+c_1\rho_1^2 + \rho_2^2} ) V_{P,c_1}({\xi}_{k}) + \frac{c_1 L}{2}\alpha^2 \|{w}_k\|^2 + \alpha^2 \tilde{P}_{11} \|{w}_k\|^2 + \mg{a\|w_k\|^2}\nonumber\\
&& \quad \mg{+ 2bc_1 (1+\nu)^2 L^2 \|x_k - x_*\|^2 + 2bc_1 \nu^2 L^2 \|x_{k-1} - x_*\|^2} \mg{+\rho_3^2 V_{P,c_1}({\xi}_{k-1}) } ,\nonumber\\
% &\leq &  (\rho_0^2  \mg{+c_1\rho_1^2 + \rho_2^2}  ) V_{P,c_1}({\xi}_{k}) + \frac{c_1 L}{2}\alpha^2 \|{w}_k\|^2 + \alpha^2 \tilde{P}_{11} \|{w}_k\|^2 + \mg{a\|w_k\|^2} \nonumber\\
%&& \quad \mg{+ \frac{4bc_1 (1+\nu)^2 L^2}{\mu} (f(x_k) - f(x_*)) +  \frac{4bc_1\nu^2 L^2}{\mu} (f(x_{k-1}) - f(x_*))} + \mg{\rho_3^2 V_{P,c_1}({\xi}_{k-1}) } \nonumber\\
&\leq &   (\rho_0^2  \mg{+c_1\rho_1^2 + \rho_2^2}  ) V_{P,c_1}({\xi}_{k}) + \frac{c_1 L}{2}\alpha^2 \|{w}_k\|^2 + \alpha^2 \tilde{P}_{11} \|{w}_k\|^2 +\mg{a\|w_k\|^2} \nonumber \\
&& \quad \mg{+  \frac{4b(1+\nu)^2 L^2}{\mu} c_1 V_{P,c_1}({\xi}_{k})  +  \frac{4b\nu^2 L^2}{\mu} c_1 V_{P,c_1}({\xi}_{k-1})} + \mg{\rho_3^2 V_{P,c_1}({\xi}_{k-1}) }, \nonumber
 %\label{ineq-Vpc-2} 
 \eeq\end{small}where $\xi_{-1}=\xi_0$. %where in the second inequality, we used the inequality \eqref{constraint-sdp} and $L$-smoothness of the gradient. %and in the last step we used strong convexity.
%If $\tilde{M}_a \succeq 0$, then we get $M_a \succeq 0$ and
%\beq V_P({\xi}_{k+1})   &\leq& \rho^2 V_P({\xi}_{k}) + \frac{L}{2}\alpha^2 \|{w}_k\|^2 + (1+a) {w}_k^T B^T P  B {w}_k \\
%	 &\leq& \rho^2 V_P({\xi}_{k}) + \frac{L}{2}\alpha^2 \|{w}_k\|^2 + (1+a)\alpha^2 \tilde{P}_{11} \|{w}_k\|^2
%\eeq
Summing the last inequality for $k=0,1,2, \dots, K$, and reorganizing the terms, %we obtain
\beq  \frac{V_{P,c_1}({\xi}_{K+1})}{1-s} &+& \sum_{k=0}^K V_{P,c_1}({\xi}_{k})\nonumber \\  &\leq&    \frac{\alpha^2\big(\frac{c_1L}{2}+ \tilde{P}_{11} \big) \mg{+ a}}{1-s }   \sum_{k=0}^K \|w_k\|^2 
+ \frac{\mg{(1 +   \frac{4b\nu^2 L^2}{\mu} c_1 + \rho_3^2}) V_{P,c_1}(\xi_0)}{1-s},
\label{ineq-sum-lyapunov}
\eeq 
\mg{provided that $s=  \rho_0^2\mg{+c_1\rho_1^2+\rho_2^2+\rho_3^2}  \mg{+ \frac{4b(\nu^2 + (1+\nu)^2) L^2}{\mu} c_1 }  <1$}. 
If $\tilde{P}_{22} \neq 0$, using Schur components, %we can write
$$ \tilde{P} = \begin{bmatrix} 
	\tilde{P}_{11} - {\tilde{P}_{12}^2}/{\tilde{P}_{22}} & 0\\
	0 & 0
\end{bmatrix} 
+ 
\begin{bmatrix}
{\tilde{P}_{12}^2}/{\tilde{P}_{22}} & \tilde{P}_{12}\\
\tilde{P}_{12} & \tilde{P}_{22}
\end{bmatrix} \succeq 
\begin{bmatrix} 
	\tilde{P}_{11} - {\tilde{P}_{12}^2}/{\tilde{P}_{22}} & 0\\
	0 & 0
\end{bmatrix}. 
$$
Otherwise, if $\tilde{P}_{22}=0$, then we have $\tilde{P}_{12}=\tilde{P}_{21}=0$ as $\tilde{P}\succeq 0$. Therefore, in any case, %we can write 
\beq ~\begin{bmatrix} {x}_k - x_* \\ {x}_{k-1} - x_*  \end{bmatrix}^T P \begin{bmatrix} {x}_k - x_* \\ {x}_{k-1} - x_*  \end{bmatrix}  &\geq& r(\tilde P)  \|{x}_k-x_*\|^2
\geq   \frac{2}{L} r(\tilde P)  \left(  f({x}_k) - f(x_*) \right),
\eeq
where in the last inequality we used the $L$-smoothness of $f$.
%$$ r(\tilde P) = \begin{cases}
% \tilde{P}_{11} - {\tilde{P}_{12}^2}/{\tilde{P}_{22}} & \mbox{if} \quad \tilde{P}_{22}\neq 0, \\
% \tilde{P}_{11} & \mbox{if} \quad \tilde{P}_{22}=  0,
% \end{cases} 
%$$ 
Consequently, 
%\beq  
 $\left( c_1 + \frac{2}{L}r(\tilde{P}) \right) [ f( x_k) - f(x_*) ] \leq V_{P,c_1}({\xi}_k).$  
%\label{ineq-from-lyap-to-fvals}
%\eeq
From the strong convexity of $f$ and the fact that $x_{-1}=x_{0}$, we have also
%\beq 
$V_{P,c_1}(\xi_0) \leq \|\tilde P\| \| \xi_0\|^2 + c_1 (f(x_0) - f(x_*)) \leq
\left( \frac{4 \|\tilde P\|}{\mu} + c_1\right)  (f(x_0) - f(x_*))$.
%\label{ineq-upper-bound-on-lyap}
%\eeq 
Combining \eqref{ineq-sum-lyapunov} with these inequalities on the Lyapunov function $V_{P,c_1}$, we obtain
 \beq \frac{V_{P,c_1}({\xi}_{K+1})}{1-s}  &+& \sum_{k=0}^K \left[ f({x}_k) - f(x_*) \right] \leq  \big(\overline{L}_{GMM}(\alpha,\beta,\nu)\big)^2 \sum_{k=0}^K  \|w_k\|^2 + H(\xi_0), \quad \mbox{with} \quad 
 \label{ineq-gmm-nonasymptotic-perf} \\
% \eeq %\eqref{def-hinfty-nonlinear} is satisfied for 
%$$\gamma^2 =\frac{1}{c_1 + \frac{2}{L}r(\tilde{P})}
%\frac{\alpha^2  \big(\frac{c_1L}{2}+ (1+a) \tilde{P}_{11}\big) \mg{+ a}}{1- \big( \rho_0^2\mg{+c_1\rho_1^2+\rho_2^2+\rho_3^2}  \mg{+ \frac{4b(\nu^2 + (1+\nu)^2) L^2}{\mu} c_1 } \big) } , %\quad h_0= \frac{1}{c_1 + \frac{2}{L}r(\tilde{P}) } \left( \frac{4 \|\tilde P\|}{\mu} + c_1\right), 
%$$
%with the choice of 
%\beq 
H(\xi_0)&=& \frac{\mg{(1 +   \frac{4b\nu^2 L^2}{\mu} c_1 + \rho_3^2}) }{1-s} \left( \frac{4 \|\tilde P\|}{\mu} + c_1\right) \frac{1}{c_1 + \frac{2}{L}r(\tilde{P})} (f(x_0) - f(x_*)).
\label{def-Hxi0-GMM}
\eeq
Then, letting $K\to \infty$ leads to $L_{2,*} \leq \overline{L}_{GMM}(\alpha,\beta,\nu)$ by the definition. %\myqed
%\beq c_{a,\rho} = &&\inf_{ \tilde{P}\succeq 0 }  {\frac{\alpha^2}{1-\rho^2} (\frac{L}{2}+ (1+a) \tilde{P}_{11})}\\
%&&\mbox{subject to } M_a \succeq 0
%\eeq
%Then, $H_\infty^2 \leq c_{a,\rho}$ for any $a\geq 0$ and $\rho \in [0,1)$. This completes the proof.
\end{proof}
%\endproof
%\vspace{-0.2in}
\begin{coro}[$\mathcal{O}(1/K)$ ergodic rate for inexact GMM]\label{coro-ergodic-rates-gmm} In the setting of Thm. \ref{thm-hinfty-bound}, when the upper bound \eqref{def-Hinfty-ub} holds, $\bar{x}_K = \frac{1}{K+1} \sum_{j=0}^K x_K$ satisfies
$ [f(\bar x_K) - f(x_*) ] \leq \frac{1}{K+1} \sum_{j=0}^K [f(x_j) - f(x_*)] \leq \big( \overline{L}_{GMM}(\alpha,\beta,\nu)\big)^2 \frac{\sum_{j=0}^{K} \|w_j\|^2}{K+1} + \frac{H(\xi_0)}{K+1} -  \frac{1}{1-s}\frac{V_{P,c_1}({\xi}_{K+1})}{K+1 } $
where $H(\xi_0)$ is as in \eqref{def-Hxi0-GMM}.
%with $H(x_0) = \frac{2}{\sqrt{\alpha\mu}} 
% \bigg(1 +  \frac{\sqrt{\alpha\mu} (1- \sqrt{\alpha\mu}) }{4 \mg{(5-2\sqrt{\alpha\mu} + \alpha\mu)} } + \frac{\sqrt{\alpha\mu} (1-\sqrt{\alpha\mu})^2}{4(4+ (1-\sqrt{\alpha\mu})^2)} \bigg) V_0$ for any $k\geq 0$. % and $w \in \ell_2(\mathbb{R}^d)$.
 %\sum_{k=1}^K f(x_k) - f(x_0)%By the convexity, ergodic iterate% $\bar{x}_K = \sum_{k=1}^K \frac{1}{K} x_k$.
%While \eqref{ineq-nag-cum-lyap-bound} has an asymptotic nature, it is also straightforward to obtain non-asymptotic bounds on the cumulative sums $\sum_{k=0}^K \[f(x_k) - f(x_*)\]$ in terms of the noise sequence $\{w_k\}_{k=0}^K$ from the recurrence \eqref{ineq-to-sum-side-by-side-nag}.
\end{coro}
\proof {\textbf{Proof}.}  This is a direct consequence of the convexity of $f$ and \eqref{ineq-gmm-nonasymptotic-perf}. %\myqed
\endproof
%\vspace{-0.2in}
\begin{rema}[MI approach recovers our explicit bounds]\label{rema-MI-explicit-bounds} Thm. \ref{thm-hinfty-bound} can recover the same explicit bounds we obtained for GD and NAG if we choose the parameters of the MI in particular ways: For GD, taking
$\tilde{P}_{11} = 1$, $\tilde{P}_{22}=\tilde{P}_{12}=c_1 = b=\rho_1 = \rho_2 = \rho_3 = 0$ and $a = \alpha^2 (\frac{1}{\rho_{GD}}-1)$ generates the $\ell_2$ gain upper bound \eqref{ineq-gd-hinfty-bd1} whereas $c_1 = 1$, $\tilde{P}=0$, $c_0=b=\rho_1 = \rho_2 = \rho_3 = 0$, $a=\frac{\alpha(|1- L\alpha|)}{2r}$ and $\rho_0^2 = 1-2\mu \alpha (1 - \frac{L\alpha}{2}) + \alpha \mu |1-\alpha L| r$ with $r=1$ when $\alpha \leq \frac{1}{L}$ and $r = \frac{2-\alpha L}{\alpha L}$ for $\alpha \in (1/L,2/L)$ generates the other upper bounds for $L_{2,*}$ we obtained for GD in the proof of Prop. \ref{prop-hinf-gd-bound}. Similarly for NAG with $\beta = \nu =  \frac{1-\sqrt{\alpha\mu}}{1+\sqrt{\alpha\mu}}$ for $\alpha \in (0,1/L]$, if we take $\rho_1^2 = \frac{2\alpha}{4s_1(1+\sqrt{\alpha\mu})^2}$, $\rho_2^2 = \frac{2}{\mu s_1}
  \frac{\alpha^2\mu^2 + 2\alpha \mu + \alpha\mu (1-\sqrt{\alpha \mu})}{4(1+\sqrt{\alpha\mu})^2}$, $\rho_3^2 = \frac{2}{\mu s_1} \frac{\alpha\mu (1-\sqrt{\alpha \mu})}{4(1+\sqrt{\alpha\mu})^2}$, $c_1 = 1, c_0 = 0$, $\rho_0^2 = \rho_{NAG}^2 = 1-\sqrt{\alpha\mu}$, 
$a=  s_1 + \frac{L\alpha^2}{2} s_2$ and $b = \frac{L\alpha^2}{2s_2}$ with $s_i = \sqrt{\alpha}\hat{s}_i$ for $i=1,2$ where $\hat{s}_i$ is as in \eqref{def-hat-r1-r2} and $\tilde{P}$ according to \eqref{def-lyap-agd}, then Thm. \ref{thm-hinfty-bound} yields the same explicit upper bound  obtained for NAG in Thm. \ref{prop-hinf-agd-bound}. For these choice of parameters for GD and NAG, the fact that the MI ($\tilde{M}_4\succeq 0$) holds follows from straightforward but tedious computations; we provide the details in Appendix \ref{sec-appendix-online-companion}. \looseness=-1
\end{rema}
\begin{rema}[Non-square-summable errors]\label{remark-non-square-summable}
In case $w_k\not\in\ell_2(\mathbb{R}^d)$ but has a finite power, i.e if $\mathcal{P}_\omega:=\sup_k \frac{\sum_{j=0}^{k-1} \|w_j\|^2}{k+1} < \infty$, then Cor. \ref{coro-ergodic-rates-gmm} implies that $$ f(\bar x_K) - f(x_*) \leq \overline{L}_{GMM}(\alpha,\beta,\nu) \mathcal{P}_\omega + \mathcal{O}(1/K).$$ For example, if $\|w_k\| \leq \delta ~\forall k$, we can take $\mathcal{P}_\omega= \delta^2$.  In other words, the $\ell_2$ gain and its estimate $\overline{L}_{GMM}(\alpha,\beta,\nu)$ is also relevant to performance for bounded errors or errors with a finite power, in which case the convergence happens only in a neighborhood of the solution due to persistence of errors. Similar conclusions can be drawn for NAG and GD (see Coro. \ref{coro-ergodic-rates} and proof of Prop. \ref{prop-hinf-gd-bound}).\looseness=-1%However, when stepsize is larger; i.e. when $\alpha = 1/L$; then our bounds grows linearly with $\sqrt{\kappa}$ whereas in the quadratic case, our bounds for $L_{2,*}$ based on the $H_\infty$ norm were independent of $\kappa$ for NAG. 
\end{rema} 

\begin{rema}[GMM rate with exact gradients]\label{remark-conv-GMM} By setting the noise term equal to zero, i.e. $w_k = 0$, we can obtain convergence guarantees for (deterministic) GMM methods without noise. In particular, if we set $c_1= 1$ and \mg{$\|w\|_{\ell_2} = \rho_1=\rho_2=\rho_3 = a=b=0$} in \eqref{ineq-evolution-Lyap-with-noise}, %\eqref{ineq-evolution-Lyap-with-noise}, 
we obtain 
%\beq 
$V_{P,1}({\xi}_{k+1})   \leq \rho_0^2 V_{P,1}({\xi}_{k})$
%\label{ineq-Vpc-last}
%\eeq
provided that $\mathcal{S}_{\rho_0}(\tilde{P})$ (the leading $3\times 3$ principal submatrix of $\tilde{M}_2
   +c_1 \tilde{M}_1 +c_0 \tilde{M_0}$) is positive semi-definite for $c_1=1$.
  % \eeq
%\beq \footnotesize ~\qquad\qquad \mathcal{S}_{\rho_0}(\tilde{P}) : = \begin{bmatrix}   -\tilde{A}^T \tilde{P} \tilde{A} + \rho_0^2 \tilde{P} + c_1\tilde{Y}_1   + c_0 2mL\tilde{C}^T\tilde{C}   &- \tilde{A}^T \tilde{P} \tilde{B}  + c_1\tilde{Y}_2   -c_0 (m+L)\tilde{C}^T \\ - \tilde{B}^T \tilde{P} \tilde{A} + c_1\tilde{Y}_2^T  -c_0(m+L)\tilde{C}   & -  \tilde{B}^T \tilde{P} \tilde{B} + c_1\tilde{Y}_3 +  2c_0  \\
%\end{bmatrix}
%%\label{def-S-rho}
%\eeq
%is positive semi-definite. 
Since, $V_{P,1}({\xi}_k) \geq f(x_k) - f(x_*)$, this %\eqref{ineq-Vpc-last}
implies the linear convergence rate result $f({x}_k) - f(x_*) \leq \rho_0^{2k} V_{P,1}({\xi}_0)$. 
In the special case when $c_0 = 0$, this would recover the deterministic convergence rate analysis provided in \cite{hu2017dissipativity} for deterministic GMM methods. Our analysis here supports more general choices of $c_0$ in the sense that we allow $c_0\geq 0$. %as long as it is non-negative.
\end{rema}
\subsection{Trading rate with robustness.}\label{subsec-trading-rate-robustness} We next describe how we can design the parameters to trade rate with robustness in a systematic fashion. To illustrate the ideas, we first focus on the NAG method (where $\beta=\nu$). Given $\rho \in [0,1)$ fixed, we next consider all choice of parameters $(\alpha,\beta)$ so that NAG can be certified to converge with rate $\rho^2$ without gradient errors. That is, we consider the following parameters that are guaranteed to satisfy 
$f(x_k)-f(x_*)\leq \rho^{2k}V_{P,1}(\xi_0)$ when $w_k = 0$ for all $k$:
	\beq \mathcal{P}({\rho}^2):= \{ (\alpha,\beta) ~:~ \alpha\geq 0, \beta\geq 0, \exists \tilde{P}\succeq 0 \mbox{ such that }\mathcal{S}_\rho(\tilde{P}) \succeq 0 \},
	\label{def-set-params-rate}
	\eeq
where\footnote{Note that if we set $w_k = 0$ $\forall k$ in the proof of Prop. \ref{prop-hinf-agd-bound}, the linear convergence rate $\rho^2$ can be certified if \eqref{def-set-params-rate} holds.} $\mathcal{S}_\rho(\tilde P)$ is defined by \eqref{def-S-rho}. In particular for NAG with $\alpha \in (0,1/L]$, $\beta = 1-\sqrt{\alpha\mu}$ and $\tilde{P}$ as in \eqref{def-lyap-agd}, we have 
$\mathcal{S}_\rho(\tilde{P}) \succeq 0$ when
 %$\tilde{P}$ is as in \eqref{def-lyap-agd}  we have $\mathcal{S}_\rho(\tilde{P}) \succeq 0$ when $\tilde{P}$ is as in \eqref{def-lyap-agd} with the rate
  $\rho^2 = 1-\sqrt{\alpha\mu}$ 
  (see the proof of Prop. \ref{prop-hinf-agd-bound}). For $\alpha=1/L$, we obtain the fastest such certified rate $\rho^2 = \rho^2_{NAG}=1-\sqrt{\frac{\mu}{L}}=1-\frac{1}{\sqrt{\kappa}}$. A natural way to trade robustness with rate is to find $(\alpha,\beta)$ that lead to the best robustness bound if we allow the convergence rate to be slower than the baseline rate $\rho_{NAG}^2$ by a certain percentage, i.e. solve % at a given rate, choose the parameters as a solution to
\beq  \mbox{Opt}_\varepsilon := \min_{(\alpha,\beta)\in \mathcal{P}{({\rho^2})}}  \overline{L}_{GMM}(\alpha,\beta,\beta) \quad \mbox{such that}\quad \rho^2 = \rho_{NAG}^2 (1+\varepsilon),
\label{opt-pbm-tradeoff}
\eeq
where $ \overline{L}_{GMM}(\alpha,\beta,\nu)$ is defined by \eqref{def-Hinfty-ub} %$\rho_{NAG}^2 = 1-\sqrt{\alpha\mu}$ is the benchmark accelerated rate and $\varepsilon\geq 0$ 
and $\varepsilon \geq 0$ is the trade-off parameter that represents the percentage rate degradation compared to the fastest rate $\rho^2_{NAG}$. By compactness of the constraint set, a minimizer $ (\alpha_\varepsilon,\beta_\varepsilon)$ of \eqref{opt-pbm-tradeoff} exists. Clearly, the interesting case is when $\varepsilon$ is not too large so that $\rho < 1$. For given $\varepsilon\geq 0$ fixed, \eqref{opt-pbm-tradeoff} is a non-convex problem. However, it is a small-dimensional problem and we can approximate its solutions by a simple grid search approach as follows: First, we grid the parameter space $(\alpha, \beta) \in (0, \frac{2}{L}) \times [0, \frac{\kappa}{\kappa-1})$, and\footnote{%We recall we assume $\kappa>1$ throughout this paper. 
This range for $\beta$ is motivated by the fact that for $\alpha=1/L$, any $\beta\in [0, \frac{\kappa}{\kappa-1})$ results in $\rho(A_Q)<1$ and also satisfies $H_\infty<\infty$ by Thm. \ref{thm-h-inf}, but by the same theorem we have also $c_\mu \to 0$ and $H_\infty\to \infty$ as $\beta \to \frac{\kappa}{\kappa-1}$ where $\kappa>1$.} for each element of the grid, we check whether it lies in the constraint set $\mathcal{P}(\rho^2)$ by solving the SDP feasability problem given in \eqref{def-set-params-rate}. If that is the case, we generate an upper bound on $L_{2,*}$ based on Thm. \ref{thm-hinfty-bound}, by another grid search over the parameters %($\rho_0,\rho_1,\rho_2,\rho_3,a,b,c_0,c_1$) 
of the $4\times 4$ matrix inequality $\tilde{M}_4\succeq 0$ from Thm. \ref{thm-hinfty-bound}. %\footnote{Note that $\mathcal{S}_\rho(\tilde P)$ is the $3\times 3$ leading principal matrix of $\tilde{M}_4$ if $\rho_0=\rho$ and $\rho_1=\rho_2=\rho_3=a=0$, by passing to the Schur complements $\tilde{M}_4\succeq 0$ iff $\mathcal{S}_\rho(\tilde P) - a $
%when $w_k = 0$, both $\mathcal{S}_\rho(\tilde P)\succeq 0$ and $\tilde{M}_4\succeq 0$ holds for the same $\til de{P}$.} 
In addition, we can use the fact that
%we can combine the grid search approach with explicit choice of parameters that %we can generate explicit approximations to the optimal value $Opt_\varepsilon$
the specific choice of parameters $\beta(\alpha) =\frac{1-\sqrt{\alpha\mu}}{1+\sqrt{\alpha\mu}}$ and $\alpha\in (0,1/L)$, which arises commonly \cite{aybat2018robust, aybat2019universally}, lie in $\mathcal{P}(\rho^2)$ for 
$\rho^2=1-\sqrt{\alpha\mu}$ and this is decreasing in $\alpha$, while the $L_{2,*}$ bound \eqref{ineq-hinf-bound-agd-by-hand} is increasing in $\alpha$. Therefore, in this particular momentum parametrization, the best upper bound for $L_{2,*}$ satisfying the rate constraints will be obtained for the stepsize $\tilde{\alpha}_\varepsilon$ with rate $1-\sqrt{{\tilde \alpha}_\varepsilon \mu}=\rho^2_{NAG}(1+\varepsilon)$, i.e. when 
\beq {\tilde\alpha}_\varepsilon = \frac{1}{\mu}\left(1 - (1+\varepsilon)(1-\frac{1}{\sqrt{\kappa}})\right)^2 \in (0,\frac{1}{L}], \quad {\tilde\beta}_\varepsilon = \frac{1-\sqrt{\tilde{\alpha}_\varepsilon\mu}}{1+\sqrt{{\tilde\alpha}_\varepsilon\mu}}, %\implies L_{2,*}^{ub}\leq \bar{L}(\alpha_\varepsilon), 
\label{def-explicit-tradeoff-NAG-params}
\eeq 
which leads to the explicit bound
$ \mbox{Opt}_\varepsilon  \leq \overline{L}_{NAG}(\tilde{\alpha}_\varepsilon)$
where $\overline{L}_{NAG}(\cdot)$ is as in \eqref{ineq-hinf-bound-agd-by-hand}. For computing the approximate solutions of \eqref{opt-pbm-tradeoff}, we compare the best parameters on the grid with the explicit choice of parameters given in \eqref{def-explicit-tradeoff-NAG-params} and choose the best one. This methodology can be applied to any strongly convex problem to design parameters to trade-off worst-case robustness with rates.  In the next section, we will provide a numerical example that follows this methodology. For GMM with more general parameters, the solutions to the analogue of the rate problem \eqref{def-set-params-rate} is also known for some parameterization of parameters \cite{can2022entropic}. Therefore, we can follow a similar approach that combines grid search with the explicit choice of parameters that satisfies the rate constraints.

%\vspace{-0.15in}
\section{Numerical Experiments.}\label{sec-num-experiments}
\begin{figure}[h!]
  \centering
    \includegraphics[height=1.9in, width=2.5in]{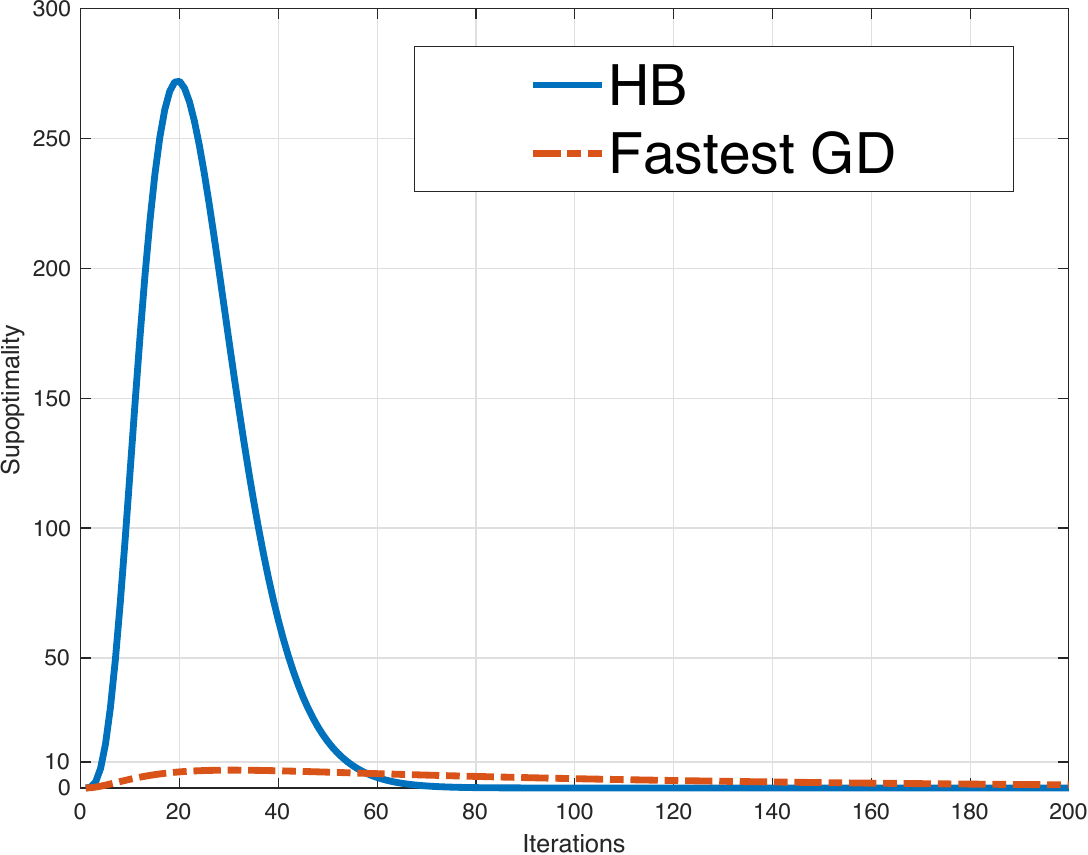}  
    \includegraphics[height=1.9in, width=2.5in]{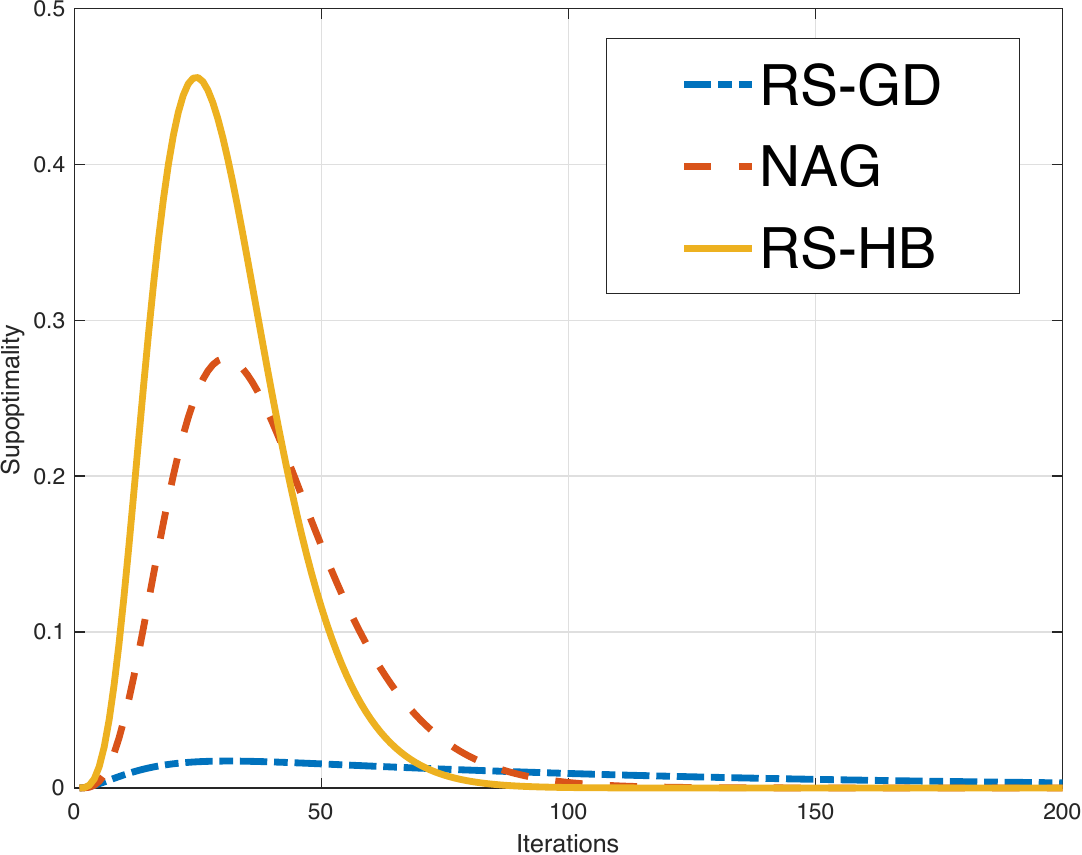}
  %  \includegraphics[height=1.672in, width=2.2in]{str_cvx_cropped.pdf}  
 % \vspace{-0.01in}
  \caption{\label{fig-quad-robustness} (Left) HB method with standard parameters, and GD with the fastest rate obtained for stepsize $\alpha = \frac{2}{L+\mu}$. (Right) RS-GD, RS-HB and NAG methods that achieve the best robustness level $L_{2,*}= H_\infty = \frac{1}{\sqrt{2\mu}}$ for quadratics.\vspace{-0.35in}}
\end{figure}

In our first set of experiments, we consider the quadratic function $f(x) = x^T Q x + \frac{c}{2}\|x\|^2 $
previously considered in \cite{Hardt-blog} for studying the robustness of accelerated methods to stochastic noise where $d=100$ and $Q$ is the Laplacian of a cyclic graph
and $ c= 0.01$ is a regularization parameter. Here, $f \in \Cml$ with $\mu=c=0.01$ and $L = 4 + c = 4.01$ with a condition number $\kappa = L/\mu = 401$ and minimum at $x_* = 0$. In Fig. \ref{fig-quad-robustness}, we compare the performance of HB, GD, AG, RS-GD and RS-HB methods starting from the optimum, i.e. we take $x_0 = x_{-1} = x_*$. For each method, we construct the worst-case noise according to the formula \eqref{eq-worst-case-noise} with $h=0.1$, normalized to have $\|w\|_{\ell_2}= 1$. On the left-hand side of Fig. \ref{fig-quad-robustness}, we compared HB and ``Fastest GD" (i.e. GD with $\alpha = \frac{2}{L+\mu}$) by plotting suboptimality vs. iterations; both methods admit a robustness
of $ L_{2,*}^2 = H_\infty^2 = \sum_{k\geq 0} \left( f(x_k) - f(x_*) \right) = \frac{\kappa}{2\mu} $  (see also Table \ref{table}). %where $\sum_{k\geq 0} \left( f(x_k) - f(x_*) \right)$ is the cumulative suboptimality over iterates.
Therefore, for both HB and ``Fastest GD", the cumulative suboptimality $ \sum_{k\geq 0} \left( f(x_k) - f(x_*) \right)$ is the same but we observe this is achieved in different ways:  %a%nd a
%the way to achieve this is different in the sense that the peak value of suboptimality $f(x_k)- f(x_*)$ and it decays for large $k$ are different. %To simulate the worst-case noise for each method, we use Prop. \ref{prop-worst-case-noise}. 
%We observe from Fig. \ref{fig-quad-robustness}  that 
HB gets significantly larger suboptimality values in the beginning % (that exceed the level 250) 
with a faster decay towards the end, whereas the suboptimality of ``Fastest GD" has a (maximum) peak that is smaller (by a factor of at least $25$) and it decays more slowly compared to HB. On the right panel of Fig. \ref{fig-quad-robustness}, we plot the most robust methods RS-GD, NAG and RS-HB from Table \ref{table} which have $L_{2,*}^2 = H_\infty^2 = \sum_{k\geq 0} \left( f(x_k) - f(x_*) \right)= \frac{1}{2\mu}$. In particular, we see that the cumulative suboptimality of these robust methods will be $\kappa$ times smaller compared to (standard) HB and ``Fastest GD"; we observe that the peak suboptimality values are also much smaller (smaller than 0.5) for these methods. We note that cumulative suboptimality of RS-GD, NAG and RS-HB will be the same; therefore it is not possible for one of these methods to have suboptimality values always strictly less than the others (uniformly) over the iterations. Here, %but achieving this cumulative sum will be different for each method. 
RS-HB has a faster rate $\rho$ than the other methods (see also Prop. \ref{prop-robust-hb}), as a consequence it exhibits faster suboptimality decay towards the end, and as such RS-HB peak suboptimality values exceed the others. %Basically, it is not possible for one method to have suboptimality values always less than the other or more than the others uniformly over all the iterations; because the cumulative suboptimality is the same in this case and only the rates $\rho$ differ. 
Comparing the left and right panels of Fig. \ref{fig-quad-robustness}, the results show that RS-HB and RS-GD are indeed robust versions of GD and HB methods: %Among the robust methods; RS-GD has a smaller peak (for suboptimality) but admits a slower convergence towards the end, whereas RS-HB is faster towards the end but with a higher peak.
% as suggested by our theoretical results. 
%if the peak suboptimality needs to be kept small, RS-GD may be a method of choice; whereas if the goal is to minimize the convergence behavior towards the end when RS-HB would be an efficient method. 
The peak values of RS-HB and RS-GD are much smaller compared to HB and GD, and the cumulative suboptimality improvement obtained by RS-HB and RS-GD is proportional to the condition number; i.e. we can get significant gains  in performance for ill-conditioned problems. %he difference will get mor%e can have significant improvements. %for problems with a larger condition number; the improvement will be larger. 
%\begin{figure}[h!]
%  \centering
%    \includegraphics[height=1.9in, width=2.5in]{str_cvx_cropped.pdf}  
 % \caption{Trading robustness with rate specified with different $\varepsilon$ levels\label{experiment-str-cvx}}
%\end{figure}
\noindent
\begin{minipage}[t]{0.55\textwidth}
%% PUT TEXT HERE
Next, we  consider the following strongly convex smooth objective from \cite{scoy-tmm-ieee} which is similar to the heavy-ball counter-example studied in \cite{lessard2016analysis}:\looseness=-1
\beq f(x) &=& \sum_{i=1}^p g(a_i^Tx - b_i) + \frac{\mu}{2}\|x\|^2 \quad \mbox{with} \\ g(x) &=& \begin{cases} \frac{1}{2}x^2 e^{-r/x} & x>0, \\
																0 & x\leq 0,
\end{cases} 
\label{obj}
\eeq
where $\mu>0$, $p\geq 1$, with $a_i \in \mathbb{R}^d$ for every $i=1,2,\dots,p$ so that the matrix $A = \begin{bmatrix} a_1, a_2, \dots, a_p \end{bmatrix} \in \mathbb{R}^{d\times p}$ and $b\in\mathbb{R}^p$. Here, we generate the entries of the data matrix $A$ and the vector $b$ randomly (from the uniform distribution on $[-\frac{1}{2}, \frac{1}{2}]$) and then the matrix $A$ is scaled so that $\|A\|=\sqrt{L-\mu}$ for some $L>\mu$. With this scaling, we have $f\in\Cml$ and $f$ admits continuous derivatives
\end{minipage}\hfill
\begin{minipage}[t]{0.4\textwidth}
\hspace{1in}
  \centering\raisebox{\dimexpr \topskip-\height}{%
  \includegraphics[width=\textwidth]{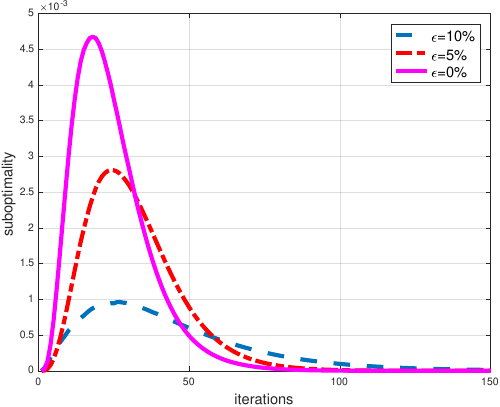}} \hspace{1in}\\
  \captionof{figure}{Trading robustness with rate specified with different $\varepsilon$ levels\label{experiment-str-cvx}}
  \label{fig1}
\end{minipage}\hfill
\\of all order. We take $\mu=\frac{1}{2}, L=30$, $d=100$, $p=5$ in \eqref{obj} and consider the problem \eqref{opt-pbm-tradeoff} for trading rate with robustness as $\varepsilon$ is varied in the set $\{0\%, 5\%, 10\%\}$. For each $\varepsilon$, we approximate the solutions $(\alpha_{\varepsilon_*}, \beta_{\varepsilon_*})$ of \eqref{opt-pbm-tradeoff} using the methodology from Sec. \ref{subsec-trading-rate-robustness} where larger $\varepsilon$ is associated with increased robustness at the expense of a slower (certified) rate.
%For more general choices of $(\alpha,\beta)$; we do a grid search. More specifically, by a grid search approach over the parameters $(\alpha,\beta)$ as well as the parameters of the $4\times 4$ matrix inequality \eqref{def-Hinfty-ub} together with the explicit bound , we first approximate the solutions $(\alpha_\varepsilon,\beta_\varepsilon)$. 
We then run the NAG method with these parameters in Fig. \ref{experiment-str-cvx} starting from the optimum $x_0 = x_{-1}=x_*$ subject to gradient noise. For given trade-off parameter $\varepsilon$, (almost) the worst-case noise sequence $\{w_k^{(h)}\}$ depends on the parameters $(\alpha_\varepsilon,\beta_\varepsilon)$; therefore it will depend on the value of $\varepsilon$ as well. In the strongly convex case, we do not have an explicit characterization of the almost worst-case noise sequence %that corresponds to the choice of parameters $(\alpha_\varepsilon, \beta_\varepsilon)$ 
as the noise sequence \eqref{eq-worst-case-noise} constructed in Prop. \ref{prop-worst-case-noise} was tailored to quadratic objectives. Therefore, to measure robustness, we considered a class of error sequences that generalizes \eqref{eq-worst-case-noise} where at every step $k$, 
we make a quadratic approximation
 $ f_{Q_k}(x) = (x-x_k)^T Q_k (x-x_k)$ with $Q_k := \nabla^2 f(x_k)$ to $f(x)$
around $x =  x_k$ and we take the worst-case noise sequence\footnote{More specifically, we let $\mu_k$ and $L_k$ to be the smallest and largest eigenvalue of $Q_k$,  we take $w_k^{(h)} =\sqrt{h (2-h)} (1-h)^k \cos(\omega_* k) u_{k}$ where $u_{k}$ is a normalized eigenvector that corresponds to the eigenvalue $\lambda_*\in \{\mu_k, L_k\}$ where $\lambda_*$ is a maximizer of \eqref{hinfty-quad} for $\mu=\mu_k$ and $L=L_k$.} that corresponds to $f_{Q_k}(x)$ where we also
%take 
% We take $\mu=\frac{1}{2}, L=30$, $d=100$, $p=5$ in \eqref{obj} and consider the problem \eqref{opt-pbm-tradeoff} as $\varepsilon$ is varied in the set $\{0\%, 5\%, 10\%\}$. 
%\beq w_{k,\varepsilon}^{(h)} =\sqrt{\delta (2-\delta)} (1-\delta)^k \cos{\omega_*(k) k} u_* \quad \mbox{for} \quad k\geq 0, 
%\label{eq-str-cvx-worst-case-noise}
%\eeq 
%and choose the worst-case noise corresponding to the quadratic $ f_{Q_k}$ at step $k$,
tune the decay parameter $h$.\footnote{To simulate worst-case noise, we also considered alternative noise structures where $w_k^{(h)}$ is generated randomly over i.i.d. trials from the uniform distribution and then the realization that maximizes $f(x_{k+1})-f(x_*)$ is selected.} For each $\varepsilon$ value, we plot the response with respect to the approximated worst-case errors. Fig. \ref{experiment-str-cvx} shows the results where we can see the effect of the trade-off parameter $\varepsilon$ on the path of the iterates. When $\varepsilon$ is the smallest, NAG is the fastest asymptotically as expected but then there is a larger peak of suboptimality at the (initial) transient phase. When $\varepsilon$ gets larger, we have a larger set of parameters to satisfy the same rate; therefore we can improve the robustness behavior by optimizing this set of parameters, trading the rate with robustness. %by trading the rate with robustness. %For a given budget $K$
% where 
 %with initialization $x_0 = x_{-1} = $ where $\begin{bmatrix} v^H & v^H \end{bmatrix}^H$ is a right eigenvector of $A_Q + B \Delta_* T$ corresponding to the eigenvalue $e^{i\omega_*}$ with $\Delta_*  = \frac{1}{H_\infty}u_* v_*^H$, 
%{\begin{small}$$
%u_* = \begin{cases} u_1 & \mbox{if } \lambda_* = \mu \\
%                      u_d & \mbox{if } \lambda_* = L     
%\end{cases}, ~
%$$ 
%\end{small}}
%\mtodo{the fig 10 section above needs further editing}
%\vspace{-0.1in}
\section{Conclusion.} %$\ell_2$ gains were never used to design parameters of momentum algorithms for achieving (approximate) Pareto-optimal rate and robustness performance despite the popularity of such approaches for designing robust control systems

%While (induced) $\ell_2$ gains from the noise input to output of interest (such as suboptimality) are popular to assess the robustness of dynamical systems in control theory, they have not been applied to optimization algorithms to design the parameters of first-order methods to make them robust. 

In this work, we considered generalized momentum methods (GMM) for strong convex smooth minimization
%he trade-offs between the convergence rate and robustness to deterministic gradient errors when designing the parameters of a first-order algorithm for minimizing a strongly convex smooth function. 
%We focus on a general class of momentum methods called 
that include GD, NAG, TMM and HB methods as special cases depending on how the GMM parameters are chosen. We consider the inexact gradient setting where gradients admit worst-case deterministic errors. We proposed a framework to design the parameters for trading off the convergence rate and robustness in a computationally tractable and systematic manner, where robustness is measured in terms of the cumulative suboptimality over iterations normalized properly by the $\ell_2$ norm of the gradient errors. We interpret this measure as the $\ell_2$ gain of an associated dynamical system corresponding to GMM iterations. %For quadratic objec
%convergence rates with robustness systematically in a computationally tractable fashion, where robustness is 
  %We introduced a robustness 
%the study of momentum-based optimization algorithms. 
%We propose the \ell_2 gain from gradient errors to suboptimality
%We propose to measure robustness in terms of the cumulative suboptimality over iterations normalized 
%by the $\ell_2$ norm of the gradient errors. 
%We interpret this robustness measure, as the $\ell_2$ gain (from input to output) of the (dynamical) system corresponding to GMM iterations where the input is the gradient error sequence and the output is a weighted distance of the iterates to the optimum. 
%argue that this robustness measure corresponds to the $\ell_2$ gain of the dynamical system corresponding to GMM iterations that admits gradient errors as the input when the output is set to a weighted distance to the optimum. 
For strong convex quadratic objectives, $\ell_2$ gain coincides with the $H_\infty$ norm for which we provide an explicit formula, characterizing all the parameters that achieve the best robustness level and constructing worst-case gradient error sequences explicitly as a function of the parameters. We also study the Pareto-optimal boundary between the convergence rate and robustness and find that %and construct worst-case gradient errors. 
%Our results show that 
HB is less robust than NAG with standard parameters despite being faster. As a remedy, we propose the robustly stable heavy ball method that is faster than NAG while being at the most robust level. In addition, we propose the robustly stable gradient descent that is the fastest version of GD while being at the best robustness level.  %We also study the Pareto-optimal boundary between the convergence rate and robustness and show that RS-GD and RS-HB lie on the Pareto-optimal boundary.  %By characterizing 
%%we can compute the $\ell_2$ gain explicitly, leveraging its representation as the $H_\infty$ norm of the GMM system in the frequency domain. %, which we refer to as $L_{2,*}$ is equivalent to the $H_\infty$ norm.
%We consider robustness of GMM subject to gradient errors that are square-summable. We measure robustness by the cumulative suboptimality normalized by the $\ell_2$ norm of the input gradient errors, which
% is equivalent to $\ell_2$ gain of a transformed system 
%for quantifying robustness of GMM to square summable worst-case gradient noise which measures the effect of gradient errors on the cumulative suboptimality over the iterations. F
%or quadratic objectives,
% this robustness measure coincides with 
% the $H_\infty$ norm of the transfer function of the GMM system 
%We also provide a lower bound on it, characterizing all the parameter choices that correspond to the best rbustness level which allows us to compare robustness of GD, NAG, HB and TMM. We
%This representation also allows us to construct gradient error sequences that correspond to the worst-case behavior. 
By explicit characterizations of the complex stability radius and real stability radius of the GMM system (which is related to $\ell_2$ gains), we also obtain new results about robustness of GMM to multiplicative gradient errors. Finally, we discuss how our analysis can be extended to general strongly convex smooth objectives where we provide non-asymptotic rate results for inexact GMM methods and derive bounds on the %induced 
$\ell_2$ gain. %where we can choose the GMM parameters to systematically trade off the rate and robustness in a computationally tractable framework. 
\section*{Acknowledgments.}
This paper is dedicated to the 70th birthday of Michael Overton (New York University), who has made key contributions to the computation of various robustness measures for linear dynamical systems including the stability radius and the $H_\infty$ norm. The author would also like to acknowledge his former student Bugra Can for pointing out a useful lemma for the inversion of block matrices and for helpful discussions about the manuscript. M.G. acknowledges partial support by the grants ONR N00014-21-1-2244 and NSF DMS-2053485.
\bibliographystyle{alpha}
\bibliography{hinfty_momentum}
% Appendix here
% Options are (1) APPENDIX (with or without general title) or 
%             (2) APPENDICES (if it has more than one unrelated sections)
% Outcomment the appropriate case if necessary
%
% \begin{APPENDIX}{<Title of the Appendix>}
% \end{APPENDIX}
%
%   or 
%
%\newpage
\appendix
% \section{<Title of Section A>}
% \section{<Title of Section B>}
% etc

\section{\label{app-proof-thm-hinf}Proof of Theorem {\ref{thm-h-inf}}.}

We will compute the $H_\infty$ norm explicitly based on the representation \eqref{eq-h-inf-tf} which says $H_\infty = \max_{\omega\in\mathbb{R}}\|G(e^{i\omega})\|$ and use the fact that $L_{2,*}=H_\infty$. Recall from Sec. \ref{sec-quad} that $Q$ has the eigenvalue decomposition $Q = U \Lambda U^T$ where $U$ is orthonormal. Therefore, the matrix $V=\Diag(U,U) \in \mathbb{R}^{2d\times 2d}$ with diagonal blocks $U$
%\beq\label{def-V-orthonormal} V = \begin{pmatrix}
%    U & 0_d \\
%    0_d & U
%\end{pmatrix}.
%\eeq
is also orthonormal and we have
  \beq A_\Lambda := V^T A_Q V =  \begin{bmatrix} 
(1+\beta)I_d -\alpha(1+\nu)\Lambda & -\beta I_d +\alpha\nu \Lambda \\
I_d & 0_d 
\end{bmatrix},
\label{eq-A-lambda}
\eeq
where $A_Q$ is as in \eqref{def-AQ}. Therefore, $A_Q = V A_\Lambda V^T$. Consequently, plugging $T$ and $B$ into \eqref{def-transfer-matrix} (where $T$ is as in \eqref{def-T} and $B=\tilde B\otimes I_d$ with $\tilde B$ as in \eqref{def: system mat for TMM}) along with this identity, we obtain
\beq 
G(z) 
%&=& \frac{1}{\sqrt{2}}\begin{bmatrix} 
%\Lambda^{1/2}U^T & 0_d 
%\end{bmatrix}   V (zI_{2d} - A_\Lambda)^{-1} V^T B  \\
%&=& \frac{1}{\sqrt{2}
%\begin{bmatrix} 
%\Lambda^{1/2}U^T & 0_d 
%\end{bmatrix}   
%\begin{bmatrix} U & 0_d \\
%0_d & U
%\end{bmatrix}
%(zI_{2d} - A_\Lambda)^{-1} 
%\begin{bmatrix} U^T & 0_d \\
%0_d & U^T
%\end{bmatrix} \begin{bmatrix} -\alpha I_d \\ 0_d \end{bmatrix}  \nonumber \\
&=&
 -\frac{\alpha}{\sqrt{2}}
\begin{bmatrix} \Lambda^{1/2} & 0_d 
\end{bmatrix}  (zI_{2d} - A_\Lambda)^{-1} \begin{bmatrix} U^T \\ 0_d \end{bmatrix} = -\frac{\alpha}{\sqrt{2}} \Lambda^{1/2}   \left[ (zI - A_\Lambda)^{-1}\right]_{(11)}  U^T, \label{eq-trans-func}
%\begin{bmatrix} 
%\Lambda_d^{1/2}, 0_d
%\end{bmatrix}
%(e^{i\omega}I_{2d}-A_{\Lambda})^{-1}\begin{bmatrix} 
%-\alpha I_d\\
%0_d
%\end{bmatrix}\\
%&= 
%-\frac{\alpha}{2} \Lambda_d^{1/2} \left[\left(e^{i\omega}I_{2d}-A_{\Lambda} \right)^{-1}\right]_{11}
\eeq
where $\left[ (zI - A_\Lambda)^{-1}\right]_{(11)}$ denotes the leading principal $d\times d$ submatrix of the $2d\times 2d$ matrix $(zI - A_\Lambda)^{-1}$. Note that when $z$ is on the unit circle, we have $z=ie^{\omega}$ for some $\omega \in \mathbb{R}$ and based on \eqref{eq-A-lambda}, we have
%\begin{eqnarray} 
$e^{i\omega}I_{2d}-A_{\Lambda}= \begin{bmatrix} 
Z_1 & Z_2 \\
Z_3 & Z_4
\end{bmatrix}$,
%\end{eqnarray}
where
%\begin{eqnarray}
$Z_1 = e^{i\omega}I_d -(1+\beta)I_d+\alpha(1+\nu)\Lambda$, $Z_2 = \beta I_d-\alpha\nu \Lambda$, $Z_3 = - I_d$, $Z_4 =  e^{i\omega}I_d$.
%\end{eqnarray}
%{\color{red} $I_d$ may need to be $-I_d$ above?}
Then, by the matrix inversion formula for block matrices \cite[Theorem 2.1]{BlockInverse}, %we have 
%Theorem \ref{thrm: matrix inverse}, 
\begin{eqnarray*} 
\left[\left(e^{i\omega}I_{2d}-A_{\Lambda}\right)^{-1}\right]_{(11)}&=& (Z_1 - Z_2 Z_4^{-1}Z_3)^{-1}\\
%&=&\left[\alpha(1+\nu)\Lambda-\left(\beta I_d-\alpha\nu\Lambda \right)e^{-i\omega}I_d\right]^{-1}\\
%&=& \left[\left(e^{i\omega}-(1+\beta)+\beta e^{-i\omega} \right)I_d +\alpha\Lambda(1+\nu-\nu e^{-i\omega})\right]^{-1}\\
&=& \underset{j=1,..,d}{\Diag} \left[\frac{1}{e^{i\omega}-(1+\beta-\beta e^{-i\omega})+\alpha\lambda_j (1+\nu-\nu e^{-i\omega})} \right].
\end{eqnarray*}
Therefore from \eqref{eq-trans-func},
\begin{equation} 
G(e^{i\omega})=\frac{1}{\sqrt{2}} \underset{j=1,...,d}{\Diag}\left[\frac{-\alpha \sqrt{\lambda_j}}{e^{i\omega}-(1+\beta-\beta e^{-i\omega})+\alpha\lambda_j (1+\nu-\nu e^{-i\omega})} \right] U^T,
\label{eq-transfer-fun-factor}
\end{equation}
so that 
%\beq 
	$\lambda_{\max} \left[G(e^{i \omega})G(e^{i \omega})^*\right]% &=& 
	%\frac{1}{2} \underset{j=1,...,d}{\max}\left[\frac{\alpha^2  \lambda_j}{\| e^{i\omega}-(1+\beta-\beta e^{-i\omega})+\alpha\lambda_j (1+\nu-\nu e^{-i\omega})\|^2 } \right]  \nonumber \\ 
	= \frac{1}{2} \underset{j=1,...,d}{\max}\left[\frac{\alpha^2  \lambda_j}{
	\|e^{2i\omega} + b_{\lambda_j} e^{iw} + c_{\lambda_j} \|^2} \right]$.
%\eeq 
%where $b_{\lambda}$ and $c_\lambda$. % are as in Thm. \ref{thm-h-inf}. %defined as in \eqref{def-b-c-lambda}. 
We conclude from \eqref{eq-h-inf-tf} that
%\beq %\|G\|_\infty^2
$(H_\infty)^2 = \sup_{\omega\in\mathbb{R}} \underset{j=1,...,d}{\max}\left[\frac{\alpha^2  \lambda_j}{2
	\|e^{2i\omega} + b_{\lambda_j} e^{iw} + c_{\lambda_j} \|^2} \right]$.
%	\eeq
%where we dropped the subscripts of $G_{\hat{w},\hat{z}}$ as it is clear from the context. 
Note that the supremum over $\omega$ is attained for some $\omega = \omega_*$ as it is the supremum of a continuous, $2\pi$-periodic function. Therefore, we can interchange the maximum over the index $j$ and the maximum over $\omega$, i.e. 
\beq\label{eq-max-lambda-hinf} %\|G\|_\infty^2 
(H_\infty)^2=  \underset{\lambda \in \{\lambda_j\}_{j=1}^d} {\max} \max_{\omega} h_\omega(\lambda), \quad  h_\omega (\lambda) := \left[\frac{\alpha^2  \lambda}{2
	\|e^{2i\omega} + b_\lambda e^{iw} + c_\lambda \|^2} \right]. \label{def-h-omega} 
\eeq
Next, we will show that the function $h_\omega(\lambda)$ is a quasi-convex function of $\lambda$ which then would directly imply that the maximum of  $h_\omega(\lambda)$ over $\lambda \in [\mu,L]$ is attained at a boundary point, i.e. attained for either $\lambda = \mu$ or $\lambda = L$. For this purpose, we write
	%\beq 
	$h_\omega (\lambda) = \left[\frac{\alpha^2 }{D_\omega(\lambda)} \right]$,
	%\eeq
with	$ D_\omega(\lambda):= \frac{2}{\lambda} \left( (\cos(2\omega)+ b_\lambda \cos(\omega) + c_\lambda)^2 + (\sin(2\omega) + b_\lambda \sin(\omega))^2\right).$
%We will show that for every $\omega$ fixed, this function is quasi-convex with respect to its parameter $\lambda$ on the inverval $(0,\infty)$. 
Since $b_\lambda$ and $c_\lambda$ are linear in $\lambda$, %the denominator of $h_\omega(\lambda)$ is quadratic in $\lambda$, 
we can also write 
$ D_\omega(\lambda) = \tilde{a}_\omega \lambda  + \frac{ \tilde{b}_\omega}{\lambda} + \tilde{c}_\omega, $
for some constants $\tilde{a}_\omega$, $\tilde{b}_\omega$, $\tilde{c}_\omega$ that depends only on the parameters $\alpha$, $\beta$, $\nu$ and $\omega$. If $\tilde{b}_w \leq 0$ then $D_\omega(\lambda)$ is concave and therefore quasi-concave. %\footnote{A real-valued function $h(\lambda)$ defined on an interval is \emph{quasi-concave} if $-h(\lambda)$ is quasi-convex.%the level sets $U_a := \{x \in \mathbb{R} : f(x)\geq a\}$ are convex for every $a\in \mathbb{R}$.
%} 
Consequently, $h_\omega(\lambda)$ is quasi-convex. %\footnote{A real-valued function $h(\lambda)$ defined on an interval is \emph{quasi-convex} if the lower level sets $L_a := \{x \in \mathbb{R} : f(x)\leq a\}$ are convex for every $a\in \mathbb{R}$.} 
If $\tilde{b}_w > 0$, then $D_\omega(\lambda)$ is convex on $(0,\infty)$ and therefore $h_\omega(\lambda)$ is convex (and therefore quasi-convex) as being the composition of two convex functions $\lambda \mapsto \alpha^2/\lambda$ and $\lambda\mapsto D_\omega(\lambda)$. In either case, $h_\omega(\lambda)$ is quasi-convex on the interval $[\mu, L]$ for any fixed $\omega$. The maximum of quasi-convex functions is also quasi-convex; therefore $\max_{\omega} h_\omega(\lambda)$ is quasi-convex. By quasi-convexity, $h_\omega(\lambda)$ attains its maximum on the interval $[\mu,L]$ either at $\lambda=\mu$ or $\lambda=L$. We conclude from \eqref{eq-max-lambda-hinf} that 
{\small
\beq % \|G\|_\infty^2 
%L_{2,*}^2 =
 (H_\infty)^2&=& \underset{\lambda \in \{\mu, L\}} {\max} \max_{\omega} h_\omega(\lambda)
= \underset{\lambda \in \{\mu, L\}} {\max} \max_{\omega\in\mathbb{R}}
\left[\frac{\alpha^2  \lambda}{2
	\|e^{2i\omega} + b_\lambda e^{iw} + c_\lambda \|^2} \right].
	\label{eq-hinf-norm-sq}
%\\ 
%	&=& \underset{\lambda_i \in \{\mu, L\}}{\max} \frac{\alpha^2 \lambda_i }{ 2\big | |1 + c_{\lambda_i} | - |b_{\lambda_i}| \big |^2}
\eeq
} 
We next provide a lemma and its proof, which will help in characterizing the right-hand side. %The proof of this lemma will be postponed to Appendix \ref{app-proof-hinf-opt-lemma}.
\begin{lemm}\label{lemma-hinf-calc} In the setting of Theorem \ref{thm-h-inf}, let $\lambda>0$ be given. We have 
{\small
\beqs \min_{\omega\in\mathbb{R}}\left[{
	\|e^{2i\omega} + b_{\lambda} e^{iw} + c_{\lambda} \|^2} \right] 
	&=& \begin{cases} (1-c_\lambda)^2 (1 - \frac{b_\lambda^2}{4c_\lambda}) & \mbox{if   } c_\lambda>0 \mbox{ and } \frac{|b_\lambda| (1+c_\lambda)}{4c_\lambda} < 1,\nonumber \\
	(|1+c_\lambda| - |b_\lambda|)^2 & \mbox{otherwise},
	\end{cases}%\\
	%&=&r_\lambda	
\eeqs
}
and the minimum is achieved at
%$$ \omega_*(\lambda):= \begin{cases} 
%-\frac{b_\lambda (1+c_\lambda)}{4c_\lambda}  & \mbox{if } c_\lambda>0 \mbox{ and } \frac{|b_\lambda| (1+c_\lambda)}{4c_\lambda} < 1,\\
%0 			& 			\mbox{if }    
%\end{cases}$$
\begin{small}
$$ \omega_*(\lambda)= \begin{cases} 
\mbox{arccos}(-\frac{b_{\lambda} (1+c_{\lambda})}{4 c_{\lambda}})  & \mbox{if}\quad c_{\lambda}>0 \mbox{ and } {|b_{\lambda} (1+c_{\lambda})|} < 4|c_{\lambda}| ,\nonumber \\
\pi & \mbox{if} \quad c_{\lambda} \leq 0 \mbox{ and }  (1-b_{\lambda} + c_{\lambda_*})^2 \leq (1+b_{\lambda} + c_{\lambda})^2, \\
\pi &  \mbox{if} \quad |b_{\lambda}(1+c_{\lambda})| \geq 4|c_{\lambda}| \mbox{ and }  (1-b_{\lambda} + c_{\lambda})^2 \leq (1+b_{\lambda} + c_{\lambda})^2, \\
0 & otherwise.
\end{cases}
$$
\end{small}
\end{lemm}	
%\begin{proof}
\proof {\textbf{Proof of Lemma $\ref{lemma-hinf-calc}$.}} After a straightforward computation, 
\beqs \|e^{2i\omega} + b_{\lambda} e^{i\omega} + c_{\lambda} \|^2  %&=&(e^{2i\omega} + b_{\lambda} e^{i\omega} + c_{\lambda})^* (e^{2i\omega} + b_{\lambda} e^{i\omega} + c_{\lambda})\\
%&=&
%1 + b_\lambda^2 + c_\lambda^2 + 2 \left[ b_\lambda \cos(\omega) + c_\lambda \cos(2\omega) + b_\lambda c_\lambda \cos(\omega) \right]\\
&=& 1 + b_\lambda^2 + c_\lambda^2 + 2 \left[ b_\lambda \cos(\omega) + c_\lambda (2\cos^2(\omega)-1) + b_\lambda c_\lambda \cos(\omega) \right].
\eeqs	
Therefore, with the change of  variable $q = \cos(w)$, we have $q\in [-1,1]$ and 
\beq 
\min_{\omega\in\mathbb{R}}\left[{
	\|e^{2i\omega} + b_{\lambda} e^{iw} + c_{\lambda} \|^2} \right]  &=& \min_{q\in [-1,1]} E(q):=1 + b_\lambda^2 + c_\lambda^2 + 2 \left[ b_\lambda q + c_\lambda (2q^2-1) + b_\lambda c_\lambda q \right].
\label{eq-transformed}
\eeq
%where $E(q):= 1 + b_\lambda^2 + c_\lambda^2 + 2 \left[ b_\lambda q + c_\lambda (2q^2-1) + b_\lambda c_\lambda q \right]$. 
If $c_\lambda = 0$,  $E(q)$ is a linear function of $q$ with $E(q)=1+b_\lambda^2 + 2 b_\lambda q$. In this case, it is straightforward to verify that
$\min_{q\in[-1,1]} E(q) = (1-|b_\lambda|)^2.$
For $c_\lambda\neq 0$, the objective $E(q)$ is a quadratic with a gradient 
$ \nabla E(q) = 2 [b_\lambda + 4 c_\lambda q + b_\lambda c_\lambda],$
which vanishes at the point $q_*:= -\frac{b_\lambda (1+c_\lambda)}{4c_\lambda}$. For $c_\lambda>0$, $E(q)$ is a strongly convex quadratic in $q$, therefore its minimum on $[-1,1]$ is attained at $q_*$ if $q_* \in (-1,1)$. Otherwise, it will be attained at a boundary point $q=-1$ or $q=1$. If $c_\lambda < 0$, then $E(q)$ is strictly concave and its minimum on $[-1,1]$ will be attained at a boundary point when $q=-1$ or $q=1$. We also compute that
%\beq 
$E(-1) = 1 + b_\lambda^2 + c_\lambda^2 + 2 \left[ -b_\lambda  + c_\lambda -b_\lambda  c_\lambda  \right] = (1-b_\lambda + c_\lambda)^2$, 
$E(1) =  1 + b_\lambda^2 + c_\lambda^2 + 2 \left[ b_\lambda  + c_\lambda + b_\lambda  c_\lambda  \right] = (1+b_\lambda + c_\lambda)^2$, %.\label{eq-opt-angle-2}\\
and $E(q_*) = (1-c_\lambda)^2 (1 - \frac{b_\lambda^2}{4c_\lambda}) \quad \mbox{if} \quad c_\lambda \neq 0$. {Therefore we conclude that if $c_\lambda > 0$ and $|q_*|= \frac{|b_\lambda (1+c_\lambda)|}{4|c_\lambda|} < 1$, the minimum of $E(q$) in \eqref{eq-transformed} is achieved at $q_*$ with a value $E(q_*)$. Otherwise it is achieved at $q=-1$ when $E(-1) \geq  E(1)$ or at $q=1$ when $E(-1)\leq E(1)$ and the minimum value is $\min(E(1), E(-1)) = (|1+c_\lambda| - |b_\lambda|)^2$. Noting the relationship $q=\arccos(w)$, we conclude. %\myqed
}
\endproof
%\vspace{-0.1in}
Equipped with this lemma, the stage is set to finalize the proof of Thm. \ref{thm-h-inf}. We start with noting that applying Lemma \ref{lemma-hinf-calc} to \eqref{eq-hinf-norm-sq}, the desired equality \eqref{hinfty-quad} follows. Next, we prove the lower bound on the $H_\infty$ norm. It follows from the proof of Lem. \ref{lemma-hinf-calc} above that if $c_\lambda>0$ and $|q_*| = \frac{|b_\lambda| (1+c_\lambda)}{4c_\lambda} < 1$, then 
$ E(q_*) = (1-c_\lambda)^2 (1 - \frac{b_\lambda^2}{4c_\lambda}) < \min (E(1), E(-1)) = (|1+c_\lambda|-|b_\lambda|)^2$. By taking the square root of both sides, we obtain
$ r_\lambda \leq  | | 1 + c_\lambda| - |b_\lambda| |,$
with strict inequality in case $c_\lambda>0$ and $|b_\lambda (1+c_\lambda)| < 4|c_\lambda |$ where $r_\lambda$ is defined by \eqref{def-r-lambda}. On the other hand, by the triangle inequality for $s_1, s_2 \in \mathbb{R}$,  we have 
$ | | s_1 + s_2| - | s_1| | \leq |s_2|$. Choosing $s_2 = \alpha \lambda $, $s_1 = -b_\lambda$; we obtain $s_1 + s_2 = 1+ c_\lambda$ and
$ | | 1 + c_\lambda| - |b_\lambda| | \leq \alpha \lambda. $
We conclude that $r_\lambda \leq \alpha \lambda$ with the equality holding only when $c_\lambda \leq 0$ or $|b_\lambda (1+c_\lambda)| \geq 4|c_\lambda |$. Then, it follows from \eqref{hinfty-quad} that
%\begin{equation}
	$H_\infty \geq  \frac{\alpha}{\sqrt{2}} \max_{\lambda \in \{ \mu, L\} }  \frac{ \sqrt{\lambda}}{\alpha \lambda} = \frac{1}{\sqrt{2\mu}}$. Also, from  \eqref{hinfty-quad}, the latter equality can hold only if $r_\mu = \alpha \mu$ and 
	 $\frac{\alpha}{\sqrt{2}}  \frac{ \sqrt{L}}{r_L} \leq \frac{1}{\sqrt{2\mu}}$. This is equivalent to the fact that parameters lie in the set $\mathcal{S}_1 \cap \mathcal{S}_2$, completing the proof of Theorem \ref{thm-h-inf}. %\looseness=-1 
	 %\myqed
%\end{proof} 
%\section{\label{app-proof-hinf-opt-lemma}Proof of Lemma \ref{lemma-hinf-calc}.}
%\begin{proof} 

%\end{proof}

%%\vspace{-0.3in}
\section{\label{app-proof-of-worst-case-noise}Proof of Proposition \ref{prop-worst-case-noise}.}
%\begin{proof}
For the proof, we resort to a frequency domain approach based on Fourier analysis techniques. Given input noise sequence $\{w_k^{(h)}\}_{k\geq 0}$ and output sequence $\{z_k^{(h)}\}_{k\geq 0}$, we consider the discrete-time Fourier transform (DTFT) of the input noise sequence and output sequence defined as
%\beq 
$W_h (e^{i\omega }) = \sum_{k=-\infty}^{\infty} w_k^{(h)} e^{-i \omega k},  \quad Z_h(e^{i\omega}) = \sum_{k=-\infty}^{\infty} z_k^{(h)} e^{-i \omega k}$,
%\label{def-dtft}
%\eeq 
 %\cite{hinrichsen1991stability,zhou1996robust}.  
(see e.g. \cite{oppenheim1999discrete}) where we use the convention that $w_k^{(h)} = z_k^{(h)} = 0$ for $k<0$. Since Fourier transforms $W_h(e^{i\omega })$ and $Z_h(e^{i\omega})$ are periodic with period $2\pi$, it suffices to consider $\omega \in [0,2\pi)$ to characterize them. It is also well known that we have
	$Z_{h}(e^{i\omega}) = G(e^{i\omega}) W_h (e^{i\omega})$
where $G(e^{i\omega})$ is the transfer function matrix defined in \eqref{def-transfer-matrix}. Furthermore, by Parseval's identity, 
 \beq \sum_{k\geq 0 }\|w_k^{(h)}\|^2  = \frac{1}{2\pi}  \int_{0}^{2\pi} \| W_h (e^{i\omega}) \|^2 d\omega, \quad \sum_{k\geq 0 }\|z_k^{(h)}\|^2  = \frac{1}{2\pi}  \int_{0}^{2\pi} \| Z_h (e^{i\omega}) \|^2 d\omega,
 \label{eq-parseval}
 \eeq
(see e.g. \cite{oppenheim1999discrete,stein2003shakarchi}). Also, from \eqref{eq-transfer-fun-factor}, we see that $G(e^{i\omega})$ has the form 
\beq \hspace{-0.25in}
G(e^{i\omega}) &=&  \underset{j=1,...,d}{\Diag}\left[ R_\omega(\lambda_j)\right] U^T, 
\eeq
with
\beq
R_{\omega}(\lambda)&:=& \frac{1}{\sqrt{2}} \frac{-\alpha \sqrt{\lambda}}{e^{i\omega}-(1+\beta-\beta e^{-i\omega})+\alpha\lambda (1+\nu-\nu e^{-i\omega})}.%$$
%Note that $R_\omega(\lambda).%\\
%&=&
\label{eq-input-to-output-fourier}
\eeq
%where 
%$$ R_{\omega}(\lambda):= \frac{1}{\sqrt{2}} \frac{-\alpha \sqrt{\lambda}}{e^{i\omega}-(1+\beta-\beta e^{-i\omega})+\alpha\lambda_i (1+\nu-\nu e^{-i\omega})}.$$
Note that $R_\omega(\lambda)$ can be complex-valued depending on the choice of $\omega$, % we could write $R_\omega(\lambda) = |R_\omega(\lambda)|e^{i \theta(R_\omega(\lambda))}$ where $\theta(R_\omega(\lambda)) \in [0,2\pi)$ is the argument of the complex number $R_\omega(\lambda)$. Note 
and we have also $|R_\omega(\lambda)| = \sqrt{h_{\omega}(\lambda)}$ where $h_\omega(\lambda)$ is as in \eqref{def-h-omega}. Therefore, we can write 
%\beq 
$$G(e^{i\omega}) =  \underset{j=1,...,d}{\Diag} \left[  \frac{R_\omega(\lambda_j)}{|R_\omega(\lambda_j)| }  \right]
\underset{j=1,...,d}{\Diag}\left[ \sqrt{h_{\omega}(\lambda_j)} \right] U^T.$$ 
%\label{def-ustar-vstar}
%\eeq
In light of Lemma \ref{lemma-hinf-calc}, note that $\lambda_*$ and $\omega_*$ are maximizers of the optimization problem
%Let $\lambda_* \in \{\mu, L\}$ and $\omega_* \in [0,\pi]$ be a point where the maximum is attained 
in \eqref{eq-hinf-norm-sq}, and by symmetry $-\omega_*$ is also a maximizer. Also, we have $H_\infty = \|G(e^{i\omega_*)}\|_2 = \sqrt{h_{\omega_*}(\lambda_*)}$ by the definition of the $H_\infty$ norm. Let 
\begin{equation}
u_* = \begin{cases} u_1 & \mbox{if } \lambda_* = \mu \\
                      u_d & \mbox{if } \lambda_* = L     
\end{cases}, \quad
v_* = \begin{cases} \frac{R_{\omega_*}(\mu)}{\| R_{\omega_*}(\mu)\|}  e_1 & \mbox{if } \lambda_* = \mu, \\
                   \frac{R_{\omega_*}(L)}{\| R_{\omega_*}(L)\|}  e_d & \mbox{if } \lambda_* = L,     
\end{cases}
\label{def-ustar-vstar}
\end{equation}
where $e_i$ is the $i$-th standard basis vector. Then, it is straightforward to check that 
%    \beq
     $\varepsilon_* G(e^{i\omega_*}) u_* = v_*$ and $\varepsilon_* v_*^H   G(e^{i\omega_*}) = u_*^H$.
   % 		\label{eq-sing-vecs}
    %\eeq
%where $^H$ denotes Hermitian transpose. 
Assume for now that the noise sequence is given by
\beq w_{k}^{(h)} =\sqrt{h (2-h)} (1-h)^k e^{i\omega_* k} u_* \quad \mbox{for} \quad k\geq 0,
\label{eq-complex-worst-case-seq}
\eeq 
which can potentially have complex entries. Using the definition of DTFT, %\eqref{def-dtft}, 
note that we can compute
 $ W_h (e^{i\omega}) = \frac{ \sqrt{h (2-h)} }{1- (1-h) e^{i(\omega_*-\omega)} } u_*.$   
Furthermore, it is easy to check that $\sum_ k \|w_k^{(h)}\|^2 = h (2-h)\|u_*\|^2 \sum_{k\geq 0}(1-h)^{2k} = 1$ for every $h$; therefore by \eqref{eq-parseval}, we have for every $h\in (0,1)$,
\beq 
 \int_{0}^{2\pi} \| W_h (e^{i\omega}) \|^2 d\omega= \int_{0}^{2\pi}\frac{{h (2-h)} }{|1- (1-h) e^{i(\omega_*-\omega)}|^2 }d\omega   =2\pi. 
 \label{def-integral-noise-freq-domain}
 \eeq
Similarly, by \eqref{eq-parseval} and using $Z_{h}(e^{i\omega}) = G(e^{i\omega}) W_h (e^{i\omega})$, %and \eqref{eq-input-to-output-fourier},
\begin{small}
   \beq  \sum_{k\geq 0}\| z_k^{(h)} \|^2  &=&  \frac{1}{2\pi}   \int_{0}^{2\pi} \| Z_h (e^{i\omega}) \|^2 d\omega %=  \frac{1}{2\pi}   \int_{0}^{2\pi} \| Z_h (e^{i\omega}) \|^2 d\omega\\
  =   \int_{0}^{2\pi}   \|G(e^{i\omega})  u_*\|^2  K_h(\omega) d\omega,
 \label{eq-l2norm-of-output-delta}
   \eeq
\end{small}where $K_h(\omega) := \frac{1}{2\pi}\frac{ {h(2-h)} }{|1- (1-h) e^{i(\omega_*-\omega)}|^2 } = 
 \frac{1}{2\pi} \frac{h(2-h)}{2+h^2 - 2 h - 2(1-h)\cos(\omega - \omega_*)}$. The family $K_h$ is a family of \emph{good kernels} on the interval $I = [0,2\pi]$ in the sense of \cite{stein2003shakarchi}, i.e. it satisfies the following three properties: (i) $K_h(\omega)\geq 0$ for any $\omega \in \mathbb{R}$ and $h \in (0,1), $
 (ii) $\int_I K_h(\omega)d\omega= 1$ based on \eqref{def-integral-noise-freq-domain},
 (iii) For every $a>0$, $\int_{I, \|w-w_*\|>a} K_h (\omega) d\omega \to 0$ as $h \to 0$. This is because for $|\omega-\omega_*| > a$, we have 
 $K_h(\omega) \leq  \frac{1}{2\pi} \frac{h(2-h)}{2+h^2 - 2 h - 2(1-h)|\cos (a)|} $
 so that 
 $\int_{I, \|w-w_*\|>a}K_h(\omega)d\omega \leq \frac{h(2-h)}{2+h^2 - 2 h - 2(1-h)|\cos (a)|}$
and the right-hand side goes to zero as $h\to 0$. 
% \end{itemize}
Then, based on standard arguments in Fourier analysis similar to \cite[Sec. 2]{stein2003shakarchi}, it follows that\footnote{The limit of $K_h$ as $h\to 0$ can be viewed as the ``Dirac delta" function at $\omega= \omega_*$ in the sense of generalized functions studied in \cite{lighthill1958introduction}.} 
$  \lim_{h \to 0 }\int_{0}^{2\pi}   \|G(e^{i\omega})  u_*\|^2  K_h(\omega) d\omega = \|G(e^{i\omega_*})  u_*\|^2  = \left(\frac{1}{\varepsilon_*}\right)^2 \|v_*\|^2 =  \left(\frac{1}{\varepsilon_*}\right)^2 = H_\infty^2.$
 We conclude from \eqref{eq-l2norm-of-output-delta} that 
    $ \lim_{h\to 0}\sum_{k\geq 0 }\|z_k^{(h)}\|^2 = H_\infty^2.$
 Noting $\sum_{k\geq 0 }\|w_k^{(h)}\|^2=1$, this shows that the noise sequence \eqref{eq-complex-worst-case-seq} attains the worst-case behavior as $h\to 0$. However, this sequence can be complex-valued when $\omega_* \in (0,\pi)$. Next, we argue that from the symmetry of the problem, we can simply take the real part of this sequence as the worst-case noise. Recall that by passing to the complex conjugates, $-\omega_*$ is also a maximizer of $G(e^{i\omega})$, i.e. $\| G(e^{i\omega_*}\| = \|G(e^{-i\omega_*})\|$. Therefore, by replacing $\omega_*$ with $-\omega_*$ in the above analysis, we can also show that the complex conjugate $(w_k^{(h)})^* = \sqrt{h (2-h)} (1-h)^k e^{-i\omega_* k} u_*$ of \eqref{eq-complex-worst-case-seq} is also an almost worst-case sequence with an output $(z_k^{(h)})^*$. % which is the complex conjugate of the output $z_k^{(h)}$).
  Finally, by linearity of the system, the average 
 $\frac{1}{2}[(w_k^{(h)})^* + (w_k^{(h)})] =  \sqrt{h (2-h)} (1-h)^k \cos(\omega_* k) u_* $ as a noise input also leads to the worst-case ($\ell_2$ gain) behavior. % with $ \lim_{h\to 0}\sum_{k\geq 0 }\|z_k^{(h)}\|^2 = H_\infty^2$. 
 We conclude. %\myqed %\looseness=-1
\section{\label{app-proof-of-real-stab-rad}Proof of Theorem \ref{thm-real-hinf}.}

We first consider part $(i)$. The fact that $r_{\mathbb{C}}(A_Q, B, T) = H_\infty^{-1} = L_{2,*}^{-1}$ is a consequence of \eqref{def-complex-stab-rad}. To show $r_{\mathbb{C}}(A_Q, B, T)=r_{\mathbb{R}}(A_Q, B, T)$, by the definitions of these quantities, it suffices to show that there exists a real matrix $\Delta_* \in \mathbb{R}^{d\times d}$ of norm $\varepsilon_*$
with $\rho(A+B\Delta_*C) = 1$. It is known that the matrix
$\Delta_* = \varepsilon_* u_* v_*^H$ satisfies $\rho(A+B\Delta_*C) = 1$ where $u_*$ and $v_*$  given in \eqref{def-ustar-vstar} are the right and left eigenvectors of the transfer matrix $G(e^{i\omega_*})$, see e.g. \cite{GGO}.  %Furthermore, $\omega_* \in {0,\pi\}$ $R_{\omega_*}(\lambda_*)$ is real-value
%o show $L_{2,*}=H_\infty$, it suffices to show that the almost worst-case sequence $w_k^{(h)}$ has real entries, i.e. $w_k^{(h)}\in\mathbb{R}^d$
Note that it follows from the proof argument of Prop. \ref{prop-worst-case-noise} that $\|G(e^{iw})\|$ is maximized on the unit circle for $\omega=\omega_*$ given in Prop. \ref{prop-worst-case-noise}. We first consider part $(i)$ where the condition $c_{\lambda_*} \leq 0   \quad \mbox{or} \quad {|b_{\lambda_*}(1 + 	c_{\lambda_*})|}\geq 4| c_{\lambda_*}|$ means that $\omega_* = 0$ or $\omega_*=\pi$, 
%In this case, it follows from the proof argument of Prop. \ref{prop-worst-case-noise} that $\|G(e^{iw_*})\|$ is maximized on the unit circle for either 
%$w_* = 0$ or $w_* = \pi$
in which case the matrix $G(e^{iw_*})$ will be a matrix with all real entries and $R_{\omega_*}(\lambda_*)$ given in \eqref{def-ustar-vstar} will be real. Therefore $v_*$ is real and since the Hessian of $f$ is symmetric, we can choose the eigenvector $u_*$ in \eqref{eq-worst-case-noise} to have real entries. Then, we conclude that $\Delta_* = u_* v_*^H\in \mathbb{R}^{d\times d}$ and this completes the proof of part $(i)$. %and in this case $w_k^{(h)} \in \mathbb{R}^d$ for every $h$ and $k$ and we conclude.% This together with
We next consider part $(ii)$. The proof of the lower bound on $r_{\mathbb{R}}(A_Q, B, T)$ is straightforward, by definition $r_{\mathbb{R}}(A_Q, B, T) \geq r_{\mathbb{C}}(A_Q, B, T)$ and the latter quantity is equal to $L_{2,*}^{-1} =H_\infty^{-1}$ by \eqref{def-complex-stab-rad}. %, implying the desired lower bound for $r_{\mathbb{R}}(A_Q, B, T)$. 
Next, we prove the upper bound. %ower bound on $r_{\mathbb{R}}(A_Q, B, T)$. 
From  \cite{qiu1995formula}, we have the representation
{\small
 \beq r_{\mathbb{R}}(A_Q, B, T) = \bigg(\sup_{\omega \in [0,2\pi]} \tau_1( G(e^{i\omega}) \bigg)^{-1},  \label{eq-formula-real-stab-rad} \eeq 
 with 
 \beqs  \tau_1(M):= \inf_{\gamma \in (0,1]} \sigma_2\left( \begin{bmatrix} \mbox{Re}(M) & -\gamma \mbox{Im}(M) \\ -\gamma^{-1} \mbox{Im}(M)  & \mbox{Re}(M) \end{bmatrix}\right),%\label{def-real-stab-rad}
 \eeqs}where $\sigma_2$ denotes the second largest singular value. Hence, evaluating $\tau_1(G(e^{i\omega}))$ requires solving a minimization problem.  %\cite{qiu1995formula} and 
A closed-form formula for it is not available %(. However, for the function $\tau_1$ there appears to be no (explicit formula in terms of the entires of the $G$ matrix) closed-form solution 
except in special cases such as when the imaginary part of $G$ is of rank one \cite{qiu1995formula}. 
That being said, by \cite[Prop. 6.7.1]{karow2003geometry}, we have the lower bound
	$\tau_1(G(e^{i\omega})) \geq \sigma_2(G(e^{i\omega})),$ 
for any $\omega$. From the analysis given in \eqref{eq-transfer-fun-factor}--\eqref{eq-max-lambda-hinf} and in the rest of the proof of Thm. \ref{thm-h-inf}, it follows that the singular values of $G(e^{i\omega})$ are equal to
$ \sqrt{h_\omega(\lambda)} = \frac{\alpha \sqrt{\lambda}}{\sqrt{2}\|e^{2i\omega} + b_{\lambda} e^{i\omega} + c_{\lambda} \|}$ for $\lambda \in \{\mu=\lambda_1, \lambda_2, \dots, \lambda_d=L\},$
where the function $h_\omega(\lambda)$ is quasi-convex on the interval $[\mu,L]$ as a function of $\lambda$ for $\omega$ fixed. For fixed $\omega$, with some abuse of notation, let $\lambda_*(\omega)$ be the largest singular value of $G(e^{i\omega})$ which coincides with the largest value of $\sqrt{h_\omega(\lambda)}$ on $[\mu, L]$. Therefore, by quasi-convexity $\lambda_*(\omega)=\lambda_1 = \mu$
or $\lambda_*(\omega)=\lambda_d = L$. In the former case, the maximum of $h_\omega(\lambda)$ on the interval $[\lambda_2, L]$ occurs at a boundary point, either at $\lambda=\lambda_2$ or $\lambda = L$ so that %. %In the latter case, 
$\sigma_2(G(e^{i\omega})) \geq \max_{\lambda \in \{\lambda_2, L\}} 
\frac{\alpha \sqrt{\lambda}}{\sqrt{2}\|e^{2i\omega} + b_{\lambda} e^{i\omega} + c_{\lambda} \|}. 
$
In the latter case when $\lambda_*(\omega)=\lambda_d=L$, similarly by the quasi-convexity of $h_\omega(\lambda)$,
$\sigma_2(G(e^{i\omega}) \geq \max_{\lambda \in \{\mu, \lambda_{d-1}\}} 
\frac{\alpha \sqrt{\lambda}}{\sqrt{2}\|e^{2i\omega} + b_{\lambda} e^{i\omega} + c_{\lambda} \|}. 
$ 
Taking pointwise minimum of both lower bounds, 
%\begin{small}
%$$ \sigma_2(G(e^{i\omega}) \geq
%\min \bigg( 
%\max_{\lambda \in \{\lambda_2, L\}} 
%\frac{\alpha \sqrt{\lambda}}{\sqrt{2}\|e^{2i\omega} + b_{\lambda} e^{i\omega} + c_{\lambda} \|} ,  \max_{\lambda \in \{\mu, \lambda_{d-1}\}} 
%\frac{\alpha \sqrt{\lambda}}{\sqrt{2}\|e^{2i\omega} + b_{\lambda} e^{i\omega} + c_{\lambda} \|}  
%\bigg),
%$$
%\end{small}
and taking supremum over $\omega$ we obtain the lower bound
$ \sup_{\omega \in [0,2\pi]} \tau_1( G(e^{i\omega})   \geq
H_\infty^{\tiny\mbox{lb}}.$
Using \eqref{eq-formula-real-stab-rad}, this implies $r_{\mathbb{R}}(A_Q, B, T) \leq \big( H_\infty^{\tiny\mbox{lb}} \big)^{-1}$ completing the proof of part $(ii)$. %\myqed
\section{\label{app-Lyap-eval-GMM}Proof of Lemma \ref{lemma-gmm-lyap-evolution}.}
%\begin{proof} %Let $\tilde{P}\in\mathbb{R}^{2\times 2}$ be a given symmetric matrix that satisfies the constraint \eqref{constraint-sdp}.  
We first provide another lemma and its proof.% This lemma will be useful for establishing the proof of Lemma \ref{lemma-gmm-lyap-evolution}.
\begin{lemm}\label{lemma tmm-noisy-suboptimality} In the setting of Lemma \ref{lemma-gmm-lyap-evolution}, we have%Consider the GMM algorithm with constant parameters$ (\alpha,\beta,nu)$ and gradient noise $w_k$ at step $k$. %Then 
%\begin{eqnarray}
%x_{k+1}&=&x_{k}-\alpha (\nabla f( y_k ) + w_k) +\beta(x_{k}- x_{k-1}), \\
%y_{k}&=&x_k+\nu (x_k- x_{k-1}), 
%\end{eqnarray}
%where $w_k$ is the gradient noise at step $k$. 
%we have
 \beq
  f(x_{k+1}) - f_* &\leq& \rho_0^2 (f(x_k) - f_*) - 
 \begin{bmatrix} 
 	{\xi}_k^c \\
 	\nabla f(y_k) 
 \end{bmatrix}^T  
 (\tilde X_0 \otimes I_d)
 \begin{bmatrix} 
 	{\xi}_k^c\\
 	\nabla f(y_k) 
 \end{bmatrix}^T + \frac{L}{2}\alpha^2 \|w_k\|^2  \nonumber \\
 &&- \alpha (1-L\alpha) \nabla f(y_k)^T w_k + L\alpha (\beta - \nu) w_k^T(x_{k-1}-x_k).
 \label{ineq-to-prove}
 \eeq
%where $X_0 =  \tilde{X}_0 \otimes I_d$ with
%$\tilde{X} = \tilde{X}_1 + \rho^2 \tilde{X}_2 + (1-\rho^2) \tilde{X}_3$ and
%\beq \tilde{X} = \tilde{X}_1 + \rho^2 \tilde{X}_2 + (1-\rho^2) \tilde{X}_3\eeq
%with
%\beq \begin{scriptsize} \quad\quad
% ~\tilde{X}_0 = \tilde{X}_1 + \rho_0^2 \tilde{X}_2 + (1-\rho_0^2) \tilde{X}_3, ~ \tilde{X}_1 = \frac{1}{2}
%\begin{bmatrix}
%-Lh^2 & Lh^2 & -(1-\alpha L)h \\
%Lh^2 & -Lh^2 & (1-L\alpha)h \\
%-(1-L\alpha)h & (1-L\alpha)h & \alpha(2-L\alpha)
%\end{bmatrix},
%\end{scriptsize}
%%\label{def-tilde-X}
%\eeq 
%\begin{align*}\begin{scriptsize} %\begin{small}
%\tilde{X}_2 = \frac{1}{2}
%\begin{bmatrix} 
%	\nu^2\mu& -\nu^2\mu & -\nu \\
%    -\nu^2\mu & \nu^2\mu & \nu  \\
%    -\nu  & \nu  & 0
%\end{bmatrix},~ 
%\tilde{X}_3 = \frac{1}{2}
%\begin{bmatrix} 
%	(1+\nu)^2\mu & -\nu(1+\nu)\mu & -(1+\nu) \\
%    -\nu(1+\nu)\mu & \nu^2\mu & \nu  \\
%    -(1+\nu)  & \nu  & 0 
%\end{bmatrix},
%\end{scriptsize} 
%\end{align*}
%and $h:= \beta - \nu$.
\end{lemm} 
\proof {\textbf{Proof of Lemma \ref{lemma tmm-noisy-suboptimality}.}} The proof of \cite[Lemma 5]{hu2017dissipativity} concerns the deterministic case for GMM when $w_k = 0$ for every $k$. In the following, we follow the same proof technique and extend it to the inexact gradient case that can allow arbitrary values of $w_k$. We introduce $X_i = \tilde{X}_i\otimes I_d \in \mathbb{R}^{2d\times 2d}$ for $i=0,1,2,3$. By the second inequality of \eqref{ineq-strcvx-smooth}, we have
 \beqs 
 f(y_{k}) - f(x_{k+1}) &\geq& \nabla f(y_k)^T (y_k - x_{k+1}) - \frac{L}{2}\|x_{k+1}-y_k\|^2 %\nonumber \\
 %&=&\nabla f(y_k)^T \left( (\beta - \nu)(x_{k-1}-x_k) + \alpha \nabla f (y_k) + \alpha w_k \right) \\
 %&& - \frac{L}{2}  \| (\beta - \nu)(x_{k-1}-x_k) + \alpha \nabla f (y_k) + \alpha w_k \|^2 \nonumber \\
 = \begin{bmatrix} 
 	\txi_k \\
 	\nabla f(y_k) 
 \end{bmatrix}^T X_1
  \begin{bmatrix} 
 	\txi_k \\
 	\nabla f(y_k) 
 \end{bmatrix} - a_k,
 \nonumber
 \eeqs
where $ a_k :=   \frac{L}{2}\alpha^2 \|w_k\|^2 - \alpha (1-L\alpha) \nabla f(y_k)^T w_k + L\alpha (\beta - \nu) w_k^T(x_{k-1}-x_k).$ Also, 
 \beqs 
 f(x_k) - f(y_k) \geq \nabla f(y_k)^T (x_* - y_k) + \frac{m}{2}\|x_* - y_k \|^2  
  %&=& -\nabla f(y_k)^T \left( (1+\nu)(x_k - x_*) - \nu (x_{k-1}-x_*)\right) + \frac{m}{2} \|(1+\nu)(x_k - x_*) - \nu (x_{k-1}-x_*)\|^2 \\
  =
  \begin{bmatrix} 
 	\txi_k \\
 	\nabla f(y_k) 
 \end{bmatrix}^T X_2
  \begin{bmatrix} 
 	\txi_k \\
 	\nabla f(y_k) 
 \end{bmatrix},
 \eeqs
%with $X_2 = \tilde{X}_2 \otimes I_d$.
where we used the first inequality of \eqref{ineq-strcvx-smooth}. Similarly,  
\beqs 
f(x_*) - f(y_k) &\geq& \nabla f(y_k)^T(x_* - y_k) + \frac{m}{2}\|x_* - y_k\|^2 
%&=& -\nabla f(y_k)^T \left((1+\nu)(x_k - x_*) - \eta (x_{k-1}-x_*) \right) + \frac{m}{2} \|(1+\nu)(x_k - x_*) - \eta (x_{k-1}-x_*) \|^2 \\
= \begin{bmatrix} 
 	\txi_k \\
 	\nabla f(y_k) 
 \end{bmatrix}^T X_3
  \begin{bmatrix} 
 	\txi_k \\
 	\nabla f(y_k) 
 \end{bmatrix}. 
\eeqs
Therefore, by summing up these inequalities we obtain 
\beqs \scriptsize  f(x_k) - f(x_{k+1}) &\geq& \begin{bmatrix} 
 	\txi_k \\
 	\nabla f(y_k) 
 \end{bmatrix}^T (X_1 + X_2)
  \begin{bmatrix} 
 	\txi_k \\
 	\nabla f(y_k) 
 \end{bmatrix} - a_k, \\
f(x_*) - f(x_{k+1})  &\geq& \begin{bmatrix} 
 	\txi_k \\
 	\nabla f(y_k) 
 \end{bmatrix}^T (X_1 + X_3)
  \begin{bmatrix} 
 	\txi_k \\
 	\nabla f(y_k) 
 \end{bmatrix}  - a_k.
\eeqs
Consequently, we have
%\begin{eqnarray*}
%\rho^2 (f(x_k) - f(x_{k+1}) + (1-\rho^2) (f(x_*) - f(x_{k+1}))
%&=&
$\rho^2 (f(x_k) - f(x_*))  + f(x_*) - f(x_{k+1}) %\\
\geq \begin{bmatrix} 
 	\txi_k \\
 	\nabla f(y_k) 
 \end{bmatrix}^T X_0
  \begin{bmatrix} 
 	\txi_k \\
 	\nabla f(y_k) 
 \end{bmatrix} - a_k,$
% \end{eqnarray*}
and this yields the inequality \eqref{ineq-to-prove}. %\myqed
%\end{proof} 
Equipped with Lemma \ref{lemma tmm-noisy-suboptimality}, we are now ready to complete the proof of Lemma \ref{lemma-gmm-lyap-evolution}. Using strong convexity $f(x_k) - f(x_*) \geq \frac{\mu}{2}\|x_k - x_*\|^2$ and Lemma \ref{lemma tmm-noisy-suboptimality},
\begin{small}
\beqs
  f({x}_{k+1}) - f_* &\leq& (\rho_0^2 + \mg{ \rho_2^2}) (f({x}_k) - f_*) - 
 \begin{bmatrix} 
 	{\xi}_k^c \\
 	\nabla f({y}_k) 
 \end{bmatrix}^T  
 \big( \mg{(\tilde X_0 + \tilde{Z}) \otimes I_d} \big)
 \begin{bmatrix} 
 	{\xi}_k^c \\
 	\nabla f({y}_k) 
 \end{bmatrix}^T + \frac{L}{2}\alpha^2 \|{w}_k\|^2  \nonumber \\
 &&- \alpha (1-L\alpha) \nabla f(y_k)^T {w}_k + L\alpha (\beta - \nu) {w}_k^T({x}_{k-1}-{x}_k)  + \mg{\rho_3^2 (f(x_{k-1}) - f(x_*))} \\
 && 
 + \mg{\rho_1^2  \begin{bmatrix} 
 	{\xi}_k^c 
 \end{bmatrix}^T P \begin{bmatrix} 
 	{\xi}_k^c 
 \end{bmatrix}}.  
 \eeqs
 \end{small}
%\mg{ where $$\tilde{Z} = \begin{bmatrix} \rho_1^2 P_{11} +\frac{\mu}{2} \rho_2^2 &  \rho_1^2 P_{12} & 0\\
% 							   				 \rho_1^2 P_{12} & \rho_1^2 P_{22} +\frac{\mu}{2}\rho_3^2 & 0 \\
% 							   				 					0   &				0			& 0
% \end{bmatrix}$$}
%Recalling P = \tilde{P}\otimes I_d, 
By straightforward computations, we also have 
\begin{eqnarray} 
%\begin{footnotesize}
({\xi}_{k+1}^c)^T P {\xi}_{k+1}^c &=& (\rho_0^2 %\mg{\rho_1^2 + 
% \rho_2^2}
 ) ({\xi}_{k}^c)^T P {\xi}_{k}^c \nonumber \\
&&
 +  \begin{bmatrix} 
 	{\xi}_k^c \\
 	\nabla f({y}_k) \\
 	{w}_k
 \end{bmatrix}^T 
 \begin{bmatrix}	A^T P A - \rho_0^2%\mg{+\rho_1^2 + \rho_2^2}
  P  & A^T P B &  A^T P B\\
 						   B^T P A   &  B^T P B &   B^T P B \\
 						   B^T P A & B^T P B   & B^T P B
 \end{bmatrix} 
 \begin{bmatrix} 
 	{\xi}_k^c \\
 	\nabla f({y}_k) \\
 	{w}_k
 	\end{bmatrix}.
 	\label{eq-evolution-quad-Lyap}
%\end{footnotesize} 	
\end{eqnarray}
%Similarly,  we have
%\beqs \| {\xi}_{k+1}^c\|^2 = \rho^2 \| {\xi}_{k}^c\|^2 +  \begin{bmatrix} 
% 	{\xi}_k^c \\
% 	\nabla f({y}_k) \\
% 	{w}_k
% \end{bmatrix}^T 
% \begin{bmatrix}	A^T A - \rho^2 I_{2d}  & A^T  B &  A^T  B\\
% 						   B^T  A   &  B^T  B &   B^T  B \\
% 						   B^T  A & B^T  B   & B^T  B
% \end{bmatrix} 
% \begin{bmatrix} 
% 	{\xi}_k^c \\
% 	\nabla f({y}_k) \\
% 	{w}_k
% 	\end{bmatrix}
%\eeqs  
By summing up these inequalities and using $P=\tilde{P}\otimes I_d$, %or $\tilde{P}\succeq 0$ and $X_0 = \tilde{X}_0\otimes I_d$ where $\tilde{X}_0$ is as in Lemma \ref{lemma tmm-noisy-suboptimality}, 
\beq V_{P,c_1}({\xi}_{k+1})   &\leq&( \rho_0^2 \mg{+c_1\rho_1^2 + \rho_2^2}) V_{P,c_1}({\xi}_{k}) \| - \begin{bmatrix} 
 	{\xi}_k^c \\
 	\nabla f({y}_k) \nonumber \\
 	{w}_k
 \end{bmatrix}^T
 \bigg((\tilde{M_2}+ c_1\tilde{M}_1 )\otimes I_d\bigg)
 \begin{bmatrix} 
 	{\xi}_k^c \\
 	\nabla f({y}_k) \\
 	{w}_k
 \end{bmatrix} \\
 && \quad + \frac{c_1L}{2}\alpha^2 \|{w}_k\|^2 +  w_k^T B^T P  B w_k \mg{+ a\|w_k\|^2 + b c_1 \|\nabla f(y_k)\|^2} + \mg{\rho_3^2 V_{P,c_1}({\xi}_{k-1}) }. 
 \label{ineq-Vpc}
 \eeq
If we %introduce $M_0 = \tilde{M}_0 \otimes I_d$ and 
use the fact that ${y}_k -x_* = C {\xi}_k^c$, then 
%\beq 
$\scriptsize \quad \begin{bmatrix} 
 	{\xi}_k^c \\
 	\nabla f({y}_k) \\
 	{w}_k
 \end{bmatrix}^T 
 (c_0\tilde{M}_0 \otimes I_d)
 \begin{bmatrix} 
 	{\xi}_k^c \\
 	\nabla f({y}_k) \nonumber \\
 	{w}_k
 \end{bmatrix} = 2c_0(\mu+L) \bigg[\frac{\mu L}{\mu+L} \|{y}_k -x_* \|^2 - (y_k - x_*)^T \nabla f(y_k) %\\
 %&&\qquad \qquad \quad  
 + \frac{1}{\mu+L} \|\nabla f(y_k)\|^2\bigg] \leq 0$,
% \label{ineq-for-gd-rate}
%\eeq  
where last inequality is due to strong convexity and smoothness  \cite[Theorem 2.1.12]{nesterov2003introductory}.
%\mg{
%We have also 
%\beq - \|\nabla f(y_k)\|^2 &\leq& -2\mu \big( f(y_k) - f(x_*) \big)  
%= -2\mu \big(f(x_k) - f(x_*)\big) - 2\mu \big(f(y_k) - f(x_k) \big)\nonumber \\
%&\leq& -2\mu \big(f(x_k) - f(x_*)\big) -2\mu \big( \langle \nabla f(x_k), y_k - x_k \rangle + \frac{\mu}{2}\|y_k - x_k\|^2\nonumber\\
%&=&-2\mu \big(f(x_k) - f(x_*)\big) -2\mu \big( \langle \nabla f(x_k)-\nabla f(y_k), y_k - x_k \rangle \\
%&& \quad + \langle \nabla f(y_k), y_k - x_k \rangle + \frac{\mu}{2}\|y_k - x_k\|^2
%\big) \nonumber\\
%&\leq& -2\mu \big(f(x_k) - f(x_*)\big) + 2\mu L \|y_k - x_k\|^2 
%-2\mu   \langle \nabla f(y_k), y_k - x_k \rangle \\
%&& \quad - \mu^2 \|y_k - x_k\|^2 \label{ineq-PL}
%\big)
%\eeq
%where in the first inequality we used PL inequality and in the second inequality we used \eqref{ineq-strcvx-smooth} and in the last inequality we used $L$-smoothness of $f$. Since $y_k - x_k = \nu(x_k - x_{k-1})$, \eqref{ineq-PL} implies
%$$ \begin{bmatrix} 
% 	{\xi}_k^c \\
% 	\nabla f({y}_k) \\
% 	{w}_k
% \end{bmatrix}^T
% \left(
% c_3  \tilde{M}_3 \otimes I_d \right) 
% \begin{bmatrix} 
% 	{\xi}_k^c \\
% 	\nabla f({y}_k) \\
% 	{w}_k
% \end{bmatrix} + 2\mu c_3 c_1  V_{P, c_1}(\xi_k) \leq 0
%$$ 
%}
Combining the last two inequalities, we conclude. %\myqed
\endproof
%\end{proof}
%\newpage
\section{Online Companion.}

\subsection{Proof of Proposition \ref{prop-robust-hb}.}
%\begin{proof}
 It is straightforward to check  that the derivative of $a(\kappa$) with respect to $\kappa$ for $\kappa \geq 32$ satisfies 
 $$ a'(\kappa) = \frac{ \frac{4\kappa-1}{\sqrt{2\kappa-1}} -1 } { 2\sqrt{\kappa}(\kappa-1)} - \frac{ \sqrt{2\kappa^2 - \kappa} - \sqrt{\kappa }}{(\kappa - 1)^2 } > 0.$$
Therefore, $a(\kappa)$ is an increasing function of $\kappa$ for $\kappa \geq 32$ and  $a(\kappa)> a(32)>5/4$ for $\kappa>32$. Also, it can be seen that
$ a(\kappa) = \sqrt{2} - \frac{1}{\sqrt{\kappa}} + \mathcal{O}(\frac{1}{\kappa})$
 as $\kappa \to \infty$. In particular,  we have the limit $a(\kappa) \to \sqrt{2}$ as $\kappa \to \infty$ and $1\leq a(\kappa)\leq \sqrt{2}$ for every $\kappa>1$.

For the HB method we have $\nu=0$, and it is well-known that the eigenvalues of the $2d\times 2d$ iteration matrix $A_Q$ matrix are given 
%On the other hand, it follows from the proof of Lemma \ref{thm-tmm-risk-formula} that the matrix $A_Q$ is similar to the block diagonal matrix with blocks $\tilde{A}^{(i)}$ for $i=1,2,\dots,d$. Therefore, the eigenvalues of the matrices $\tilde{A}^{(i)}$ coincide with the eigenvalues of $A_Q$. For the heavy-ball method, $\nu = 0$ in which case eigenvalues of $\tilde{A}^{(i)}$  are given 
by all the roots of the following quadratic equations:
$$ x^2 - (1+\beta - \alpha \lambda_i) x  +\beta = 0 \quad \mbox{for} \quad i=1,2,\dots,d.
$$ 
(see e.g. \cite{can2019accelerated,wright2022optimization}). In particular, for 
\beq 
\beta \geq \max_{i=1,2,\dots,d} (1-\sqrt{\alpha\lambda_i})^2 = \max((1-\sqrt{\alpha\mu})^2, (1-\sqrt{\alpha L})^2),
\label{ineq-beta-lower-bound}
\eeq
all the roots are complex numbers, each with magnitude $\sqrt{\beta}$. In this case, the spectral radius of the iteration matrix $A_Q$ is $\rho = \sqrt{|\beta|}$. The rest of the proof will follow by verifying that our choice of parameters lies in the set $\mathcal{S}_1 \cap \mathcal{S}_2$ that was characterized in Theorem \ref{thm-h-inf}. First, by our assumption on the stepsize, it can be checked that  $\alpha \leq \frac{4}{(\sqrt{L}+\sqrt{\mu})^2}$ and \eqref{ineq-beta-lower-bound} holds with $\beta = (1-\sqrt{\alpha\mu})^2$. Furthermore, it is also straightforward to check that $\alpha \leq  \frac{1+\beta}{L} $.
Recall from Theorem \ref{thm-h-inf} that
$$r_\lambda = \begin{cases} | 1-c_\lambda | \sqrt{1 - \frac{b_\lambda^2}{4c_\lambda}} & \mbox{if   } c_\lambda>0 \mbox{ and } \frac{|b_\lambda| (1+c_\lambda)}{4c_\lambda} < 1,\nonumber \\
	\big| |1+c_\lambda| - |b_\lambda| \big| & \mbox{otherwise}.
	\end{cases} $$ 
Note that for the heavy-ball method with general parameters $(\alpha,\beta)$, we have $\nu=0$, $b_{\lambda} =\alpha\lambda -(1+\beta), \quad c_\lambda = \beta$ where $b_{\lambda}$ and $c_{\lambda}$ are as in Theorem \ref{thm-h-inf}. Since  $\alpha \leq 1/\mu$ it can also be shown that we have 
$
\frac{|b_\mu| (1+ c_\mu)}{4 c_\mu} = (1+\beta - \alpha \mu) \frac{(1+\beta)}{4\beta} 
$
so that 
$r_\mu = | |1 + c_\mu| - |b_\mu|| = \alpha \mu$. Similarly, we compute 
$$r_L = \begin{cases} | 1-\beta | \sqrt{1 - \frac{(1+\beta - \alpha L)^2}{4\beta}} & \mbox{if   }  \frac{(1+\beta - \alpha L) (1+\beta)}{4\beta} < 1,\nonumber \\
	\alpha L  & \mbox{otherwise}.
	\end{cases} 
$$  
Therefore, from Theorem \ref{thm-h-inf}, we see that $L_{2,*}=H_\infty = 1/\sqrt{2\mu}$ provided that
$\frac{\alpha}{\sqrt{2}} \frac{\sqrt{L}}{r_L} \leq  \frac{\alpha}{\sqrt{2}} \frac{\sqrt{\mu}}{r_\mu}= \frac{1}{ \sqrt{2\mu}}.$
This inequality is satisfied if
$$\alpha \leq  \frac{1}{\sqrt{L\mu}} | 1-\beta | \sqrt{1 - \frac{(1+\beta - \alpha L)^2}{4\beta}}
\quad \mbox{whenever} \quad \frac{(1+\beta - \alpha L) (1+\beta)}{4\beta} < 1,$$
where we recall that $\beta = (1 - \sqrt{\alpha \mu})^2 =  (1- \frac{a(\kappa)}{\kappa})^2$.
% Assume $\alpha = a^2/L$ with $1\leq a^2 \leq 1+\beta = 1 + (1-\frac{a}{\sqrt{\kappa}})^2$. This holds when
%\beq 1 \leq a \leq A(\kappa)= \frac{\sqrt{\frac{4}{\kappa} + 8(1- \frac{1}{\kappa}) } - \frac{2}{\sqrt{\kappa}} }{2(1- \frac{1}{\kappa})} = \frac{\sqrt{\kappa} (\sqrt{2\kappa -1 } -1 ) }{\kappa - 1}
%\label{ineq-step-bound-a}
%\eeq
%Note also that throughout the paper we assume  $L>\mu$, therefore $\kappa>1$ and the denominator of $A(\kappa)$ cannot vanish. 
Notice that, with our stepsize choice, we have $\alpha \geq \frac{1}{L} \frac{(1-\beta)^2}{1+\beta}$ so that the inequality
$\frac{(1+\beta - \alpha L) (1+\beta)}{4\beta} < 1$ holds. Therefore, if we can show that
\beq \alpha \leq  \frac{1}{\sqrt{L\mu}} | 1-\beta | \sqrt{1 - \frac{(1+\beta - \alpha L)^2}{4\beta}},
\eeq
this will imply $H_\infty = 1/\sqrt{2\mu}$. Plugging in $\beta = (1- \frac{a(\kappa)}{\kappa})^2$, this is equivalent to

\beq a^2(\kappa) \leq \sqrt{\kappa} \left( 1 - \left(1 - \frac{a(\kappa)}{\sqrt{\kappa}}\right)^2\right) 
\sqrt{ 1 -  \frac{ \left(1 +(1-\frac{a(\kappa)}{\sqrt{\kappa}})^2 - a^2(\kappa)\right )^2 }{4 \left(1-\frac{a(\kappa)}{\sqrt{\kappa}}\right)^2} 
}.
\label{ineq-to-prove-hb}
\eeq 
 Using the identities,
$$ \left( 1 - (1 - \frac{a(\kappa)}{\sqrt{\kappa}})^2\right) 
= \frac{a(\kappa)}{\sqrt{\kappa}} \left(2 - \frac{a(\kappa)}{\sqrt{\kappa}}\right),
 $$ 
\beqs 
\small
1 -  \tfrac{ \left(1 +(1-\frac{a(\kappa)}{\sqrt{\kappa}})^2 - a^2(\kappa)\right )^2 }{4 (1-\tfrac{a}{\sqrt{\kappa}})^2} &=& 
\left( 1 -  \tfrac{ \left(1 +(1-\tfrac{a(\kappa)}{\sqrt{\kappa}})^2 - a^2(\kappa)\right ) }{2 (1-\tfrac{a(\kappa)}{\sqrt{\kappa}})} 
\right)
\left( 1 + \tfrac{ \left(1 +(1-\tfrac{a(\kappa)   }{\sqrt{\kappa}})^2 - a^2(\kappa)\right ) }{2 (1-\tfrac{a(\kappa)}{\sqrt{\kappa}})} 
\right)\\
&=&
\frac{  a^2(\kappa)  (1- \frac{1}{\kappa})
}{2 (1-\frac{a}{\sqrt{\kappa}})} 
\frac{  4 - a^2(\kappa)  - \frac{4a(\kappa) }{\sqrt{\kappa}}    + \frac{a(\kappa) }{\kappa}
}{2 (1-\frac{a}{\sqrt{\kappa}})},
\eeqs
the inequality \eqref{ineq-to-prove-hb} becomes
\beq 
a^2(\kappa)\leq \left(2 - \frac{a(\kappa)}{\sqrt{\kappa}}\right)
\frac{  a^2(\kappa)  \sqrt{1- \frac{1}{\kappa}}
}{2 (1-\frac{a(\kappa)}{\sqrt{\kappa}})} 
  \sqrt{4 - a^2(\kappa) - \frac{4a(\kappa)}{\sqrt{\kappa}}    + \frac{a(\kappa)}{\kappa}, 
}
\label{ineq-inter-hb}
\eeq 
which holds if and only if 
\beq
2 (1-\frac{a(\kappa)}{\sqrt{\kappa}}) \leq (2 - \frac{a(\kappa)}{\sqrt{\kappa}})
  \sqrt{1- \frac{1}{\kappa}}
  \sqrt{4 - a^2(\kappa) - \frac{4a(\kappa)}{\sqrt{\kappa}}    + \frac{a(\kappa)}{\kappa}.
}
\label{ineq-hb-to-prove-robust-hb}
\eeq
For $\kappa < 32$, $a(\kappa)=1$ and it is straightforward to check that the inequality \eqref{ineq-hb-to-prove-robust-hb} is satisfied. Otherwise, for $\kappa\geq 32$, the function  $g(a):=\sqrt{4 - a^2 - \frac{4a}{\sqrt{\kappa}}    + \frac{a}{\kappa} }$ is decreasing on the interval $1\leq a \leq \sqrt{2}$ with
$ g(a)\geq g(\sqrt{2}) = \sqrt{2 - \frac{4\sqrt{2}}{\sqrt{\kappa}} +  \frac{\sqrt{2}}{\kappa }} >0. 
$
Therefore, it suffices to show that for $\kappa \geq 32$, we have
\beq
2 (1-\frac{a(\kappa)}{\sqrt{\kappa}}) \leq (2 - \frac{a(\kappa)}{\sqrt{\kappa}})
  \sqrt{(1- \frac{1}{\kappa})}
 g(\sqrt{2}).
 \label{ineq-robust-hb-to-show}
\eeq
Noticing that we have $ \sqrt{1- \frac{1}{\kappa}}
 g(\sqrt{2}) \geq 1$ for $\kappa\geq 32$, the inequality \eqref{ineq-robust-hb-to-show} trivially holds. We conclude that $L_{2,*}=H_\infty = \frac{1}{\sqrt{2\mu}}$.  Furthermore, the convergence rate satisfies $ \rho = \sqrt{|\beta|} = 1 - \sqrt{\alpha \mu} = 1 - \frac{a(\kappa)}{\sqrt{\kappa}}$ as claimed. This completes the proof.
 %For $1\leq \kappa <32$, we can check that $a=1$ satisfies \eqref{ineq-to-prove-hb}.
%\end{proof} 
%\begin{proof} Assume that $f\in \Cml$ is a quadratic function and that the multiplicity of the eigenvalues $\mu$ and $L$ of the Hessian of $f$ is one. Assume also that $c_{\lambda_*}>0$ and $\frac{|b_{\lambda_*}| (1 + 	c_{\lambda_*})}{4c_{\lambda_*}} < 1$ where  $\lambda_* \in \{\mu, L\}$ is the maximizer of the optimization problem in \eqref{eq-hinf-quad} and that $\frac{\sqrt{\mu}}{r_\mu} \neq \frac{\sqrt{L}}{r_L}$. Then, $L_{2,*} < H_\infty$.
%\end{proof} 
\subsection{Recovering explicit bounds for GD and NAG with the MI approach.}\label{app-suppmat-MI-implies-explicit-bounds}
%In the following discussion, we argue that for GD and NAG, we can choose the parameters of the matrix inequality $\tilde{M}_2 +c_1 \tilde{M_1} + c_0 \tilde{M}_0 \succeq 0$ 
In the following discussion, we argue that the matrix inequality GMM approach is able to recover the explicit bounds obtained for GD and NAG in Prop. \ref{prop-hinf-gd-bound} and Prop. \ref{prop-hinf-agd-bound} in the sense that if we choose the parameters $c_0, c_1, a,b, \rho_0, \rho_1,\rho_2,\rho_3$ from Theorem \ref{thm-hinfty-bound} in a particular fashion (that we will make precise), the matrix inequality 
$\tilde{M}_2 +c_1 \tilde{M_1} + c_0 \tilde{M}_0 \succeq 0$ given in Theorem \ref{thm-hinfty-bound} is feasible and for these parameters, Theorem \ref{thm-hinfty-bound} yields bounds for the $\ell_2$ gain that match those obtained in Prop. \ref{prop-hinf-gd-bound} and Prop. \ref{prop-hinf-agd-bound} for GD and NAG. \looseness=-1%This discussion, 
%This discussion justifies why the constants $c_0, c_1, a,b, \rho_0, \rho_1,\rho_2,\rho_3$ are needed in Theorem \ref{thm-hinfty-bound} 
%In the following discussion, we argue that for GD and NAG, we can choose the parameters from Theorem \ref{thm-hinfty-bound} explicitly so that the matrix inequality $\tilde{M}_2 +c_1 \tilde{M_1} + c_0 \tilde{M}_0 \succeq 0$ is satisfied. Our discussion also shows that the matrix inequality GMM approach is able to recover the explicit bounds obtained for GD and NAG in Prop. \ref{prop-hinf-gd-bound} and Prop. \ref{prop-hinf-agd-bound} in the sense that for particular choice of parameters for the matrix inequality given in Prop. \ref{prop-hinf-gd-bound}, we can show that the matrix inequality is feasible and  Prop. \ref{prop-hinf-gd-bound} yields the same bounds as those obtained in Prop. \ref{prop-hinf-gd-bound} and Prop. \ref{prop-hinf-agd-bound}.
% for particular choice of parameters that we will precise for the matrix inequality  given in Prop. \ref{prop-hinf-gd-bound}.
%\begin{rema} 
\subsubsection{GD case.}
GD is the special case of the GMM method when $\beta=\nu = 0$. We will argue that Theorem \ref{thm-hinfty-bound} can recover the result obtained in Theorem \ref{thm-h-inf} in this special case. For GD, we take $\tilde{P}_{22}=\tilde{P}_{12}=0$, $\tilde{P}_{11} = 1$, $c_1 = 0$  in which case the Lyapunov function becomes $V_{P,c_1}({\xi}_k) = \|{x}_k - x_*\|^2$. Furthermore, the second rows and the second columns of the matrices in the matrix inequality $\tilde{M}_2 +c_1 \tilde{M_1} + c_0 \tilde{M}_0 \succeq 0$  become zero. Therefore, omitting the second rows and columns, and taking $b=\rho_1 = \rho_2 = \rho_3 = 0$ in Theorem \ref{thm-hinfty-bound},  this matrix inequality is equivalent to the $3\times 3$ matrix inequality
\beq
\begin{bmatrix} 
(\rho_0^2 - 1)   & \alpha & \alpha \\
\alpha       &  -\alpha^2 & -\alpha^2 \\
\alpha       &    -\alpha^2 & a
\end{bmatrix} 
+ c_0 
\begin{bmatrix} 
2mL & -(m+L) & 0 \\
-(m+L) & 2 & 0 \\
0 & 0 & 0      
\end{bmatrix}
\succeq 0.
\label{ineq-matrix-3by3-gd}
\eeq
From the standard analysis of gradient descent methods, it is known that there exists a positive constant $c_3$ such that 
\beq
\begin{bmatrix} 
(\rho_{GD}^2 - 1)   & \alpha  \\
\alpha       &  -\alpha^2 
\end{bmatrix} 
+ c_3 
\begin{bmatrix} 
2mL & -(m+L)  \\
-(m+L) & 2     
\end{bmatrix}
\succeq 0,
\label{ineq-standard-gd-lmi}
\eeq
(see  \cite[Section 4.4]{lessard2016analysis}) which, after reorganizing, is equivalent to 

\beq
\begin{bmatrix} 
-1  & \alpha  \\
\alpha       &  -\alpha^2 
\end{bmatrix} 
\succeq 
\begin{bmatrix}
-\rho_{GD}^2 & 0\\
0 & 0 
\end{bmatrix} 
- c_3 
\begin{bmatrix} 
2mL & -(m+L)  \\
-(m+L) & 2     
\end{bmatrix}.
\label{ineq-matrix-reorganized}
\eeq
It can also be checked that for any positive constant $c_4$, 
\beqs
\begin{bmatrix} 
0 & 0 & \alpha \\
0 & 0 &  -\alpha^2 \\
\alpha & -\alpha^2 & 0 
\end{bmatrix} 
&\succeq& 
\begin{bmatrix}
-c_4 & \alpha c_4 & 0 \\
\alpha c_4 & -c_4 \alpha^2 & 0 \\
0 & 0 & -\alpha^2 / c_4
\end{bmatrix}\\
&\succeq&
\begin{bmatrix}
-c_4 \rho_{GD^2}& 0  & 0 \\
0 & 0 & 0 \\
0 & 0 & -\alpha^2 / c_4
\end{bmatrix}  - c_3 c_4 \begin{bmatrix} 
2mL & -(m+L)  & 0\\
-(m+L) & 2     & 0 \\
0 & 0 & 0
\end{bmatrix},
\eeqs
where we used \eqref{ineq-matrix-reorganized} in the latter matrix inequality. Based on this inequality, a sufficient condition for \eqref{ineq-matrix-3by3-gd} to hold is 
the existence of $c_0$ such that
\begin{small}
\beqs
\begin{bmatrix} 
(\rho_0^2 - c_4 \rho_{GD}^2)- 1   & \alpha & 0 \\
\alpha       &  -\alpha^2 & 0 \\
0      &  0 & a - \alpha^2/c_4
\end{bmatrix} 
+ (c_0 - c_3 c_4)
\begin{bmatrix} 
2mL & -(m+L) & 0 \\
-(m+L) & 2 & 0 \\
0 & 0 & 0      
\end{bmatrix}
\succeq 0.
\eeqs
\end{small}
With the choice of $c_4 = \frac{1}{\rhogd}- 1$, $a= \alpha^2/c_4$, $c_0 = c_3 (1+c_4)$ this is equivalent to 
\beqs
\begin{bmatrix} 
(\rho_0^2 - c_4 \rho_{GD}^2)- 1   & \alpha  \\
\alpha       &  -\alpha^2 
\end{bmatrix} 
+ c_3
\begin{bmatrix} 
2mL & -(m+L)  \\
-(m+L) & 2 
\end{bmatrix}
\succeq 0.
\eeqs
From \eqref{ineq-standard-gd-lmi}, we see that this inequality holds if 
$\rho_0^2 - c_4 \rhogd^2 = \rhogd^2$ or equivalently for 
$$ \rho_0^2 = (1+c_4) \rhogd^2 = \rhogd.$$
Plugging this identity together with $c_1 = 0, a= \alpha^2/c_4, r(\tilde P) = 1$ leads to $L_{2,*} \leq \left( \frac{L}{2} \frac{\alpha^2}{(1-\rhogd)^2} \right)^{1/2}$ which is equivalent to the inequality \eqref{ineq-gd-hinfty-bd1} obtained in Proposition \ref{prop-hinf-gd-bound}. 
%\mg{On the other hand, if we take $c_1 = 1, c_0 = 0, b= 0, a=\frac{\alpha( 1- L\alpha)}{2r}$ for some $r>0$ that will be specified later, $P=0$, $c_3 = 0 $}, then for GD where $\nu=\beta=0$, \eqref{constraint-sdp} is equivalent to
%Since $\begin{bmatrix} 0 &   \frac{1}{2}\alpha (1-\alpha L)  \\
% \frac{1}{2}\alpha (1-\alpha L)  & 0 
%\end{bmatrix}  \succeq \begin{bmatrix}  -\frac{1}{2}\alpha |1-\alpha L|r  & 0 \\
% 0 & -\frac{1}{2}\alpha |1-\alpha L|/r,   
%\end{bmatrix}$

\mg{On the other hand, given GD stepsize $\alpha \in (0,2/L)$, assume the inequalities $1 > \rho_0^2 \geq 1-2\mu \alpha (1 - \frac{L\alpha}{2}) + \alpha \mu |1-\alpha L| r \geq 0 $ hold for some $\rho_0\geq 0 $ and $r>0$ that we will specify below. Then, we have 
\beq \begin{bmatrix}  \frac{\mu}{2}(1-\rho_0^2) & 0 & -\frac{(1-\rho_0^2) }{2}& 0\\ 
0   &  0  & 0 & 0 \\
 -\frac{(1-\rho_0^2) }{2}& 0 &  \alpha (1- \frac{\alpha L}{2}) -\frac{1}{2}\alpha |1-\alpha L|r   &  0 \\
 0 & 0 & 0 &  0
\end{bmatrix} \succeq 0,
\eeq
which, based on the matrix inequality $\begin{bmatrix}  \frac{1}{2}\alpha |1-\alpha L|r   &   \frac{1}{2}\alpha (1-\alpha L)  \\
 \frac{1}{2}\alpha (1-\alpha L)  & \frac{1}{2}\alpha |1-\alpha L|/r  
\end{bmatrix}  \succeq 0 $, 
implies
\beq \begin{bmatrix}  \frac{\mu}{2}(1-\rho_0^2) & 0 & -\frac{(1-\rho_0^2) }{2}& 0\\ 
0   &  0  & 0 & 0 \\
 -\frac{(1-\rho_0^2) }{2}& 0 &  \alpha (1- \frac{\alpha L}{2})  &  \frac{1}{2}\alpha (1-\alpha L)  \\
 0 & 0 & \frac{1}{2}\alpha (1-\alpha L) &  \frac{1}{2}\alpha |1-\alpha L|/r
\end{bmatrix} \succeq 0.
\label{ineq-sdp-gd}
\eeq
}
\mg{
This matrix inequality is equivalent to $\tilde{M}_2 +c_1 \tilde{M_1} + c_0 \tilde{M}_0 \succeq 0$ with $c_1 = 1, c_0 = 0, b= 0, a=\frac{\alpha(|1- L\alpha|)}{2r}$, $P=0$, $\rho_1 = \rho_2 = \rho_3 = 0$. Choosing $\rho_0^2 = 1-2\mu \alpha (1 - \frac{L\alpha}{2}) + \alpha \mu |1-\alpha L| r$ with $r=1$ when $\alpha \leq \frac{1}{L}$ and $r = \frac{2-\alpha L}{\alpha L}$ for $\alpha \in (1/L,2/L)$, Theorem \ref{thm-hinfty-bound} implies the other $\ell_2$ gain bounds obtained in Prop. \ref{prop-hinf-agd-bound}. Hence, we conclude that particular choices of parameters in Theorem \ref{thm-hinfty-bound} recovers the same $L_{2,*}$ bound previously obtained in Proposition \ref{prop-hinf-gd-bound} for GD.
}
%H_\infty \leq \frac{1}{1-\rhogd} \frac{(1+a)}{} 
%\end{rema}
%\begin{rema} 
\subsubsection{NAG case.}
We will show that Theorem \ref{thm-hinfty-bound} can recover the $L_{2,*}$ bound we obtained in Proposition \ref{prop-hinf-agd-bound} for the NAG method with $\beta = \nu =  \frac{1-\sqrt{\alpha\mu}}{1+\sqrt{\alpha\mu}}$ for $\alpha \in (0,1/L]$. 
More specifically, in Theorem \ref{thm-hinfty-bound}, we choose 
$$\rho_1^2 = \frac{2\alpha}{4s_1(1+\sqrt{\alpha\mu})^2}, \quad \rho_2^2 = \frac{2}{\mu s_1}
  \frac{\alpha^2\mu^2 + 2\alpha \mu + \alpha\mu (1-\sqrt{\alpha \mu})}{4(1+\sqrt{\alpha\mu})^2}, \quad \rho_3^2 = \frac{2}{\mu s_1} \frac{\alpha\mu (1-\sqrt{\alpha \mu})}{4(1+\sqrt{\alpha\mu})^2},$$ 
$c_1 = 1, c_0 = 0, \rho_0^2 = \rho_{NAG}^2 = 1-\sqrt{\alpha\mu}$, 
$a=  s_1 + \frac{L\alpha^2}{2} s_2$ and $b = \frac{L\alpha^2}{2s_2}$ with $s_i = \sqrt{\alpha}\hat{s}_i$ for $i=1,2$ where $\hat{s}_i$ is as in \eqref{def-hat-r1-r2}. We also take the $\tilde{P}$ matrix  according to \eqref{def-lyap-agd}. It suffices to show that the matrix inequality $\tilde{M}_2 +c_1 \tilde{M_1} + c_0 \tilde{M}_0 \succeq 0$ holds, because then for these parameters Theorem \ref{thm-hinfty-bound}
implies directly the bound \eqref{ineq-hinf-bound-agd-by-hand}. To show that the matrix inequality $\tilde{M}_2 +c_1 \tilde{M_1} + c_0 \tilde{M}_0 \succeq 0$ is satisfied, we first observe that by Cauchy-Schwarz, 
 \beq 2w_k^T B^T P A \xi_k^c  \leq \frac{1}{s_1} \|B^T PA  \xi_k^c \|^2  + s_1 \|w_k\|^2 = \frac{I_k(\alpha) }{s_1 4(1+\sqrt{\alpha\mu})^2}  + s_1 \|w_k\|^2, 
 \eeq
for any vectors $w_k$ and $\xi_k^c$. Using \eqref{ineq-I-k-alpha}, we obtain %we obtain the inequality %from the proof of Proposition \ref{prop-hinf-agd-bound}. T
%These bounds are equivalent to the matrix inequality of the form  
  $$ z_k^T
\left( 
  \left[ \begin{array}{@{}c|c@{}}
  \huge{0_{3\times 3}}
      & \begin{matrix}   -\tilde{A}^T \tilde{P} \tilde{B} \\ 0 \end{matrix} \\
   \cmidrule[0.4pt]{1-2}
   \begin{matrix} - \tilde{B}^T \tilde{P} \tilde{A}  &0 \end{matrix}  & 0\\
\end{array} \right]    
+ 
 \begin{bmatrix} \rho_1^2 \tilde{P}_{11} +\frac{\mu}{2} \rho_2^2 &  \rho_1^2 \tilde{P}_{12} & 0 & 0 \\
 							   				 \rho_1^2 \tilde{P}_{12} & \rho_1^2 \tilde{P}_{22} +\frac{\mu}{2}\rho_3^2 & 0 & 0\\
 							   				 					0   & 0						& 0 & 0 \\
 							   				 					0 & 0 & 0 & s_1 
 \end{bmatrix} \right)  z_k \succeq 0,
$$
for any $ z_k = \begin{bmatrix} 
 	{\xi}_k^c \\
 	\nabla f({y}_k) \\
 	{w}_k
 \end{bmatrix}$. Summing this inequality with the following inequality, 
$$ \frac{L\alpha^2}{2} \begin{bmatrix} 0 &  0 & 0 & 0 \\
 							   				 0 & 0 & 0 & 0\\
 							   				 					0   & 0						& \frac{1}{s_2} & -1  \\
 							   				 					0 & 0 & -1& s_2
 \end{bmatrix} \succeq 0, 
 $$
we obtain  
%$$ \begin{bmatrix}  \frac{L\alpha^2}{2s_2} & - \frac{L\alpha^2}{2} \\
%						-\frac{L\alpha^2}{2} & \frac{L\alpha^2}{2}s_2 
%\end{bmatrix} \succeq 0,
%$$ 
$$\tilde{M}_5:=  \left[ \begin{array}{@{}c|c@{}}
 \huge{0_{3\times 3}}
      & \begin{matrix}   -\tilde{A}^T \tilde{P} \tilde{B} \\ 0 \end{matrix} \\
   \cmidrule[0.4pt]{1-2}
   \begin{matrix} - \tilde{B}^T \tilde{P} \tilde{A}  &0 \end{matrix}  & 0\\
\end{array} \right]    
+ 
 \begin{bmatrix} \rho_1^2 \tilde{P}_{11} +\frac{\mu}{2} \rho_2^2 &  \rho_1^2 \tilde{P}_{12} & 0 & 0 \\
 							   				 \rho_1^2 \tilde{P}_{12} & \rho_1^2 \tilde{P}_{22} +\frac{\mu}{2}\rho_3^2 & 0 & 0\\
 							   				 					0   & 0						& b &  \frac{-L\alpha^2}{2}  \\
 							   				 					0 & 0 & \frac{-L\alpha^2}{2} & a 
 \end{bmatrix}\succeq 0.
 $$
Using $c_0=0$, $c_1=1$ and $\frac{-L\alpha^2}{2} = -\tilde{B}^T \tilde{P}\tilde{B} + \frac{\alpha(1-L\alpha)}{2}$, we obtain
\beq \tilde{M}_4 = \tilde{M}_2 + c_1 \tilde{M}_1 + c_0 \tilde{M}_0 = \tilde{M}_2 + \tilde{M}_1 = \tilde{M}_5 + \left[ \begin{array}{@{}c|c@{}}
 \huge{ \mathcal{S}_{\rho_0}(\tilde P)}
      & \begin{matrix}  0_{2\times 1} \\ 0 \end{matrix} \\
   \cmidrule[0.4pt]{1-2}
   \begin{matrix} 0_{1\times 2} &0 \end{matrix}  & 0\\
\end{array} \right]   \succeq 0,
\eeq
where the last inequality follows as $\tilde{M}_5 \succeq 0$ and we have $\mathcal{S}_{\rho_0}(\tilde{P})\succeq 0$ for our choice of $\beta = \nu =  \frac{1-\sqrt{\alpha\mu}}{1+\sqrt{\alpha\mu}}$ and $\alpha \in (0,1/L]$ (see the proof of Prop. \ref{prop-hinf-agd-bound}). This shows that the desired inequality $\tilde{M}_4=\tilde{M}_2 +c_1 \tilde{M_1} + c_0 \tilde{M}_0 \succeq 0$ holds.
%Combining everything together, we have
%\beq% \mg{c_3  \tilde{M}_3} + \tilde{M}_2(\tilde{P},a) +c_1 \tilde{M_1} + c_0 \tilde{M}_0 &=&  
%\quad \left( 
%  \left[ \begin{array}{@{}c|c@{}}
%  \huge{\mathcal{S}_{\rho_0}(\tilde{P})}
%      & \begin{matrix}   -\tilde{A}^T \tilde{P} \tilde{B} \\ 0 \end{matrix} \\
%   \cmidrule[0.4pt]{1-2}
%   \begin{matrix} - \tilde{B}^T \tilde{P} \tilde{A}  &0 \end{matrix}  & 0_{d\times d}\\
%\end{array} \right]    
%+ 
% \begin{bmatrix} \rho_1^2 \tilde{P}_{11} +\frac{\mu}{2} \rho_2^2 &  \rho_1^2 \tilde{P}_{12} & 0 & 0 \\
% 							   				 \rho_1^2 \tilde{P}_{12} & \rho_1^2 \tilde{P}_{22} +\frac{\mu}{2}\rho_3^2 & 0 & 0\\
% 							   				 					0   & 0						& 0 & 0 \\
% 							   				 					0 & 0 & 0 & s_1 
% \end{bmatrix} \right)   \succeq 0,
%\eeq
%where the left-hand side is equivalent to $c_3  \tilde{M}_3 + \tilde{M}_2 +c_1 \tilde{M_1} + c_0 \tilde{M}_0$ so 
Furthermore, by the identities \eqref{def-r1-r2} and \eqref{def-hat-r1-r2}, we have $\rho_1^2 + \rho_2^2 + \rho_3^2 = \frac{\sqrt{\alpha\mu}}{4}$ and $\frac{4b(\nu^2 + (1+\nu)^2) L^2}{\mu} c_1 = \frac{\sqrt{\alpha\mu}}{4}$. Then, with this choice of parameters, Theorem \ref{thm-hinfty-bound} implies the $L_{2,*}$ bound obtained in Proposition \ref{prop-hinf-agd-bound}.
%\end{rema}

% Acknowledgments here

% Enter the text of acknowledgments here

% References here (outcomment the appropriate case) 

% CASE 1: BiBTeX used to constantly update the references 
%   (while the paper is being written).
%\bibliographystyle{informs2014} % outcomment this and next line in Case 1
%\bibliography{<your bib file(s)>} % if more than one, comma separated

% CASE 2: BiBTeX used to generate mypaper.bbl (to be further fine tuned)
%\input{mypaper.bbl} % outcomment this line in Case 2
\subsection{Further illustrations for the spectral value sets of HB and TMM.}
%\mtodo{Add a picture about spectral sets here.}
\begin{figure}[h]
  \centering
    \includegraphics[width=0.47\linewidth, height=0.42\linewidth]{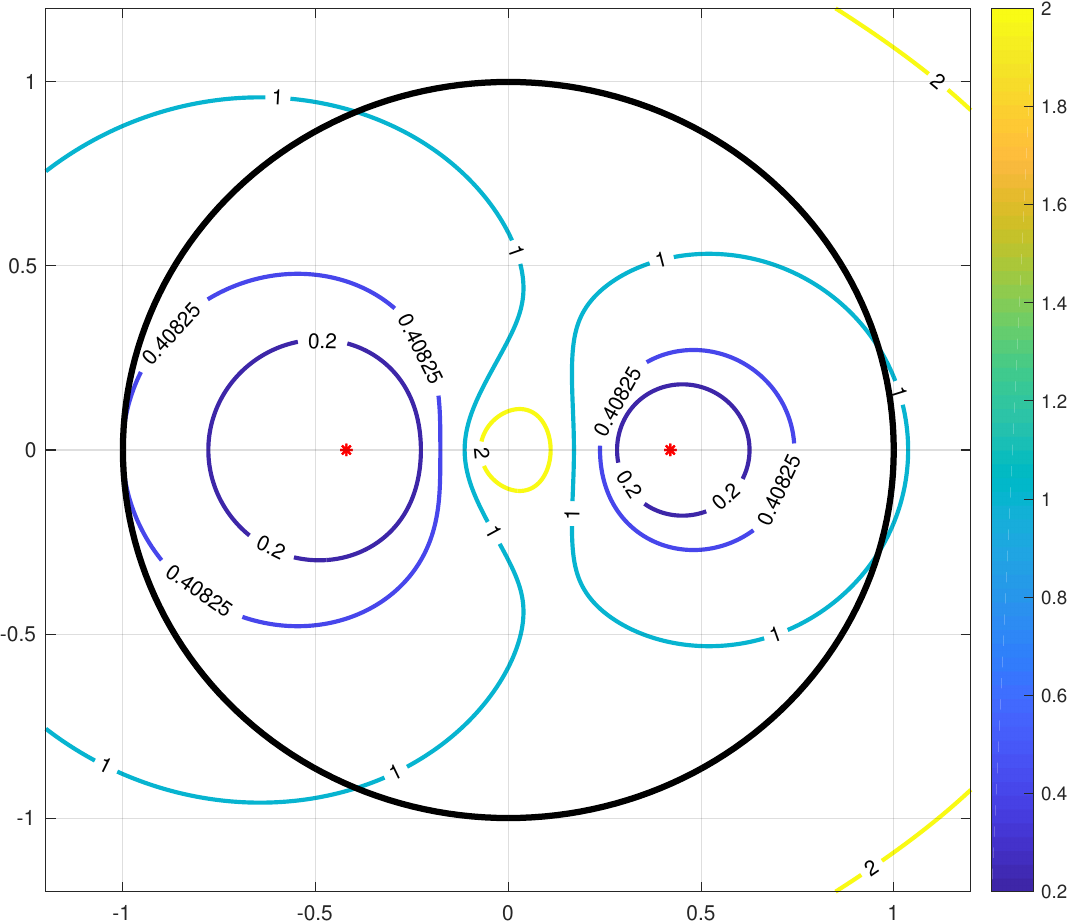}
      \includegraphics[width=0.47\linewidth, height=0.42\linewidth]{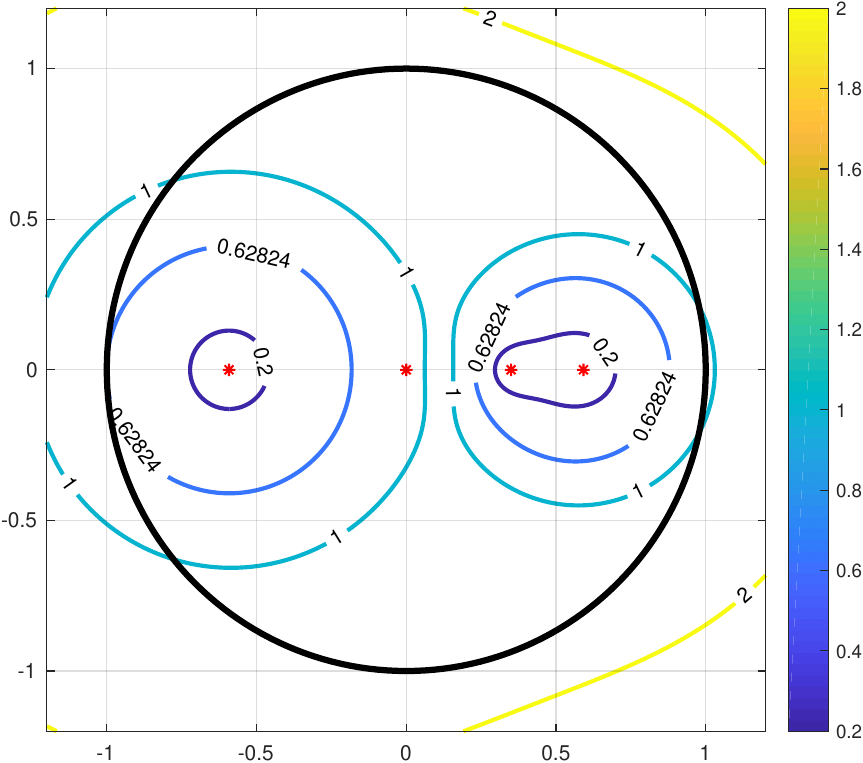}%{nag_spectral_sets_L3_muhalf_highbeta_v2-cropped.pdf}
        \caption[caption]{\label{fig-spec-val-set-hb-tmm} Boundary of the spectral value sets for $L=3$, $\mu=1/2$ and $\alpha = 1/L$. (Left panel) Heavy-ball with optimal (for rate) parameters  $\alpha=  \frac{4}{(\sqrt{L} + \sqrt{\mu})^2}$, $\beta = (\frac{\sqrt{\kappa}-1}{\sqrt{\kappa}+1})^2,\nu=0$.  (Right panel) TMM with parameters proposed in \cite{scoy-tmm-ieee}.}
\end{figure}
In Sec. \ref{subsec-multiplicative-noise}, we have provided illustrations of the spectral value sets of NAG; here in this Apsection, we provide additional similar illustrations for the HB and TMM methods. In Figure \ref{fig-spec-val-set-hb-tmm}, we plotted the spectral value sets for the HB method (left panel) and for TMM (right panel) based on the standard choice of parameters (see Table \ref{table}) for $L=3, \mu=\frac{1}{2}$ and $\kappa = 6$. We followed a similar approach to the one we used in Sec. \ref{subsec-multiplicative-noise} where we plot the boundary of the spectral value sets $\Lambda_\varepsilon$ as $\varepsilon$ is varied. We observe that for HB, $A_Q$ has two eigenvalues (each with multiplicity two) with the complex stability radius $\varepsilon_*  = \frac{1}{\sqrt{6}}\approx 0.40825$ which corresponds to $L_{2,*}=H_\infty = \frac{1}{\varepsilon_*}= \sqrt{6}$. This value also matches the formula $L_{2,*}=H_\infty = \frac{\sqrt{\kappa}}{\sqrt{2\mu}}$ given in Table \ref{table}. For TMM, the eigenvalues of $A_Q$ are real and simple; we observe that the complex stability radius $\varepsilon_* \approx 0.62284$ and $L_{2,*}=H_\infty = 2 - \frac{1}{\sqrt{6}}\approx \frac{1}{0.62284}$ as expected, based on the $H_\infty$ formula given in Table \ref{table} and the formula \eqref{def-complex-stab-rad}. We can conclude that TMM has better robustness (smaller $L_{2,*}$) compared to HB in this case. These plots illustrate further that choosing parameters to yield a smaller $\ell_2$ gain $L_{2,*}$ (or equivalently a smaller $H_\infty$ norm) enables better robustness to the deterministic relative noise $w_k$ satisfying \eqref{ineq-multiplicative-noise} and illustrate the fact that the multiplicative inverse of the $H_\infty$ norm is equal to the norm of the smallest perturbation matrix $\Delta$ such that the relative noise $w_k= \Delta z_k$ destabilizes the TMM iterations.

\end{document}